# KATO-MILNE COHOMOLOGY GROUP OVER RATIONAL FUNCTION FIELDS IN CHARACTERISTIC 2, I

AHMED LAGHRIBI AND TRISHA MAITI

## Contents

ABSTRACT. Let $\mathcal{F}$ be a field of characteristic 2. In this paper we determine the Kato-Milne cohomology of the rational function field $\mathcal{F}(x)$ in one variable $x$. This will be done by proving an analogue of the Milnor exact sequence [4] in the setting of Kato-Milne cohomology. As an application, we answer the open case of the norm theorem for Kato-Milne cohomology that concerns separable irreducible polynomials in many variables. This completes a result of Mukhija [17, Theorem A.3] that gives the norm theorem for inseparable polynomials.

## 1. Introduction

In 1970 Milnor gave a description of $W(F(x))$, the Witt group of symmetric bilinear forms of the rational function field $F(x)$ in one variable $x$, where $F$ is a field of characteristic not 2 [15]. More precisely, he proved a split exact sequence that expresses $W(F(x))$ in terms of the Witt groups of the residue fields given by all discrete valuations of $F(x)$. His result extends in the same way to symmetric bilinear forms in characteristic 2. Later, in 2006, Aravire and Jacob described the Witt group of nonsingular quadratic forms $W_q(F(x))$, where $F$ has characteristic 2 [4]. To this end, they developed an analogue of the Milnor exact sequence done in characteristic not 2 using some tricky adaptations to characteritic 2. Since the Witt group of nonsingular quadratic forms in characteristic 2 is related to the Kato-Milne cohomology by a celebrated isomorphism due to Kato [11], it is natural to ask for an analogue of the Milnor's split exact sequence in the case of Kato-Milne cohomology group of a rational function field in one variable. Our main objective in this paper is to get this analogue. We will limit ourselves to the case of Milnor's exact sequence. In a forthcoming paper, we will extend this sequence to an analogue of the Milnor-Scharlau exact sequence that uses transfer maps [19, Theorem 3.5, page 215]. In order to state our main result, we need to recall some definitions and results related to quadratic forms and Kato-Milne cohomology.

Let $F$ be a field of characteristic 2. Any $F$-quadratic form $\varphi$ decomposes as follows:

$$(1) \qquad \varphi \simeq [a_1, b_1] \perp \cdots \perp [a_r, b_r] \perp \langle c_1 \rangle \perp \cdots \perp \langle c_s \rangle$$

---

where $[a,b]$ denotes the binary quadratic form $aX^2 + XY + bY^2$, and $\langle c \rangle$ denotes the 1-dimensional quadratic form $cX^2$. The quadratic form $\varphi$ is called nonsingular when $s = 0$.

We denote by $W(F)$ (*resp.* $W_q(F)$) the Witt ring of non-degenerate symmetric bilinear forms over $F$ (*resp.* the Witt group of nonsingular quadratic forms over $F$). It is well known that $W_q(F)$ is endowed, in a natural way, with a $W(F)$-module structure induced by tensor product [8, Page 51]. Let $I(F)$ be the ideal of $W(F)$ of even dimensional bilinear forms over $F$. Clearly, this ideal is additively generated by $\langle\langle a \rangle\rangle_b := \langle 1, a \rangle_b$ for $a \in F^*$. For any integer $m \geq 1$, let $I^m(F)$ denote the $m$-th power of $I(F)$, this ideal is additively generated by the forms $\langle\langle a_1, \ldots, a_m \rangle\rangle_b := \langle 1, a_1 \rangle_b \otimes \ldots \otimes \langle 1, a_m \rangle_b$, where $a_i \in F^*$. Such a bilinear form is called an $m$-fold bilinear Pfister form (we take $I^0(F) = W(F)$). Similarly, for an integer $m \geq 0$, let $I_q^{m+1}(F)$ denote the $W(F)$-submodule of $W_q(F)$ given by $I^m(F) \otimes W_q(F)$. This submodule is additively generated by the quadratic forms $\langle\langle a_1, \ldots, a_m; b]] := \langle\langle a_1, \ldots, a_m \rangle\rangle_b \otimes [1, b]$, where $b \in F$ and $a_i \in F^*$ for $i = 1, \ldots, m$. These forms are called $(m+1)$-fold quadratic Pfister forms (we have $I_q^1(F) = W_q(F)$). We take $\overline{I}^m(F)$ (*resp.* $\overline{I}_q^{m+1}(F)$) the quotient $I^m(F)/I^{m+1}(F)$ (*resp.* the quotient $I_q^{m+1}(F)/I_q^{m+2}(F)$).

The space of 1-differential forms, denoted by $\Omega_F^1$, is the $F$-vector space generated by the symbols $\mathrm{d}\,a$ subject to the relations: $\mathrm{d}(a + b) = \mathrm{d}\,a + \mathrm{d}\,b$ and $\mathrm{d}(ab) = b\,\mathrm{d}\,a + a\,\mathrm{d}\,b$ for all $a, b \in F$. For any integer $m \geq 1$, let $\Omega_F^m = \wedge^m \Omega_F^1$ be the $m$-th exterior product of $\Omega_F^1$, that we call the space of $m$-differential forms (we take $\Omega_F^0 = F$ and $\Omega_F^m = 0$ for $m < 0$). Clearly, as $F$-vector space, $\Omega_F^m$ is generated by the symbols $\mathrm{d}\,a_1 \wedge \mathrm{d}\,a_2 \wedge \ldots \wedge \mathrm{d}\,a_m$, for $a_i \in F$.

Let $\mathcal{B} = \{e_i \mid i \in I\}$ be a 2-basis of $F$, where $I$ is an ordered set. Let us consider the set of maps $\Sigma_{m,F} = \{\sigma : \{1, \ldots, m\} \to I \mid \sigma(j) < \sigma(k) \text{ when } j < k\}$. It is well known that $\left\{\frac{\mathrm{d}e_{\sigma(1)}}{e_{\sigma(1)}} \wedge \frac{\mathrm{d}e_{\sigma(2)}}{e_{\sigma(2)}} \wedge \ldots \wedge \frac{\mathrm{d}e_{\sigma(m)}}{e_{\sigma(m)}} \mid \sigma \in \Sigma_{m,F}\right\}$ is an $F$-basis of $\Omega_F^m$. The Frobenius homomorphism is the map $\Phi : \Omega_F^m \to \Omega_F^m / \mathrm{d}\Omega_F^{m-1}$ defined by:

$$\Phi\left(\sum_\sigma c_\sigma \frac{\mathrm{d}e_{\sigma(1)}}{e_{\sigma(1)}} \wedge \ldots \wedge \frac{\mathrm{d}e_{\sigma(m)}}{e_{\sigma(m)}}\right) = \sum_\sigma c_\sigma^2 \frac{\mathrm{d}e_{\sigma(1)}}{e_{\sigma(1)}} \wedge \ldots \wedge \frac{\mathrm{d}e_{\sigma(m)}}{e_{\sigma(m)}} + \mathrm{d}\Omega_F^{m-1}, \text{ where } c_\sigma \in F.$$

The Artin-Schreier operator $\wp : \Omega_F^m \to \Omega_F^m / \mathrm{d}\Omega_F^{m-1}$ is the homomorphism given by $\wp = \Phi - \mathrm{Id}$. By a result of Cartier, one knows that $\wp$ does not depend on the 2-basis of $F$. The cokernel of $\wp$ is called the Kato-Milne cohomology group (of degree $m$) and denoted by $H_2^{m+1}(F)$. It is clear that we have $H_2^{m+1}(F) = \Omega_F^m / \left(\mathrm{d}\Omega_F^{m-1} + \wp(\Omega_F^m)\right)$.

For each positive integer $m$, we have the Kato's isomorphism $f_{m+1} : H_2^{m+1}(F) \to \overline{I}_q^{m+1}(F)$, defined by:

(2) $$f_{m+1}\left(\overline{x\frac{\mathrm{d}x_1}{x_1} \wedge \frac{\mathrm{d}x_2}{x_2} \wedge \ldots \wedge \frac{\mathrm{d}x_m}{x_m}}\right) = \overline{\langle\langle x_1, \ldots, x_m; x]]}.$$

Let $v_p : F \to \mathbb{Z}$ be a discrete valuation of $F$ and $p$ an uniformizer. Let $\overline{F_p}$ and $F_p$ be the residue field and the completion of $F$, respectively. We have the Teichmüller lifting[1] $\alpha : \overline{F_p} \hookrightarrow F_p$. Naturally, it induces a map $\alpha' : H_2^{m+1}(\overline{F_p}) \to H_2^{m+1}(F_p)$, defined on the generators as follows:

$$\alpha'\left(\overline{\overline{a}\frac{\mathrm{d}\overline{a_1}}{\overline{a_1}} \wedge \ldots \wedge \frac{\mathrm{d}\overline{a_m}}{\overline{a_m}}}\right) = \overline{\alpha(\overline{a})\frac{\mathrm{d}\alpha(\overline{a_1})}{\alpha(\overline{a_1})} \wedge \ldots \wedge \frac{\mathrm{d}\alpha(\overline{a_m})}{\alpha(\overline{a_m})}}, \quad \text{for all } \overline{a}, \overline{a_i} \in \overline{F}_p.$$

Now let us define
$$W_1(H_2^{m+1}(F_p)) := \mathrm{coker}\left(H_2^{m+1}(\overline{F_p}) \xrightarrow{\alpha'} H_2^{m+1}(F_p)\right).$$

In Corollary 2.6 we prove that this cokernel does not depend on the choice of Teichmüller lifting. Let $\beta_1 : H_2^{m+1}(F) \to H_2^{m+1}(F_p)$ be the map induced by the inclusion $F \hookrightarrow F_p$, and $\beta_2 : H_2^{m+1}(F_p) \to W_1(H_2^{m+1}(F_p))$

---
[1] We refer to the Appendix for the construction of this homomorphism.

induced by projection. Then, we define the (second) residue map $\partial_p : H_2^{m+1}(F) \to W_1(H_2^{m+1}(F_p))$ as the composition $\partial_p = \beta_2 \circ \beta_1$. As in [4], here also this residue map $\partial_p : H_2^{m+1}(F) \to W_1(H_2^{m+1}(F_p))$ plays an important role. Finally, our main result is the following:

**Theorem 1.1.** *Let $\mathcal{F}$ be a field of characteristic 2 and $F = \mathcal{F}(x)$ the rational function field in one variable $x$ over $\mathcal{F}$. Then, the following sequence is short exact:*

$$0 \to H_2^{m+1}(\mathcal{F}) \xrightarrow{i_3} H_2^{m+1}(F) \xrightarrow{\oplus_p \partial_p \oplus \partial'_{\frac{1}{x}}} \oplus_p W_1(H_2^{m+1}(F_p)) \oplus \left( W_1(H_2^{m+1}(F_{\frac{1}{x}}))/H_2^m(\mathcal{F}) \wedge \overline{\frac{dx}{x}} \right) \to 0,$$

*where $i_3$ is induced by the inclusion $\mathcal{F} \hookrightarrow F$, and the direct sum is taken over $\frac{1}{x}$ and all monic irreducible polynomials $p \in \mathcal{F}[x]$. Moreover, $\partial'_{\frac{1}{x}}$ is the composition of $\partial_{\frac{1}{x}}$ and the projection map.*

To prove this theorem, we will adapt to differential forms many arguments from the proof of Milnor exact sequence done by Aravire and Jacob in the case of Witt groups in characteristic 2 [4, Theorem 6.2]. Therefore, we are adopting several notations from their paper.

Now we summarize the contents of the paper. In the next section we compute the group $H_2^{m+1}(k)$, where $k$ is a field of characteristic 2 complete with respect to a discrete valuation $v : k \to \mathbb{Z}$. In Theorem 2.5, we prove the decomposition $H_2^{m+1}(k) = H_2^{m+1}(\overline{k}) \oplus \mathcal{R} \oplus H_2^m(\overline{k}) \wedge \overline{\frac{d\pi}{\pi}}$, where $\pi$ is an uniformizer, $\overline{k}$ is the residue field of $k$ and the group $\mathcal{R}$ will be defined later (see Definition 2.1). To prove this, we will have to fix a 2-basis of $k$ and a Teichmüller lifting. From this theorem, we deduce $W_1(H_2^{m+1}(k)) \simeq \mathcal{R} \oplus H_2^m(\overline{k}) \wedge \overline{\frac{d\pi}{\pi}}$. Consequently, we can compute the value of the residue map $\partial_\pi$.

In Section 3 we consider $F = \mathcal{F}(x)$, the rational function field in one variable $x$. Here we make a filtration $L_0 \subset L_1 \subset L_2 \subset \ldots$ of $H_2^{m+1}(F)$, and then in Section 4 we prove that

$$\bigoplus_{\deg p = d} \partial_p : L_d/L_{d-1} \to \bigoplus_{\deg p = d} W_1(H_2^{m+1}(F_p))$$

is an isomorphism. From this isomorphism, and using induction on $d$, we will get the short exact sequence

(3) $$0 \to L_0 \xrightarrow{i} H_2^{m+1}(F) \xrightarrow{\oplus_p \partial_p} \bigoplus_p W_1(H_2^{m+1}(F_p)) \to 0,$$

where $p$ varies over all monic irreducible polynomials. Next, using the definition of $L_0$, we get the following exact sequence

(4) $$0 \to H_2^{m+1}(\mathcal{F}) \xrightarrow{i_1} L_0 \xrightarrow{i_2} W_1(H_2^{m+1}(F_{1/x})) \xrightarrow{\eta} H_2^m(\mathcal{F}) \to 0.$$

Finally from the exact sequences (3) and (4), we get Theorem 1.1.

In Section 5 we apply Theorem 1.1 to prove the norm theorem for Kato-Milne cohomology. This norm theorem was done before by Mukhija only for normed irreducible polynomials which are inseparable using a different method based on transfer [17]. So in this paper we treat the remaining case that concerns irreducible separable polynomials. We finish this paper by giving an overview of the construction of the Teichmüller lifting which will play an important role.

## 2. The group $W_1(H_2^{m+1}(F))$

Let $F$ be an arbitrary field of characteristic 2 and $B := \{t_i \mid i \in \mathcal{I}\}$ a 2-basis of $F$, which means that the elements $\prod_{i \in \mathcal{I}} t_i^{\epsilon_i}$ such that $\epsilon_i \in \{0, 1\}$ and $\epsilon_i = 0$ for almost all $i \in \mathcal{I}$, is a basis of $F$ as $F^2$-vector space. We may suppose that $\mathcal{I}$ is endowed with an ordering.

Let $T = \oplus_{|\mathcal{I}|}\mathbb{Z}/2\mathbb{Z}$ be such that we consider it as an ordered set with lexicographic ordering with respect to $0 < 1$. Namely, for any $I = (I_i)_{i \in \mathcal{I}}$ and $J = (J_i)_{i \in \mathcal{I}} \in T$, we have $I < J$ if there exists an element $k \in \mathcal{I}$ such that $I_k < J_k$ and $I_i = J_i$ for all $i > k$ (i.e., we compare the components of $I$ and $J$ from right to left). We will use the notation $t^I$ to express $\prod_{i \in \mathcal{I}} t_i^{I_i}$, where $t_i \in B$ and $I = (I_i)_{i \in \mathcal{I}} \in T$. Therefore, each $f \in F$ can be uniquely written as
$$f = \sum_{I \in T} t^I f_I^2,$$
where $f_I \in F$ and $f_I = 0$ for almost all $I \in T$. Suppose that $F$ is complete with respect to a discrete valuation $v : F \to \mathbb{Z}$, and let $\pi$ be an uniformizer. Let $A_F$ and $\overline{F}$ be its valuation ring and residue field, respectively. Now consider a subset $\{t_i \mid i \in \mathcal{I}'\} \subset A_F$ such that $\{\overline{t_i} \mid i \in \mathcal{I}'\} \subset \overline{F}$ is a 2-basis of $\overline{F}$, where $\mathcal{I}'$ is an ordered set. Then, $C := \{t_i \mid i \in \mathcal{I}'\} \cup \{\pi\}$ will be a 2-basis of $F$. Set $\mathcal{I} = \mathcal{I}' \cup \{\infty\}$ and $t_\infty = \pi$. Now the ordering of $\mathcal{I}'$ can be extended to an ordering of $\mathcal{I}$ by taking $i < \infty$ for all $i \in \mathcal{I}'$. We take $T = \oplus_{|\mathcal{I}|}\mathbb{Z}/2\mathbb{Z}$. Therefore, in this case for $I \in T$, $t^I$ denotes the expression $\prod_{i \in \mathcal{I}'} t_i^{I_i} \pi^{I_\infty}$, where $I_i$ and $I_\infty$ are the components of $I$. Consequently, $v(t^I) = I_\infty$ which lies in $\{0, 1\}$. Now, from Theorem 6.6, we know that there will be a unique lifting monomorphism $\alpha : \overline{F} \hookrightarrow F$ with $\alpha(\overline{t_i}) = t_i$ for all $i \in \mathcal{I}'$. Therefore, we consider each $t_i$ as an element of $\overline{F}$ taking as a subfield of $F$.

Because of the completeness of $F$, we can assume that $F = \overline{F}((\pi))$. Let $R$ denotes the ring $\overline{F}[\pi^{-1}]$. Therefore, each element of $R$ has non-positive valuation. For $I = (I_i)_{i \in \mathcal{I}} \in T$, we define the sets $\tilde{I} = \{i \in \mathcal{I} \mid I_i \neq 0\}$ and $(T)_m = \{I \in T \mid |\tilde{I}| = m\}$, where $m \in \mathbb{N}$. One should note that for each $I \in T$, the set $\tilde{I}$ is finite.

**Definition 2.1.** *Let $F$ be a field as described above and $m \in \mathbb{N}$. We define $\mathcal{R}$ the subgroup of $H_2^{m+1}(F)$ generated by the forms:*
$$\sum_{I \in (T)_m} \sum_{J \in T, \, J+I > I} \overline{t^J r_{I,J}^2 \frac{dt_I}{t_I}} \quad \text{with} \quad r_{I,J} \in \pi^{-1}R \text{ and } r_{I,J} \neq 0 \text{ for at most finite number of } I \text{ and } J,$$
*where* $\dfrac{dt_I}{t_I} = \dfrac{dt_{i_1}}{t_{i_1}} \wedge \ldots \wedge \dfrac{dt_{i_m}}{t_{i_m}}$ *such that* $\tilde{I} = \{i_1, i_2, \ldots, i_m\}$ *and* $i_1 < i_2 < \ldots < i_m$.

**Remark 2.2.** *(1) If $a \in F$ and $v(a) > 0$, then by Hensel's lemma, $a \in \wp(F) := \{x^2 + x \mid x \in F\}$.*
*(2) If $v(a) < 0$, then in this case $a \notin \wp(F)$ if one of the following conditions holds:*
*(a) $v(a)$ is odd.*
*(b) The leading coefficient of $a$ with respect to $\pi$ is not a square in $\overline{F}$.*
*One can see that the the separable parts of the generators of $\mathcal{R}$ satisfy one of the conditions (a) and (b). So they are not in $\wp(F)$.*

**Lemma 2.3.** *Let $F$ be a field of characteristic 2, and $\{t_i \mid i \in \mathcal{I}\}$ an ordered 2-basis of $F$. Then, any element $\varphi \in H_2^{m+1}(F)$ can be written as follows*
$$\varphi = \sum_{I \in (T)_m} \overline{A_I \frac{dt_I}{t_I}},$$
*where* $A_I = \displaystyle\sum_{J \in T, \, J+I \geq I} t^J r_{I,J}^2$ *such that $r_{I,J} \in F$ and $r_{I,J} = 0$ for almost all $I$ and $J$.*

*Proof.* Since $\{t_i \mid i \in \mathcal{I}\}$ is a 2-basis of $F$, then $\varphi \in H_2^{m+1}(F)$ is a finite sum of elements $\overline{t^J b_{I,J}^2 \dfrac{dt_I}{t_I}}$, where $I \in (T)_m, J \in T$ and $b_{I,J} \in F^*$.

Let us consider a nonzero element $\overline{t^J b_{I,J}^2 \frac{dt_I}{t_I}}$ in the sum such that $I + J < I$. The condition $I + J < I$ implies: $\tilde{I} \cap \tilde{J}$ is nonempty, $Max\tilde{J} \leq Max\tilde{I}$ and $Max(\tilde{J} \cap \tilde{I}) > Max(\tilde{J} \setminus \tilde{I})$. Moreover, $\tilde{J} \nsubseteq \tilde{I}$, otherwise we would get $\overline{t^J b_{I,J}^2 \frac{dt_I}{t_I}} = 0$. Let $i_l = Max(\tilde{J} \cap \tilde{I})$. Now we compute:

$$\overline{t^J b_{I,J}^2 \frac{dt_I}{t_I}} = \overline{t^J b_{I,J}^2 \frac{dt_{i_1}}{t_{i_1}} \wedge \ldots \wedge \frac{d\,t_{i_{l-1}}}{t_{i_{l-1}}} \wedge \frac{d\,t_{i_l}}{t_{i_l}} \wedge \frac{d\,t_{i_{l+1}}}{t_{i_{l+1}}} \wedge \ldots \wedge \frac{dt_{i_m}}{t_{i_m}}} \quad \text{(where } \tilde{I} = \{i_1, \ldots, i_m\} \text{ with } i_1 < \ldots < i_m\text{)}$$

$$= \overline{t^J b^2{}_{I,J} \frac{dt_{i_1}}{t_{i_1}} \wedge \ldots \wedge \frac{d\,t_{i_{l-1}}}{t_{i_{l-1}}} \wedge \frac{dt^I}{t^I} \wedge \frac{d\,t_{i_{l+1}}}{t_{i_{l+1}}} \wedge \ldots \wedge \frac{dt_{i_m}}{t_{i_m}}} \quad \left(\text{using the relation } \frac{d(ab)}{ab} = \frac{d\,a}{a} + \frac{d\,b}{b}\right)$$

$$= \overline{t^J b_{I,J}^2 \frac{dt_{i_1}}{t_{i_1}} \wedge \ldots \wedge \frac{d\,t_{i_{l-1}}}{t_{i_{l-1}}} \wedge \frac{dt^{I+J}}{t^{I+J}} \wedge \frac{d\,t_{i_{l+1}}}{t_{i_{l+1}}} \wedge \ldots \wedge \frac{dt_{i_m}}{t_{i_m}}}$$

$$= \sum_{k \in \tilde{J} \setminus \tilde{I}} \overline{t^J b_{I,J}^2 \frac{dt_{i_1}}{t_{i_1}} \wedge \ldots \wedge \frac{d\,t_{i_{l-1}}}{t_{i_{l-1}}} \wedge \frac{dt_k}{t_k} \wedge \frac{d\,t_{i_{l+1}}}{t_{i_{l+1}}} \wedge \ldots \wedge \frac{dt_{i_m}}{t_{i_m}}}$$

$$= \sum_{I'} \overline{t^J b_{I,J}^2 \frac{dt_{I'}}{t_{I'}}} \quad \text{(where } \tilde{I}' = \{i_1, \ldots, i_{l-1}, i_{l+1}, \ldots, i_m\} \cup \{k\} \text{ such that } k < i_l\text{)}.$$

We prove that these $I'$ satisfy the relation $I' + J > I'$.

- If $i_l = i_m$, then $i_m \in \tilde{J} \setminus \tilde{I}'$ and $Max\,\tilde{I}' < i_m$. Therefore, $J + I' > I'$.
- If $i_l = i_{m-1}$, then $i_m \in \tilde{I}' \setminus \tilde{J}$, $i_{m-1} \in \tilde{J} \setminus \tilde{I}'$ and $k < i_{m-1}$ for all $k \in \tilde{J} \setminus \tilde{I}$. Therefore, $J + I' > I'$.
- $\vdots$
- If $i_l = i_1$, then $i_j \in \tilde{I}' \setminus \tilde{J}$ for $j = 2, \ldots m$, $i_1 \in \tilde{J} \setminus \tilde{I}'$, and $k < i_1$ for all $k \in \tilde{J} \setminus \tilde{I}$. Therefore, $J + I' > I'$.

Hence, $\overline{t^J b_{I,J}^2 \frac{dt_I}{t_I}}$ is in the desired form. This proves the lemma. $\square$

From this lemma, we get the following remark that will be used several times.

**Remark 2.4.** *Let $F$ be a field of characteristic $2$ and $\{t_i \mid i \in I\}$ an ordered $2$-basis of $F$. Let $A = \sum_{I \in (T)_m} \sum_{J \in T, J \neq 0} \overline{t^J r_{I,J}^2 \frac{dt_I}{t_I}}$, where $T, (T)_m$ are as defined previously. Then, in the similar way as in the previous lemma, we can change $I$ to $I'$ such that $J + I' > I'$. One can note that in this process the separable part will not be changed. Therefore, we can say that*

$$A = \sum_{I' \in (T)_m} \sum_{J \in T,\ J + I' > I'} \overline{t^J r_{I',J}^2 \frac{dt_{I'}}{t_{I'}}}, \quad \text{for some } r_{I',J} \in F.$$

In the next theorem, we compute the group $W_1(H_2^{m+1}(F))$.

**Theorem 2.5.** *Let $F$ be a field of characteristic $2$ which is complete with respect to a discrete valuation $v : F \to \mathbb{Z}$. Let $\pi$ be an uniformizer and $\overline{F}$ the residue field. If $\overline{F} \subset F$, then each element $\varphi \in H_2^{m+1}(F)$ decomposes as follows:*

$$\varphi = \varphi_1 + \psi + \varphi_2 \wedge \overline{\frac{d\pi}{\pi}},$$

where $\varphi_1 \in H_2^{m+1}(\overline{F})$, $\varphi_2 \in H_2^m(\overline{F})$ and $\psi = \sum_{I \in (T)_m} \sum_{J \in T, J+I>I} \overline{t^J r_{I,J}^2 \frac{dt_I}{t_I}} \in \mathcal{R}$ with $r_{I,J} \in \pi^{-1}R$.

Here $H_2^{m+1}(\overline{F})$ and $H_2^m(\overline{F}) \wedge \overline{\frac{d\pi}{\pi}}$ are considered as subgroups of $H_2^{m+1}(F)$. The elements $\varphi_1$, $\varphi_2$ and the $r_{I,J}$ are unique for each $\varphi$. Moreover, $\varphi_1$ and $\varphi_2$ are unique over $\overline{F}$ and $F$ both. Furthermore, we have a split exact sequence

$$0 \to H_2^{m+1}(\overline{F}) \to H_2^{m+1}(F) \to \mathcal{R} \oplus H_2^m(\overline{F}) \wedge \overline{\frac{d\pi}{\pi}} \to 0.$$

*Proof.* We start by proving the existence of the decomposition, and then we treat the uniqueness. We keep the same 2-basis $C$ of $F$ as constructed in page 4. Recall that since $F$ is complete, we have $F = \overline{F}((\pi))$. Let $R = \overline{F}[\pi^{-1}]$ and $\varphi \in H_2^{m+1}(F)$.

**1. Existence of the decomposition:** By Lemma 2.3, we have

$$\varphi = \sum_{I \in (T)_m} \sum_{J \in T, J+I \geq I} \overline{t^J b_{I,J}^2 \frac{dt_I}{t_I}}, \quad b_{I,J} \in F.$$

We have $b_{I,J} = r_{I,J} + \overline{f}_{I,J} + b'_{I,J}$ for some $r_{I,J} \in \pi^{-1}R$, $\overline{f}_{I,J} \in \overline{F}$, and $v(b'_{I,J}) > 0$ for all $I, J$. Since $v(t^J b'^2_{I,J}) > 0$, it follows that $t^J b'^2_{I,J} \in \wp(F)$. Consequently, $\overline{t^J b'^2_{I,J} \frac{dt_I}{t_I}} = 0 \in H_2^{m+1}(F)$. Now we treat two cases: $J = 0$ or $J \neq 0$.

**Case 1:** Suppose $J = 0$. Then, we have:

$$\overline{t^J b_{I,J}^2 \frac{dt_I}{t_I}} = \overline{b_{I,0} \frac{dt_I}{t_I}} = \overline{(r_{I,0} + \overline{f}_{I,0} + b'_{I,0}) \frac{dt_I}{t_I}} = \overline{(r_{I,0} + \overline{f}_{I,0}) \frac{dt_I}{t_I}}.$$

Moreover, $\overline{\overline{f}_{I,0} \frac{dt_I}{t_I}} \in H_2^{m+1}(\overline{F})$ or $H_2^m(\overline{F}) \wedge \overline{\frac{d\pi}{\pi}}$ according to $I_\infty = 0$ or 1. Since $r_{I,0} = \sum_{i=1}^{l} \frac{c_i}{\pi^i}$ for suitable $c_i \in \overline{F}$ and $l \in \mathbb{N}_0$, it suffices to prove that $\overline{\frac{c}{\pi^k} \frac{dt_I}{t_I}} \in \mathcal{R}$ for each $c \in \overline{F}$ and $k \geq 1$. We will do this by induction on $k$ after treating the case where $k$ is odd.

**Step 1:** Suppose $k$ odd. Using the 2-basis $C$, we have $c = \sum_{J' \in T, J'_\infty = 0} b_{J'}^2 t^{J'}$ for some $b_{J'} \in \overline{F}$. Therefore,

$$\overline{\frac{c}{\pi^k} \frac{dt_I}{t_I}} = \sum_{J' \in T, J'_\infty = 0} \overline{\frac{b_{J'}^2 t^{J'}}{\pi^k} \frac{dt_I}{t_I}} = \sum_{J' \in T, J'_\infty = 0} \overline{\frac{b_{J'}^2 t^{J'} \pi}{\pi^{k+1}} \frac{dt_I}{t_I}}$$

$$= \sum_{J'' \in T, J'' \neq 0} \overline{b_{J''}^2 t^{J''} \frac{dt_I}{t_I}} \quad \text{(where } t^{J''} = t^{J'} \pi \text{ and } b_{J''} = \frac{b_{J'}}{\pi^{\frac{(k+1)}{2}}} \in \pi^{-1}R\text{)}$$

$$= \sum_{I' \in (T)_m} \sum_{\substack{J'', \\ J''+I'>I'}} \overline{b_{J''}^2 t^{J''} \frac{dt_{I'}}{t_{I'}}} \in \mathcal{R}.$$

For the last equality, as all $J'' \neq 0$ we can change $I$ to $I'$ such that $J'' + I' > I'$, as in Lemma 2.3. Note that in this process $b_{J''} \in \pi^{-1}R$ will not change.

**Step 2:** Now we start induction on $k$. The case $k = 1$ was considered in Step 1. Let us suppose that $\overline{\frac{c}{\pi^k} \frac{dt_I}{t_I}} \in \mathcal{R}$ for all $k \leq l \in \mathbb{N}_0$. Now we will prove this for $k = l + 1$. By Step 1, we may suppose that

$l + 1 = 2i$ for some $i \geq 1$. Now putting $c = \sum_{J' \in T, J'_\infty = 0} b^2_{J'} t^{J'}$, we get:

$$\overline{\frac{c}{\pi^{l+1}} \frac{dt_I}{t_I}} = \sum_{J' \in T, J' \neq 0} \overline{t^{J'} \frac{b^2_{J'}}{\pi^{2i}} \frac{dt_I}{t_I}} + \overline{\frac{b^2_0}{\pi^{2i}} \frac{dt_I}{t_I}} = \sum_{I' \in (T)_m} \sum_{\substack{J' \\ J'+I'>I'}} \overline{t^{J'} \left(\frac{b_{J'}}{\pi^i}\right)^2 \frac{dt_{I'}}{t_{I'}}} + \overline{\frac{b_0}{\pi^i} \frac{dt_I}{t_I}} \in \mathcal{R}$$

To get the last equality, in the first sum as all $J' \neq 0$ we can change $I$ to $I'$ such that $J' + I' > I'$. The first sum is in $\mathcal{R}$ by definition, and the second one is in $\mathcal{R}$ by induction hypothesis.

Finally, taking Steps 1 and 2, $\overline{r_{I,0} \frac{dt_I}{t_I}} \in \mathcal{R}$ for all $r_{I,0} \in \pi^{-1} R$. Therefore, $\overline{b^2_{I,0} \frac{dt_I}{t_I}}$ can be written as $\varphi_1 + \psi + \varphi_2 \wedge \overline{\frac{d\pi}{\pi}}$, where $\varphi_1 \in H^{m+1}_2(\overline{F}), \varphi_2 \in H^m_2(\overline{F})$ and $\psi \in \mathcal{R}$ for any $b_{I,0} \in F$.

**Case 2:** If $J \neq 0$, then $I$ and $J$ satisfy the relation $J + I > I$. Now we have two possibilities: $J_\infty = 0$ or $J_\infty = 1$.

(i) Let $J_\infty = 0$. In this case $t^J \overline{f}^2_{I,J} \in \overline{F}$. Therefore, $\overline{t^J \overline{f}^2_{I,J} \frac{dt_I}{t_I}} \in H^{m+1}_2(\overline{F})$ or $H^m_2(\overline{F}) \wedge \overline{\frac{d\pi}{\pi}}$ according to $I_\infty = 0$ or $1$. Finally, we get

$$\overline{t^J b^2_{I,J} \frac{dt_I}{t_I}} = \overline{t^J r^2_{I,J} \frac{dt_I}{t_I}} + \overline{t^J \overline{f}^2_{I,J} \frac{dt_I}{t_I}} \in H^{m+1}_2(\overline{F}) + \mathcal{R} + H^m_2(\overline{F}) \wedge \overline{\frac{d\pi}{\pi}}.$$

(ii) Let $J_\infty = 1$. In this case $v(t^J \overline{f}^2_{I,J}) > 0$, which implies $t^J \overline{f}^2_{I,J} \in \wp(F)$. Consequently, $\overline{t^J \overline{f}^2_{I,J} \frac{dt_I}{t_I}} = 0 \in H^{m+1}_2(F)$. So we have,

$$\overline{t^J b^2_{I,J} \frac{dt_I}{t_I}} = \overline{t^J r^2_{I,J} \frac{dt_I}{t_I}} \in \mathcal{R}.$$

Thus, for each $\varphi \in H^{m+1}_2(F)$ we have $\varphi = \varphi_1 + \psi + \varphi_2 \wedge \overline{\frac{d\pi}{\pi}}$, where $\varphi_1 \in H^{m+1}_2(\overline{F}), \varphi_2 \in H^m_2(\overline{F})$ and $\psi \in \mathcal{R}$.

**2. Uniqueness of the decomposition:**

Let $\varphi = \varphi_1 + \psi + \varphi_2 \wedge \overline{\frac{d\pi}{\pi}} = 0$ be such that $\varphi_1 \in H^{m+1}_2(\overline{F}), \varphi_2 \in H^m_2(\overline{F})$ and $\psi = \sum_{I \in (T)_m} \sum_{J, J+I>I} \overline{t^J r^2_{I,J} \frac{dt_I}{t_I}} \in \mathcal{R}$, where $r_{I,J} \in \pi^{-1} R$.

Let $L/F$ be a finite separable extension such that $(\varphi_1)_L = 0$ and $(\varphi_2)_L = 0$. The discrete valuation of $F$ extends to $L$ uniquely, and $L$ will be also complete as $F$ is complete. Since the extension $L/F$ is separable, the 2-basis of $F$ remains a 2-basis of $L$ by [1, Lemma 2.1]. Clearly, we have $R \subset R_L$. By [4, Proposition 1.5], we can say that $\sum_{J, J+I>I} t^J r^2_{I,J} = 0$ over $L$, and thus $\sum_{J, J+I>I} t^J r^2_{I,J} = 0$ over $F$. Since $\{t_i \mid i \in \mathcal{I}\}$ is a 2-basis of $F$, it follows that $r_{I,J} = 0$ for all $I \in (T)_m$ and for all $J \in T$ such that $J + I > I$. Hence, $\varphi_1 + \varphi_2 \wedge \overline{\frac{d\pi}{\pi}} = 0$ over $F$. Now the residue map $\zeta$, given in Definition 6.12, concludes that

$$\zeta\left(\varphi_1 + \varphi_2 \wedge \overline{\frac{d\pi}{\pi}}\right) = \varphi_2 = 0 \text{ in } H^m_2(\overline{F}).$$

By $\overline{F} \subset F$, we get $\varphi_2 = 0 \in H^m_2(F)$, and thus $\varphi_1 = 0 \in H^{m+1}_2(F)$. Thus, the uniqueness is proved. □

**Corollary 2.6.** *Let $F$ be a field of characteristic 2 which is complete with respect to a discrete valuation $v : F \to \mathbb{Z}$. Let $\pi$ be an uniformizer and $\overline{F}$ the residue field of $F$. Then, $\mathcal{R}$ as a subgroup of $W_1(H_2^{m+1}(F))$ is unique up to isomorphism irrespective of different liftings of the 2-basis $\{\overline{t_i} \mid i \in \mathcal{I}\}$ of $\overline{F}$. Consequently, the group*

$$W_1(H_2^{m+1}(F)) := coker(H_2^{m+1}(\overline{F}) \to H_2^{m+1}(F))$$

*does not depend on the lifting of the 2-basis.*

*Proof.* Let $A = \{t_i \mid i \in \mathcal{I}\}$ and $A_1 = \{t'_i \mid i \in \mathcal{I}\}$ be two different liftings in $F$ of the 2-basis $\{\overline{t_i} \mid i \in \mathcal{I}\}$ of $\overline{F}$. By Theorem 6.6, we know that there are two unique monomorphisms of fields $\alpha : \overline{F} \hookrightarrow F$ with $\alpha(\overline{t_i}) = t_i$ and $\beta : \overline{F} \hookrightarrow F$ with $\beta(\overline{t_i}) = t'_i$, for all $i \in \mathcal{I}$. Hence $\alpha(\overline{F}) \simeq \beta(\overline{F})$. So $H_2^{m+1}(\alpha(\overline{F})) \simeq H_2^{m+1}(\beta(\overline{F}))$ and $\alpha(\overline{F})[\pi^{-1}] \simeq \beta(\overline{F})[\pi^{-1}]$. Therefore, the isomorphic class of $R = \overline{F}[\pi^{-1}] \subseteq F$ does not depend on the lifting of the 2-basis of $\overline{F}$. By the previous theorem each element of $\mathcal{R} \subset W_1(H_2^{m+1}(F))$ has a unique decomposition. Moreover, the elements of $\mathcal{R}$ only depend upon $R$ and it follows that the isomorphic class of $\mathcal{R}$ is same even after different liftings of the 2-basis. Hence, the corollary is proved.

*Another Proof:* By definition, $W_1(H_2^{m+1}(F)) := coker\ (\alpha' : H_2^{m+1}(\overline{F}) \to H_2^{m+1}(F))$, where $\alpha'$ is induced from a Teichmüller lifting $\alpha$. Moreover, Corollary 6.10 suggests that $\alpha'$ does not depend on the lifting of the 2-basis of $\overline{F}$. Therefore, $W_1(H_2^{m+1}(F))$ does not depend on the lifting of the 2-basis. □

**Remark 2.7.** *Let $F$ be a field of characteristic 2 which is complete with respect to a valuation $v : F \to \mathbb{Z}$, and $\overline{F}$ the residue field. Then, Theorem 2.5 implies that*

$$0 \to H_2^{m+1}(\overline{F}) \xrightarrow{i_4} H_2^{m+1}(F) \xrightarrow{\partial_v} W_1(H_2^{m+1}(F)) \to 0$$

*is a split exact sequence, where $i_4$ is induced by the Teichmüller lifting, and $\partial_v : H_2^{m+1}(F) \to W_1(H_2^{m+1}(F))$, is the residue map defined in page 3.*

## 3. The filtration of $H_2^{m+1}(\mathcal{F}(x))$

Let $\mathcal{F}$ be a field of characteristic 2 with a fixed 2-basis $\{t_i \mid i \in \mathcal{I}\}$, where $\mathcal{I}$ is an ordered set. We are taking $\mathcal{F}$ as an imperfect field, because if $\mathcal{F}$ is perfect, then the exact sequence of Theorem 4.10 will be obvious. Let $F = \mathcal{F}(x)$ be the rational function field in one variable $x$ over $\mathcal{F}$. In this section we make a filtration $(L_d)_{d \in \mathbb{N}}$ of $H_2^{m+1}(\mathcal{F}(x))$, as was done by Milnor for Witt group in characteristic not 2 [15]. For $d \in \mathbb{N}_0$, let us define the sets $\mathcal{F}[x]_{\leq d} := \{f(x) \in \mathcal{F}[x] \mid \deg f(x) \leq d\}$ and $\mathcal{F}[x]_{<d} := \{f(x) \in \mathcal{F}[x] \mid \deg f(x) < d\}$.

**Definition 3.1.** *Let $d \in \mathbb{N}$. We define the subgroup $L_d$ of $H_2^{m+1}(\mathcal{F}(X))$ as follows:*

*(1) For $d > 0$, $L_d$ is the subgroup generated by the elements $\overline{\dfrac{h}{u^e} \dfrac{df_1}{f_1} \wedge \ldots \wedge \dfrac{df_m}{f_m}}$, with $h \in \mathcal{F}[x], e \geq 0$ and $f_i, u \in \mathcal{F}[x]_{\leq d} \setminus \{0\}$.*

*(2) $L_0$ is the subgroup generated by $\overline{h \dfrac{dc_1}{c_1} \wedge \ldots \wedge \dfrac{dc_{m-1}}{c_{m-1}} \wedge \dfrac{df_m}{f_m}}$, where $c_1, \cdots, c_{m-1} \in \mathcal{F}^*$, $f_m \in \mathcal{F}^* \cup \{x\}$ and $h \in \mathcal{F}[x]$ (resp. $h \in x\mathcal{F}[x]$) if $f_m \in \mathcal{F}^*$ (resp. if $f_m = x$).*

Clearly, $H_2^{m+1}(\mathcal{F}(x)) = \cup_{d \geq 0} L_d$, as any generator of $H_2^{m+1}(\mathcal{F}(x))$ lies in $L_d$ for some large $d$. Moreover, $L_0 \subset L_1 \subset L_2 \subset \ldots$ is an increasing chain of subgroups. To give some properties on the groups $L_d$, we need the following relation:

(5) $$\dfrac{d(a+b)}{(a+b)} = \dfrac{a}{a+b} \dfrac{da}{a} + \dfrac{b}{a+b} \dfrac{db}{b}.$$

**Lemma 3.2.** (i) If $0 \neq f, h \in \mathcal{F}[x]$ and $c_1, \cdots, c_{m-1} \in \mathcal{F}^*$, then $\overline{fh \dfrac{dc_1}{c_1} \wedge \ldots \wedge \dfrac{dc_{m-1}}{c_{m-1}} \wedge \dfrac{df}{f}} \in L_0$.

(ii) For any $d \in \mathbb{N}_0$, the elements $\overline{\dfrac{h_1}{u_1^e} \dfrac{dc_1}{c_1} \wedge \ldots \wedge \dfrac{dc_m}{c_m}}$ and $\overline{\dfrac{h_2}{u_2^{e'}} \dfrac{de_1}{e_1} \wedge \ldots \wedge \dfrac{de_{m-1}}{e_{m-1}} \wedge \dfrac{dx}{x}}$ generate $L_d$, where $e, e' \in \mathbb{N}$, $c_i, e_j \in \mathcal{F}^*$, $h_1, h_2 \in \mathcal{F}[x]$ and $u_1, u_2$ are products of nonzero polynomials of degree $\leq d$.

(iii) Let $g \in \mathcal{F}[x]$, $p \in \mathcal{F}[x]^*$ and $f, u, f_i \in \mathcal{F}[x]_{\leq d-1} \setminus \{0\}$ for $d \geq 1$. Then, $\overline{\dfrac{pg}{u^e} \dfrac{df_1}{f_1} \wedge \ldots \wedge \dfrac{df_{m-1}}{f_{m-1}} \wedge \dfrac{d(pf)}{(pf)}} \in L_{d-1}$.

*Proof.* (i) Let $f = \sum_{i=0}^{k} a_i x^i \in \mathcal{F}[x]$. Then, using the relation (5), we get:

$$\overline{fh \dfrac{dc_1}{c_1} \wedge \ldots \wedge \dfrac{dc_{m-1}}{c_{m-1}} \wedge \dfrac{df}{f}} = \sum_{i=0}^{k} \overline{ha_i x^i \dfrac{dc_1}{c_1} \wedge \ldots \wedge \dfrac{dc_{m-1}}{c_{m-1}} \wedge \dfrac{d(a_i x^i)}{a_i x^i}}$$

$$= \sum_{i=0}^{k} \overline{ha_i x^i \dfrac{dc_1}{c_1} \wedge \ldots \wedge \dfrac{dc_{m-1}}{c_{m-1}} \wedge \left(\dfrac{da_i}{a_i} + \dfrac{dx^i}{x^i}\right)}$$

$$= \sum_{i=0}^{k} \overline{ha_i x^i \dfrac{dc_1}{c_1} \wedge \ldots \wedge \dfrac{dc_{m-1}}{c_{m-1}} \wedge \dfrac{da_i}{a_i}} + \sum_{\substack{i \in \{1, \cdots, k\} \\ i \text{ odd}}} \overline{ha_i x^i \dfrac{dc_1}{c_1} \wedge \ldots \wedge \dfrac{dc_{m-1}}{c_{m-1}} \wedge \dfrac{dx}{x}} \in L_0.$$

(ii) Let $u_1, u_2$ be products of nonzero polynomials of degree $\leq d$. By partial fraction method we have $\dfrac{h_1}{u_1^e} = \sum_i \dfrac{h_{i,1}}{u_{i,1}^{e_{i,1}}}$, where $u_{i,1} \in \mathcal{F}[x]_{\leq d}$, and $e_{i,1} \geq 0$. We do a similar decomposition for $\dfrac{h_2}{u_2^{e'}}$. Now by Definition 3.1, we can say that the elements $\overline{\dfrac{h_1}{u_1^e} \dfrac{dc_1}{c_1} \wedge \ldots \wedge \dfrac{dc_m}{c_m}}$ and $\overline{\dfrac{h_2}{u_2^{e'}} \dfrac{de_1}{e_1} \wedge \ldots \wedge \dfrac{de_{m-1}}{e_{m-1}} \wedge \dfrac{dx}{x}}$ belong to $L_d$.

Conversely, let $\overline{\dfrac{h}{u^e} \dfrac{df_1}{f_1} \wedge \ldots \wedge \dfrac{df_m}{f_m}}$ be a generator of $L_d$, for some $f_i, u \neq 0 \in \mathcal{F}[x]_{\leq d}$. Let $f_m = \sum_{i=0}^{d} a_i x^i$, where $a_i \in \mathcal{F}$. Then, using the relation (5), we get:

$$\overline{\dfrac{h}{u^e} \dfrac{df_1}{f_1} \wedge \ldots \wedge \dfrac{df_m}{f_m}} = \sum_{i=0}^{d} \overline{\dfrac{ha_i x^i}{f_m u^e} \dfrac{df_1}{f_1} \wedge \ldots \wedge \dfrac{df_{m-1}}{f_{m-1}} \wedge \dfrac{d(a_i x^i)}{a_i x^i}}$$

$$= \sum_{i=0}^{d} \overline{\dfrac{h'_i}{(uf_m)^e} \dfrac{df_1}{f_1} \wedge \ldots \wedge \dfrac{df_{m-1}}{f_{m-1}} \wedge \left(\dfrac{da_i}{a_i} + \dfrac{dx^i}{x^i}\right)}$$

where $h'_i = ha_i x^i f_m^{e-1}$. Then, we will have to repeat the process for $f_{m-1}, \ldots f_1$. At last the denominator of the separable part will be $(uf_1 \ldots f_m)^e$, where $\deg u, \deg f_i \leq d$. Moreover, $dx^i = 0$ if $i$ is even and $\dfrac{dx^i}{x^i} = \dfrac{dx}{x}$ if $i$ is odd and $dx \wedge dx = 0$.

Therefore, $L_d$ can be generated by the forms $\overline{\dfrac{h_1}{u_1^e} \dfrac{dc_1}{c_1} \wedge \ldots \wedge \dfrac{dc_m}{c_m}}$ and $\overline{\dfrac{h_2}{u_2^{e'}} \dfrac{de_1}{e_1} \wedge \ldots \wedge \dfrac{de_{m-1}}{e_{m-1}} \wedge \dfrac{dx}{x}}$, where $h_1, h_2 \in \mathcal{F}[x], c_i, e_j \in \mathcal{F}^*$ and $u_1, u_2$ are product of nonzero polynomials of degree $\leq d$.

(iii) Let us write $pf = \sum_i b_i x^i \in \mathcal{F}[x]$. Then, using the relation (5), we obtain:

$$\overline{\frac{pg}{u^e} \frac{df_1}{f_1} \wedge \ldots \wedge \frac{df_{m-1}}{f_{m-1}} \wedge \frac{d(pf)}{(pf)}} = \sum_i \overline{\frac{pgb_i x^i}{pfu^e} \frac{df_1}{f_1} \wedge \ldots \wedge \frac{df_{m-1}}{f_{m-1}} \wedge \frac{d(b_i x^i)}{b_i x^i}}$$

$$= \sum_i \overline{\frac{h'_i}{(fu)^e} \frac{df_1}{f_1} \wedge \ldots \wedge \frac{df_{m-1}}{f_{m-1}} \wedge \left(\frac{db_i}{b_i} + \frac{dx^i}{x^i}\right)}$$

where $h'_i = gb_i x^i f^{e-1}$. We repeat the process for each $f_i$, $1 \leq i \leq m-1$, and using the same arguments as in (ii), we get that $\overline{\frac{pg}{u^e} \frac{df_1}{f_1} \wedge \ldots \wedge \frac{df_{m-1}}{f_{m-1}} \wedge \frac{d(pf)}{(pf)}}$ is a sum of these two types of elements:

$$\overline{\frac{h_1}{u_1^e} \frac{dc_1}{c_1} \wedge \ldots \wedge \frac{dc_m}{c_m}} \quad \text{and} \quad \overline{\frac{h_2}{u_2^{e'}} \frac{de_1}{e_1} \wedge \ldots \wedge \frac{de_{m-1}}{e_{m-1}} \wedge \frac{dx}{x}},$$

where $h_1, h_2 \in \mathcal{F}[x], c_i, e_j \in \mathcal{F}^*$ and $u_1, u_2$ are product of nonzero polynomials of degree $\leq d-1$. Hence, by part (ii) of this lemma, $\overline{\frac{pg}{u^e} \frac{df_1}{f_1} \wedge \ldots \wedge \frac{df_{m-1}}{f_{m-1}} \wedge \frac{d(pf)}{pf}} \in L_{d-1}$. □

The aim of the next lemma is to prove that for $p \in \mathcal{F}[x]$ monic irreducible of degree $d$, the residue map $\partial_p$ extends to $L_d/L_{d-1}$.

**Lemma 3.3.** *Let $p \in \mathcal{F}[x]$ be a monic irreducible polynomial of degree $d$, $v_p$ the $p$-adic valuation on $\mathcal{F}(x)$ and $s, r_1, \cdots, r_m \in \mathcal{F}(x)$.*

*(1) If $v_p(s) \geq 0$ and $v_p(r_i) = 0$ for each $i$, then $\partial_p\left(\overline{s \frac{dr_1}{r_1} \wedge \ldots \wedge \frac{dr_m}{r_m}}\right) = 0$.*

*(2) We have $\partial_p(\varphi) = 0$ for any $\varphi \in L_{d-1}$, and thus $\partial_p$ extends to a homomorphism on $L_d/L_{d-1}$.*

*Proof.* We keep the same notations and hypotheses as in the lemma.

(1) Recall that $F$ denotes $\mathcal{F}(x)$. Moreover, $F_p$ and $\overline{F}_p$ are the completion and the residue field of $F$, respectively. Since $F_p \simeq \overline{F}_p((p))$, we have $r_i = r_{i,0} + pr_{i,1}$, where $r_{i,0} \neq 0 \in \overline{F}_p$ and $r_{i,1} \in F_p$ with $v_p(r_{i,1}) \geq 0$. Since $r_1$ is a unit, we get $\frac{sr_{1,0}}{r_1} = s_{1,0} + ps_{1,1} \in F_p$ for some $s_{1,0} \in \overline{F}_p$ and $s_{1,1} \in F_p$ such that $v_p(s_{1,1}) \geq 0$. We define the recurrence relation for $i = 2, 3, \ldots, m$:

$$\frac{s_{i-1,0} r_{i,0}}{r_i} = s_{i,0} + ps_{i,1}, \text{ where } s_{i,0} \in \overline{F}_p \text{ and } s_{i,1} \in F_p \text{ with } v_p(s_{i,1}) \geq 0.$$

Set $s_{0,0} = s$. We have $v_p\left(\frac{s_{i-1,0} pr_{i,1}}{r_i}\right) > 0$ and $v_p(ps_{i,1}) > 0$ for $i = 1, \ldots, m$. Hence, $\frac{s_{i-1,0} pr_{i,1}}{r_i}, ps_{i,1} \in \wp(F_p)$ for $i = 1, \ldots, m$. Using the above calculations, we get:

$$\overline{s \frac{dr_1}{r_1} \wedge \frac{dr_2}{r_2} \wedge \ldots \wedge \frac{dr_m}{r_m}} = \overline{\frac{sr_{1,0}}{r_1} \frac{dr_{1,0}}{r_{1,0}} \wedge \frac{dr_2}{r_2} \wedge \ldots \wedge \frac{dr_m}{r_m}} +$$

$$\overline{\frac{spr_{1,1}}{r_1} \frac{d(pr_{1,1})}{(pr_{1,1})} \wedge \frac{dr_2}{r_2} \wedge \ldots \wedge \frac{dr_m}{r_m}} \quad (\text{using } r_1 = r_{1,0} + pr_{1,1})$$

$$= \overline{\frac{sr_{1,0}}{r_1} \frac{dr_{1,0}}{r_{1,0}} \wedge \frac{dr_2}{r_2} \wedge \ldots \wedge \frac{dr_m}{r_m}} \quad \left(\text{as } \frac{spr_{1,1}}{r_1} \in \wp(F_p)\right).$$

Putting $\dfrac{sr_{1,0}}{r_1} = s_{1,0} + ps_{1,1}$, we obtain :

$$\overline{s\dfrac{dr_1}{r_1} \wedge \dfrac{dr_2}{r_2} \wedge \ldots \wedge \dfrac{dr_m}{r_m}} = \overline{(s_{1,0} + ps_{1,1})\dfrac{dr_{1,0}}{r_{1,0}} \wedge \dfrac{dr_2}{r_2} \wedge \ldots \wedge \dfrac{dr_m}{r_m}}$$

$$= \overline{s_{1,0}\dfrac{dr_{1,0}}{r_{1,0}} \wedge \dfrac{dr_2}{r_2} \wedge \ldots \wedge \dfrac{dr_m}{r_m}} \quad \text{(as } ps_{1,1} \in \wp(F_p)\text{)}.$$

Repeating the process for $r_2, \ldots r_m$, finally we get in $H_2^{m+1}(F_p)$ :

$$\overline{s\dfrac{dr_1}{r_1} \wedge \ldots \wedge \dfrac{dr_m}{r_m}} = \overline{s_{m,0}\dfrac{dr_{1,0}}{r_{1,0}} \wedge \ldots \wedge \dfrac{dr_{m,0}}{r_{m,0}}} \in \mathrm{Im}(H_2^{m+1}(\overline{F}_p) \to H_2^{m+1}(F_p)),$$

as $s_{m,0}, r_{1,0}, \ldots, r_{m,0} \in \overline{F}_p = \mathcal{F}(p)$. Therefore, $\partial_p\left(\overline{s\dfrac{dr_1}{r_1} \wedge \ldots \wedge \dfrac{dr_m}{r_m}}\right) = 0$.

(2) (i) Suppose $d > 1$. If $\varphi \in L_{d-1}$, then $\varphi$ is a finite sum of the following type elements $\overline{\dfrac{h}{u^e}\dfrac{df_1}{f_1} \wedge \ldots \wedge \dfrac{df_m}{f_m}}$, where $h \in \mathcal{F}[x]$ and $f_i, u \in \mathcal{F}[x]_{<d} \setminus \{0\}$. Hence $v_p(u) = v_p(f_i) = 0$ and $v_p\left(\dfrac{h}{u^e}\right) \geq 0$. Therefore, by part (1), we get $\partial_p(\varphi) = 0$.

(ii) Suppose $d = 1$ and let $\varphi \in L_0$. Then, $\varphi$ is a finite sum of the following type elements $\overline{h_1\dfrac{dc_1}{c_1} \wedge \ldots \wedge \dfrac{dc_m}{c_m}}$ and $\overline{xh_2\dfrac{de_1}{e_1} \wedge \ldots \wedge \dfrac{de_{m-1}}{e_{m-1}} \wedge \dfrac{dx}{x}}$, where $c_i, e_i \in \mathcal{F}^*, h_1, h_2 \in \mathcal{F}[x]$.

- If $p \neq x$, then $v_p(c_i) = v_p(e_i) = v_p(x) = 0$ and $v_p(h_1), v_p(xh_2) \geq 0$ . Hence, by part (1), we get $\partial_p(\varphi) = 0$.

- If $p = x$, then by part (1), $\partial_x\left(\overline{h_1\dfrac{dc_1}{c_1} \wedge \ldots \wedge \dfrac{dc_m}{c_m}}\right) = 0$ and $v_x(xh_2) > 0$, so $xh_2 \in \wp(F_x)$. Hence, we get
$$\partial_x\left(\overline{xh_2\dfrac{de_1}{e_1} \wedge \ldots \wedge \dfrac{de_{m-1}}{e_{m-1}} \wedge \dfrac{dx}{x}}\right) = 0. \text{ Thus, } \partial_p(L_{d-1}) = 0 \text{ when } \deg p = d.$$

Hence the lemma is proved. □

**Definition 3.4.** *(1) For each monic and irreducible polynomial $p$, let $S_p$ be the subgroup of $H_2^{m+1}(\mathcal{F}(x))$ generated by the elements $\overline{\dfrac{h}{p^s}\dfrac{dc_1}{c_1} \wedge \ldots \wedge \dfrac{dc_m}{c_m}}$, with $c_i \in \mathcal{F}^*, h \in \mathcal{F}[x]$ and $s \geq 0$.*

*(2) Similarly, let $S_{\frac{1}{x}}$ be the subgroup of $H_2^{m+1}(\mathcal{F}(x))$ generated by the elements $\overline{hx\dfrac{dc_1}{c_1} \wedge \ldots \wedge \dfrac{dc_m}{c_m}}$, with $c_i \in \mathcal{F}^*$ and $h \in \mathcal{F}[x]$.*

In the following lemma we find some other generators of $L_d$ using the groups $S_p$.

**Lemma 3.5.** *(1) For each $d \geq 1$, we have $L_d = \displaystyle\sum_{p,\,\deg p \leq d}\left(S_p + S_p \wedge \overline{\dfrac{dx}{x}}\right)$, where $S_p \wedge \overline{\dfrac{dx}{x}}$ is the subgroup of $H_2^{m+1}(\mathcal{F}(x))$ generated by the forms $\overline{\dfrac{h}{p^s}\dfrac{dc_1}{c_1} \wedge \ldots \wedge \dfrac{dc_{m-1}}{c_{m-1}} \wedge \dfrac{dx}{x}}$ such that $h \in \mathcal{F}[x], s \geq 0$ and $c_i \in \mathcal{F}^*$.*

(2) $H_2^{m+1}(\mathcal{F}(x)) = \sum_p \left( S_p + S_p \wedge \overline{\frac{dx}{x}} \right)$, where $p$ varies over all monic irreducible polynomials of $\mathcal{F}[x]$.

(3) For each $d \geq 1$, $L_d = \sum_{p, \deg p = d} \left( S_p + S_p \wedge \overline{\frac{dx}{x}} \right) + L_{d-1}$.

*Proof.* (1) It is clear that $S_p \subseteq L_d$ and $S_p \wedge \overline{\frac{dx}{x}} \subseteq L_d$ for $\deg p \leq d$. Conversely, by Lemma 3.2(ii), $L_d$ is generated by the following two types of forms

$$\overline{\frac{h_1}{u_1^e} \frac{dc_1}{c_1} \wedge \ldots \wedge \frac{dc_m}{c_m}} \text{ and } \overline{\frac{h_2}{u_2^{e'}} \frac{de_1}{e_1} \wedge \ldots \wedge \frac{de_{m-1}}{e_{m-1}} \wedge \frac{dx}{x}},$$

where $c_i, e_i \in \mathcal{F}^*$ and $u_1, u_2$ are products of nonzero polynomials of degree $\leq d$. Now decomposing $\frac{h_1}{u_1^e}$ using partial fractions, we conclude that $\overline{\frac{h_1}{u_1^e} \frac{dc_1}{c_1} \wedge \ldots \wedge \frac{dc_m}{c_m}}$ belongs to $\sum_{p, \deg p \leq d} S_p$. Similarly, the form $\overline{\frac{h_2}{u_2^{e'}} \frac{de_1}{e_1} \wedge \ldots \wedge \frac{de_{m-1}}{e_{m-1}} \wedge \frac{dx}{x}}$ belongs to $\sum_{p, \deg p \leq d} S_p \wedge \overline{\frac{dx}{x}}$. Hence, we get $L_d = \sum_{p, \deg p \leq d} \left( S_p + S_p \wedge \overline{\frac{dx}{x}} \right)$.

(2) This is a consequence of statement (1) and the fact that $H_2^{m+1}(\mathcal{F}(x)) = \cup_{d \geq 0} L_d$.

(3) From statement (1), we have

$$L_d = \sum_{p, \deg p \leq d} \left( S_p + S_p \wedge \overline{\frac{dx}{x}} \right) = \sum_{p, \deg p = d} \left( S_p + S_p \wedge \overline{\frac{dx}{x}} \right) + \sum_{p, \deg p < d} \left( S_p + S_p \wedge \overline{\frac{dx}{x}} \right)$$

$$= \sum_{p, \deg p = d} \left( S_p + S_p \wedge \overline{\frac{dx}{x}} \right) + L_{d-1}.$$

□

Now our aim is to generate $L_d$ using $S_p$ and $S_p \wedge \overline{\frac{dp}{p}}$, and its reason will be visible in Section 4. To this end, we find the relations between $S_p \wedge \overline{\frac{dx}{x}}$ and $S_p \wedge \overline{\frac{dp}{p}}$. Let us give some notations.

**Notations:** Let $p = \sum_{i=0}^{d} p_i x^{d-i} \in \mathcal{F}[x]$ with $p_0 = 1$. Since $\{t_i \mid i \in \mathcal{I}\}$ is a 2-basis of $\mathcal{F}$, we have $p_j = \sum_{I \in T} t^I p_{j,I}^2$ for suitable $p_{j,I} \in \mathcal{F}$, $0 \leq j \leq d$. Note that $p_{0,0} = 1$ and $p_{0,I} = 0$ for all $I \neq 0 \in T$.

We consider the set $\mathcal{I}_1 = \{i \in \mathcal{I} \mid \exists I \in T \text{ such that } i \in \tilde{I} \text{ and } p_{j,I} \neq 0 \text{ for some } j\}$. Then, $\mathcal{I}_1$ is a finite subset of $\mathcal{I}$ and we put $\mathcal{I}_1 = \{j_1, \ldots, j_n\}$ with $j_1 < j_2 < \ldots < j_n$.

Let $T_1 = \{0, 1\}^n$. This set can be treated as a subgroup of $T = \bigoplus_{|\mathcal{I}|} \{0, 1\} \simeq \bigoplus_{|\mathcal{I}|} (\mathbb{Z}/2\mathbb{Z})$. Moreover, $T_1$ is lexicographically ordered like $T$, as introduced in page 4. Therefore, the elements of $T_1$ lie in the following order

$$(0, 0, \ldots, 0) < (1, 0, \ldots, 0) < (0, 1, 0, \ldots, 0) < \ldots$$

If $K \in T_1$ we take $t^K$ the element $t_{j_1}^{K_1} \ldots t_{j_n}^{K_n}$. Now for all $K \in T_1$, let us define :

$$P_K = \begin{cases} \overline{t^K \left(p_{1,K}^2 x^{d-1} + p_{3,K}^2 x^{d-3} + \ldots + p_{d-1,K}^2 x\right)} & \text{if } d \text{ is even}, \\ \overline{t^K \left(p_{0,K}^2 x^d + p_{2,K}^2 x^{d-2} + \ldots + p_{d-1,K}^2 x\right)} & \text{if } d \text{ is odd}. \end{cases}$$

Now we define a matrix $M$ over $\mathcal{F}(p)$, whose $(K_1, K_2)$-th element is $P_{K_1+K_2}$ for all $K_1, K_2 \in T_1$. Note that the order of $M$ is $2^n \times 2^n$.

The next corollary is proved in [4, Corollary 2.7]. But for completeness of the paper, we are providing the proof.

**Corollary 3.6.** *Let $p$ be a monic, irreducible and separable polynomial over $\mathcal{F}$ such that $p \neq x$. Then, the matrix $M$, defined above, is invertible.*

*Proof.* Suppose that $M$ is not invertible. Then, by [4, lemma 2.6], we have $\sum_{K \in T_1} P_K = 0 \in \mathcal{F}(p)$. Since $p$ is separable, the elements of the set $\{t^K \mid K \in T_1\}$ are linearly independent over $\mathcal{F}(p)^2$. By definition, each $P_K$ is of the form $\overline{xt^K(\mathcal{F}(p))^2}$ and $\overline{x}$ is a unit in $\mathcal{F}(p)$. Hence $P_K = 0 \in \mathcal{F}(p)$ for all $K \in T_1$. Now we have two cases depending on $\deg p = d$ is even or odd.

(1) Suppose that $d$ is even. Since $P_K$ is a polynomial of degree $d - 1$, it follows that $p_{j,K} = 0$ for all $K \in T_1$ and $j = 1, 3, \ldots, d - 1$, which implies that $p_j = 0$ for $j = 1, 3, \ldots, d - 1$. Then, $p$ will not be a separable polynomial, a contradiction.

(2) Suppose that $d$ is odd. So $P_K = 0 \in \mathcal{F}(p)$ implies that $\overline{p_{0,K}^2 x^{d-1} + p_{2,K}^2 x^{d-3} + \ldots + p_{d-1,K}^2} = 0$ for all $K \in T_1$. Then, similarly as in the previous case, we have $p_{j,K} = 0$ for all $K \in T_1$ and $j = 0, 2, \ldots, d - 1$. But we know that $p_{0,0} = 1$, which contradicts our assumption.

Consequently, $M$ is invertible over $\mathcal{F}(p)$. □

**Definition 3.7.** Let $p = \sum_{j=0}^{d} p_j x^{d-j}$ be a monic, irreducible and inseparable polynomial of degree $d$. Then, $d$ is even and $p_j = 0$ when $j$ is odd. Since $\sum_{j=0}^{d} \overline{p_j x^{d-j}} = 0 \in \mathcal{F}(p)$, we can find some $l \in \{1, 2, \ldots, n\}$ such that $t_{j_l} \in \mathcal{F}(p)^2(t_{j_1}, t_{j_2}, \ldots, t_{j_{l-1}}, t_{j_{l+1}}, \ldots, t_{j_n})$. Without loss of generality, we may suppose $l = n$. Let $\tilde{S}_p$ be the subgroup of $H_2^{m+1}(F)$ generated by the forms $\overline{\dfrac{h}{p^e} \dfrac{dt_I}{t_I}}$, where $\dfrac{dt_I}{t_I} = \dfrac{dt_{i_1}}{t_{i_1}} \wedge \ldots \wedge \dfrac{dt_{i_m}}{t_{i_m}}$ and $\{t_{i_1}, \ldots, t_{i_m}\} \subseteq \{t_i \mid i \in I \setminus \{j_n\}\} \cup \{x\}, h \in \mathcal{F}[x], e \geq 0$.

The aim of the following proposition is to establish some relations among $S_p + S_p \wedge \overline{\dfrac{dx}{x}}$, $S_p + S_p \wedge \overline{\dfrac{dp}{p}}$ and $\tilde{S}_p + \tilde{S}_p \wedge \overline{\dfrac{dp}{p}}$ in $H_2^{m+1}(\mathcal{F}(x))$. Note that whenever we use $S_p \wedge \overline{\dfrac{dp}{p}}$ or $S_p \wedge \overline{\dfrac{dx}{x}}$, then $S_p$ will be in $H_2^m(\mathcal{F}(x))$ to balance the total wedge.

**Proposition 3.8.** *Let $p \in \mathcal{F}[x]$ be monic and irreducible. We have the following:*

*(1)* $S_p + S_p \wedge \overline{\dfrac{dp}{p}} \subseteq S_p + S_p \wedge \overline{\dfrac{dx}{x}}$.

*(2) If $p$ is separable, then* $S_p + S_p \wedge \overline{\dfrac{dx}{x}} + L_0 \wedge \overline{\dfrac{dx}{x}} = S_p + S_p \wedge \overline{\dfrac{dp}{p}} + L_0 \wedge \overline{\dfrac{dx}{x}}$.

*(3) If $p$ is inseparable, then* $S_p + S_p \wedge \overline{\dfrac{dx}{x}} + L_0 + L_0 \wedge \overline{\dfrac{dx}{x}} = \tilde{S}_p + \tilde{S}_p \wedge \overline{\dfrac{dp}{p}} + L_0 + L_0 \wedge \overline{\dfrac{dx}{x}}$.

*Proof.* (1) We take $p = \sum_{j=0}^{d} p_j x^{d-j}$ with $p_j \in \mathcal{F}$ and $p_0 = 1$. Moreover, $S_p \wedge \overline{\frac{dp}{p}}$ is generated by the following forms

$$\overline{\frac{h}{p^s} \frac{dc_1}{c_1} \wedge \frac{dc_2}{c_2} \wedge \ldots \wedge \frac{dc_{m-1}}{c_{m-1}} \wedge \frac{dp}{p}},$$

where $h \in \mathcal{F}[x], c_i \in \mathcal{F}^*, s \geq 0$. In $H_2^{m+1}(F)$, we have :

$$\overline{\frac{h}{p^s} \frac{dc_1}{c_1} \wedge \frac{dc_2}{c_2} \wedge \ldots \wedge \frac{dc_{m-1}}{c_{m-1}} \wedge \frac{dp}{p}} = \sum_{j=0}^{d} \overline{\frac{p_j x^{d-j} h}{p^{s+1}} \frac{dc_1}{c_1} \wedge \ldots \wedge \frac{dc_{m-1}}{c_{m-1}} \wedge \frac{d(p_j x^{d-j})}{p_j x^{d-j}}}.$$

For $j \in \{0, \ldots, d\}$ with $d-j$ even, we have $\frac{d(p_j x^{d-j})}{p_j x^{d-j}} = \frac{dp_j}{p_j}$. Moreover, if $d-j$ is odd, then $\frac{d(p_j x^{d-j})}{p_j x^{d-j}} = \frac{d(p_j x)}{p_j x}$. Therefore, using the definition of $S_p$, we obtain :

$$\overline{\frac{h}{p^s} \frac{dc_1}{c_1} \wedge \frac{dc_2}{c_2} \wedge \ldots \wedge \frac{dc_{m-1}}{c_{m-1}} \wedge \frac{dp}{p}} \equiv \sum_{d-j \text{ odd}} \overline{\frac{p_j x^{d-j} h}{p^{s+1}} \frac{dc_1}{c_1} \wedge \frac{dc_2}{c_2} \wedge \ldots \wedge \frac{dc_{m-1}}{c_{m-1}} \wedge \frac{d(p_j x)}{p_j x}} \pmod{S_p}.$$

Using $p_j = \sum_{K \in T_1} p_{j,K}^2 t^K$, we get :

$$\overline{\frac{h}{p^s} \frac{dc_1}{c_1} \wedge \frac{dc_2}{c_2} \wedge \ldots \wedge \frac{dc_{m-1}}{c_{m-1}} \wedge \frac{dp}{p}} \equiv \sum_{d-j \text{ odd}} \sum_{K \in T_1} \overline{\frac{t^K p_{j,K}^2 x^{d-j} h}{p^{s+1}} \frac{dc_1}{c_1} \wedge \ldots \wedge \frac{dc_{m-1}}{c_{m-1}} \wedge \frac{d(t^K x)}{t^K x}} \pmod{S_p}$$

$$\equiv \sum_{K \in T_1} \overline{\frac{\left(\sum_{d-j \text{ odd}} t^K p_{j,K}^2 x^{d-j} h\right)}{p^{s+1}} \frac{dc_1}{c_1} \wedge \frac{dc_2}{c_2} \wedge \ldots \wedge \frac{dc_{m-1}}{c_{m-1}} \wedge \frac{d(t^K x)}{t^K x}}.$$

Putting $P_K = \sum_{d-j \text{ odd}} t^K p_{j,K}^2 x^{d-j}$, we have :

$$\overline{\frac{h}{p^s} \frac{dc_1}{c_1} \wedge \frac{dc_2}{c_2} \wedge \ldots \wedge \frac{dc_{m-1}}{c_{m-1}} \wedge \frac{dp}{p}} \equiv \sum_{K \in T_1} \overline{\frac{P_K h}{p^{s+1}} \frac{dc_1}{c_1} \wedge \frac{dc_2}{c_2} \wedge \ldots \wedge \frac{dc_{m-1}}{c_{m-1}} \wedge \frac{d(t^K x)}{(t^K x)}} \pmod{S_p}.$$

Since $t^K \in \mathcal{F}$, the term of right side is in $S_p + S_p \wedge \overline{\frac{dx}{x}}$. Hence

$$\overline{\frac{h}{p^s} \frac{dc_1}{c_1} \wedge \frac{dc_2}{c_2} \wedge \ldots \wedge \frac{dc_{m-1}}{c_{m-1}} \wedge \frac{dp}{p}} \in S_p + S_p \wedge \overline{\frac{dx}{x}}, \text{ where } h \in \mathcal{F}[x] \text{ and } s \geq 0.$$

Therefore, $S_p + S_p \wedge \overline{\frac{dp}{p}} \subseteq S_p + S_p \wedge \overline{\frac{dx}{x}}$ for all monic and irreducible polynomials $p$.

(2) We suppose that $p$ is separable. By statement (1), it suffices to prove that $S_p \wedge \overline{\frac{dx}{x}} \subset S_p + S_p \wedge \overline{\frac{dp}{p}} + L_0 \wedge \overline{\frac{dx}{x}}$. This is obvious for $p = x$. So suppose $p \neq x$. The group $S_p \wedge \overline{\frac{dx}{x}}$ is generated by the elements $\overline{\frac{g}{p^s} \frac{dc_1}{c_1} \wedge \ldots \wedge \frac{dc_{m-1}}{c_{m-1}} \wedge \frac{dx}{x}}$, where $g \in \mathcal{F}[x], c_i \in \mathcal{F}^*$ and $s \geq 0$.

For any $s \in \mathbb{N}$, let $S_{p,s}$ be the subgroup of $S_p$ generated by the forms $\overline{\frac{g}{p^i} \frac{dc_1}{c_1} \wedge \ldots \wedge \frac{dc_m}{c_m}}$, where $g \in \mathcal{F}[x]$ and $0 \leq i \leq s$. Before we prove the inclusion, we do some calculations.

**Calculations:** Let $s \geq 1$ and $h_J \in \mathcal{F}[x]$ for all $J \in T_1$. Using similar process as in (1), we obtain:

$$\overline{\frac{h_J}{p^{s-1}} \frac{dc_1}{c_1} \wedge \ldots \wedge \frac{dc_{m-1}}{c_{m-1}} \wedge \frac{d(pt^J)}{pt^J}} \equiv \sum_{K \in T_1} \overline{\frac{P_K h_J}{p^s} \frac{dc_1}{c_1} \wedge \ldots \wedge \frac{dc_{m-1}}{c_{m-1}} \wedge \frac{d(t^{K+J}x)}{t^{K+J}x}} \pmod{S_p}.$$

Now we add the previous equation for all $J \in T_1$ to get:

$$\sum_{J \in T_1} \overline{\frac{h_J}{p^{s-1}} \frac{dc_1}{c_1} \wedge \ldots \wedge \frac{dc_{m-1}}{c_{m-1}} \wedge \frac{d(pt^J)}{pt^J}} \equiv \sum_{K,J \in T_1} \overline{\frac{P_K h_J}{p^s} \frac{dc_1}{c_1} \wedge \ldots \wedge \frac{dc_{m-1}}{c_{m-1}} \wedge \frac{d(t^{K+J}x)}{t^{K+J}x}} \pmod{S_p}.$$

Putting $L = J + K$, we get:

$$\sum_{J \in T_1} \overline{\frac{h_J}{p^{s-1}} \frac{dc_1}{c_1} \wedge \ldots \wedge \frac{dc_{m-1}}{c_{m-1}} \wedge \frac{d(pt^J)}{pt^J}} \equiv \sum_{L \in T_1} \overline{\frac{\sum_{J \in T_1}(P_{L+J}h_J)}{p^s} \frac{dc_1}{c_1} \wedge \ldots \wedge \frac{dc_{m-1}}{c_{m-1}} \wedge \frac{d(t^L x)}{t^L x}} \pmod{S_p}.$$

For any $g \in \mathcal{F}[x]_{<d}$ we will have to find $h_J \in \mathcal{F}[x]_{<d}$ for all $J \in T_1$ such that these will satisfy the following equations in $\overline{F}_p$

$$\sum_{J \in T_1} P_J h_J = g \quad \text{and} \quad \sum_{J \in T_1} P_{L+J} h_J = 0 \text{ for } L \neq 0 \in T_1.$$

Equivalently, we can say that the column matrix $(h_J)_{J \in T_1}$ will satisfy the following equation over $\overline{F}_p$

$$M.(h_J) = (g, 0, \ldots, 0)^T.$$

Here $M$ is the $2^n \times 2^n$ matrix with $(I, J)$-th element as $P_{I+J}$, for all $I, J \in T_1$. By Corollary 3.6, $M$ is invertible in $\overline{F}_p$. Consequently, we can find the elements $h_J \in \mathcal{F}[x]_{<d}$ such that

$$\sum_{J \in T_1} P_J h_J = g + pk \quad \text{and} \quad \sum_{J \in T_1} P_{L+J} h_J = pk_L, \text{ where } L \neq 0 \in T_1 \text{ for some } k, k_L \in \mathcal{F}[x].$$

This finishes the calculations.

Now for any $g \in \mathcal{F}[x]_{<d}$, we can find some $h_J \in \mathcal{F}[x]_{<d}$ for all $J \in T_1$, such that modulo $\left(S_p + S_{p,s-1} \wedge \overline{\frac{dx}{x}}\right)$, we have:

(6) $$\overline{\frac{g}{p^s} \frac{dc_1}{c_1} \wedge \ldots \wedge \frac{dc_{m-1}}{c_{m-1}} \wedge \frac{dx}{x}} \equiv \sum_{J \in T_1} \overline{\frac{h_J}{p^{s-1}} \frac{dc_1}{c_1} \wedge \ldots \wedge \frac{dc_{m-1}}{c_{m-1}} \wedge \frac{d(t^J p)}{t^J p}}.$$

Now we will start doing induction on $s$ to prove the desired inclusion. Clearly, we have $S_{p,0} \wedge \overline{\frac{dx}{x}} \subseteq L_0 \wedge \overline{\frac{dx}{x}}$. So suppose we have $S_{p,s-1} \wedge \overline{\frac{dx}{x}} \subseteq S_p + S_p \wedge \overline{\frac{dp}{p}} + L_0 \wedge \overline{\frac{dx}{x}}$ for some $s \geq 1$.

Let $g \in \mathcal{F}[x]$, then $\frac{g}{p^s} = \sum_{i=0}^{s} \frac{g_i}{p^i}$ for suitable $g_0 \in \mathcal{F}[x]$ and $g_i \in \mathcal{F}[x]_{<d}$, where $i = 1, \ldots, s$.

Therefore, we get $\overline{\frac{g}{p^s} \frac{dc_1}{c_1} \wedge \ldots \wedge \frac{dc_{m-1}}{c_{m-1}} \wedge \frac{dx}{x}} = \sum_{i=0}^{s} \overline{\frac{g_i}{p^i} \frac{dc_1}{c_1} \wedge \ldots \wedge \frac{dc_{m-1}}{c_{m-1}} \wedge \frac{dx}{x}}$. By induction hypothesis, $\sum_{i=0}^{s-1} \overline{\frac{g_i}{p^i} \frac{dc_1}{c_1} \wedge \ldots \wedge \frac{dc_{m-1}}{c_{m-1}} \wedge \frac{dx}{x}} \in S_p + S_p \wedge \overline{\frac{dp}{p}} + L_0 \wedge \overline{\frac{dx}{x}}$. From (6) we have:

$$\overline{\frac{g_s}{p^s} \frac{dc_1}{c_1} \wedge \ldots \wedge \frac{dc_{m-1}}{c_{m-1}} \wedge \frac{dx}{x}} \in S_p + S_p \wedge \overline{\frac{dp}{p}} + S_{p,s-1} \wedge \overline{\frac{dx}{x}} \subseteq S_p + S_p \wedge \overline{\frac{dp}{p}} + L_0 \wedge \overline{\frac{dx}{x}}.$$

Therefore,
$$\overline{\frac{g}{p^s}\frac{dc_1}{c_1}} \wedge \ldots \wedge \frac{dc_{m-1}}{c_{m-1}} \wedge \frac{dx}{x} \in S_p + S_p \wedge \overline{\frac{dp}{p}} + L_0 \wedge \overline{\frac{dx}{x}} \quad \text{for all } g \in \mathcal{F}[x].$$

By induction, $S_{p,s} \wedge \overline{\frac{dx}{x}} \subseteq S_p + S_p \wedge \overline{\frac{dp}{p}} + L_0 \wedge \overline{\frac{dx}{x}}$ for all $s \geq 0$. Consequently, we prove that

$$S_p \wedge \overline{\frac{dx}{x}} \subseteq S_p + S_p \wedge \overline{\frac{dp}{p}} + L_0 \wedge \overline{\frac{dx}{x}}.$$

Hence, statement (1) and the previous equation prove statement (2).

(3) Let $p \in \mathcal{F}[x]$ be monic, irreducible and inseparable. We have $p = \sum\limits_{i=0}^{d} p_i x^{d-i}$, where $p_i \in \mathcal{F}$ and $p_0 = 1$. Moreover, $d$ is even and $p_i = 0$ for $i$ odd as $p$ is inseparable. Using the facts of Definition 3.7, we write $p_i = p_{i,0} + p_{i,1} t_{j_n}$, for some $p_{i,0}, p_{i,1} \in \mathcal{F}^2(t_{j_1}, \ldots, t_{j_{n-1}})$. Again we take $p_{i,1} = \sum\limits_{J \in T_1,\, J_n = 0} t^J p_{i,1,J}^2$, where $p_{i,1,J} \in \mathcal{F}$. To prove statement (3), we need a matrix $\tilde{M}$ on $\overline{F}_p$ similar to $M$. Furthermore, this matrix $\tilde{M}$ should be invertible when $p$ is inseparable. To do so we define in $\overline{F}_p$

$$\tilde{P}_K = \overline{t^K(p_{0,1,K}^2 x^d + p_{2,1,K}^2 x^{d-2} + \ldots + p_{d,1,K}^2)}, \quad \text{for all } K \in T_1 \text{ such that } K_n = 0.$$

Since $t_{j_n} \in \mathcal{F}(p)^2(t_{j_1}, \ldots, t_{j_{n-1}})$, we have $\sum\limits_{K \in T_1,\, K_n = 0} \tilde{P}_K = \frac{\partial p}{\partial t_{j_n}} \neq 0$. Let $\tilde{M}$ be the $2^{n-1} \times 2^{n-1}$ matrix whose $(I, J)$-th element as $\tilde{P}_{I+J t_{j_n}}$ for all $I, J \in T_1$ with $I_n = J_n = 0$. Since, $\sum\limits_{K \in T_1,\, K_n = 0} \tilde{P}_K \neq 0$, it follows that $\tilde{M}$ is invertible, see [4, Lemma 2.6]. Now we prove that $S_p \subseteq \tilde{S}_p + \tilde{S}_p \wedge \overline{\frac{dp}{p}} + L_0$. Recall that $S_p$ is generated by the elements $\overline{\frac{h}{p^e} \frac{dt_{i_1}}{t_{i_1}}} \wedge \ldots \wedge \frac{dt_{i_m}}{t_{i_m}}$, where $h \in \mathcal{F}[x]$ and $e \geq 0, \{i_1, \ldots, i_m\} \subseteq \mathcal{I}$. Now if $e = 0$, then $h \frac{dt_{i_1}}{t_{i_1}} \wedge \ldots \wedge \frac{dt_{i_m}}{t_{i_m}} \in L_0$, which implies that $S_{p,0} \subseteq L_0$. Before we continue our induction on $e$, we do some calculations.

**Calculations:** Let $\frac{dt_{\tilde{I}}}{t_{\tilde{I}}}$ denote $\frac{dt_{i_1}}{t_{i_1}} \wedge \frac{dt_{i_2}}{t_{i_2}} \wedge \ldots \wedge \frac{dt_{i_{m-1}}}{t_{i_{m-1}}}$, where $\tilde{I} = \{i_1, \ldots, i_{m-1}\} \subseteq \mathcal{I} \setminus \{j_n\}$ such that $i_1 < i_2 < \ldots < i_{m-1}$. For any $e \geq 0$ and $L \in T_1$ such that $L_n = 0$, we have:

$$\overline{\frac{h_L}{p^e}\frac{dt_{\tilde{I}}}{t_{\tilde{I}}}} \wedge \frac{d(t^L p)}{t^L p} = \sum_{i=0}^{d} \overline{\frac{p_i x^{d-i} h_L}{p^{e+1}}\frac{dt_{\tilde{I}}}{t_{\tilde{I}}}} \wedge \frac{d(t^L p_i)}{t^L p_i} \quad (\text{since } p_i = 0 \text{ for } i \text{ odd}).$$

Putting $p_i = p_{i,0} + p_{i,1}t_{j_n}$ we get:

$$\overline{\frac{h_L\,\mathrm{d}t_I}{p^e\,t_I}} \wedge \overline{\frac{\mathrm{d}(t^L p)}{t^L p}} = \sum_{i=0}^{d} \overline{\frac{x^{d-i}h_L p_{i,0}\,\mathrm{d}t_I}{p^{e+1}\,t_I}} \wedge \overline{\frac{\mathrm{d}(t^L p_{i,0})}{t^L p_{i,0}}} + \sum_{i=0}^{d} \overline{\frac{x^{d-i}h_L p_{i,1}t_{j_n}\,\mathrm{d}t_I}{p^{e+1}\,t_I}} \wedge \overline{\frac{\mathrm{d}(t^L p_{i,1}t_{j_n})}{t^L p_{i,1}t_{j_n}}}$$

$$\equiv \sum_{i=0}^{d} \overline{\frac{p_{i,1}t_{j_n}x^{d-i}h_L\,\mathrm{d}t_I}{p^{e+1}\,t_I}} \wedge \overline{\frac{\mathrm{d}(t^L p_{i,1}t_{j_n})}{t^L p_{i,1}t_{j_n}}} \quad (\mathrm{mod}\ \tilde{S}_p)$$

$$\equiv \sum_{i=0}^{d} \sum_{K \in T_1,\ K_n=0} \overline{\frac{t^K p_{i,1,K}^2 t_{j_n} x^{d-i} h_L\,\mathrm{d}t_I}{p^{e+1}\,t_I}} \wedge \overline{\frac{\mathrm{d}(t^{K+L}t_{j_n})}{t^{K+L}t_{j_n}}} \quad (\text{using value of } p_{i,1})$$

$$\equiv \sum_{K \in T_1,\ K_n=0} \overline{\frac{\left(t^K \sum_{i=0}^{d} p_{i,1,K}^2 x^{d-i}\right) t_{j_n} h_L\,\mathrm{d}t_I}{p^{e+1}\,t_I}} \wedge \overline{\frac{\mathrm{d}(t^{K+L}t_{j_n})}{t^{K+L}t_{j_n}}}$$

$$\equiv \sum_{K \in T_1,\ K_n=0} \overline{\frac{\tilde{P}_K t_{j_n} h_L\,\mathrm{d}t_I}{p^{e+1}\,t_I}} \wedge \overline{\frac{\mathrm{d}(t^{K+L}t_{j_n})}{t^{K+L}t_{j_n}}} \quad (\text{by putting the value of } \tilde{P}_K).$$

Therefore, taking sum on $L \in T_1$ such that $L_n = 0$, we get:

$$\sum_{L \in T_1,\ L_n=0} \overline{\frac{h_L\,\mathrm{d}t_I}{p^e\,t_I}} \wedge \overline{\frac{\mathrm{d}(t^L p)}{t^L p}} \equiv \sum_{L \in T_1,\ L_n=0} \sum_{K \in T_1,\ K_n=0} \overline{\frac{\tilde{P}_K t_{j_n} h_L\,\mathrm{d}t_I}{p^{e+1}\,t_I}} \wedge \overline{\frac{\mathrm{d}(t^{L+K}t_{j_n})}{t^{L+K}t_{j_n}}} \quad (\mathrm{mod}\ \tilde{S}_p).$$

Taking $H = L + K$, we obtain:

$$\sum_{L \in T_1,\ L_n=0} \overline{\frac{h_L\,\mathrm{d}t_I}{p^e\,t_I}} \wedge \overline{\frac{\mathrm{d}(t^L p)}{t^L p}} \equiv \sum_{H \in T_1,\ H_n=0} \overline{\frac{\left(\sum_{L \in T_1,\ L_n=0} \tilde{P}_{H+L} t_{j_n} h_L\right)\,\mathrm{d}t_I}{p^{e+1}\,t_I}} \wedge \overline{\frac{\mathrm{d}(t^H t_{j_n})}{t^H t_{j_n}}} \quad (\mathrm{mod}\ \tilde{S}_p).$$

For any $g \in \mathcal{F}[x]_{<d}$ we will have to find $2^{n-1}$ polynomials $h_L \in \mathcal{F}[x]$, which satisfy the following equations in $\overline{F}_p$:

$$g = \sum_{\{L \in T_1\ |\ L_n=0\}} \tilde{P}_L t_{j_n} h_L \quad \text{and} \quad \sum_{\{L \in T_1\ |\ L_n=0\}} \tilde{P}_{H+L} t_{j_n} h_L = 0 \text{ for } H \neq 0 \in T_1.$$

Equivalently, the column matrix $(h_L)_{\{L \in T_1\ |\ L_n=0\}}$ satisfies the following relation:

$$\tilde{M}.(h_L)_{\{L \in T_1\ |\ L_n=0\}} = (g, 0, \ldots, 0)^T.$$

Since $\tilde{M}$ is invertible, there exist $h_L \in \mathcal{F}[x]_{<d}$, which satisfy the above equations in $\overline{F}_p$. Hence, for any $g \in \mathcal{F}[x]_{<d}$ and $I \in \mathcal{T}$ with $\tilde{I} \subseteq \mathcal{I} \setminus \{j_n\}$, we can find $h_L \in \mathcal{F}[x]_{<d}$, for all $L \in T_1$ with $L_n = 0$, such that

(7) $$\overline{\frac{g\,\mathrm{d}t_I}{p^{e+1}\,t_I}} \wedge \overline{\frac{\mathrm{d}t_{j_n}}{t_{j_n}}} \equiv \sum_{\{L \in T_1\ |\ L_n=0\}} \overline{\frac{h_L\,\mathrm{d}t_I}{p^e\,t_I}} \wedge \overline{\frac{\mathrm{d}(t^L p)}{(t^L p)}} \quad (\mathrm{mod}\ \tilde{S}_p + S_{p,e}).$$

This finishes the calculations.

Now we know that $S_{p,0} \subseteq L_0$. For induction let us suppose $S_{p,j} \subseteq \tilde{S}_p + \tilde{S}_p \wedge \overline{\frac{\mathrm{d}p}{p}} + L_0$ for $0 \leq j \leq e$.

Moreover, for any $g \in \mathcal{F}[x]$ we write $\frac{g}{p^{e+1}} = \sum_{i=0}^{e+1} \frac{g_i}{p^i}$, where $g_0 \in \mathcal{F}[x]$ and $g_i \in \mathcal{F}[x]_{<d}$ for $1 \leq i \leq e+1$.

Therefore,

$$\overline{\frac{g}{p^{e+1}}\frac{\mathrm{d}t_I}{t_I}} \wedge \frac{\mathrm{d}t_{j_n}}{t_{j_n}} = \overline{\frac{g_{e+1}}{p^{e+1}}\frac{\mathrm{d}t_I}{t_I}} \wedge \frac{\mathrm{d}t_{j_n}}{t_{j_n}} + \sum_{i=0}^{e} \overline{\frac{g_i}{p^i}\frac{\mathrm{d}t_I}{t_I}} \wedge \frac{\mathrm{d}t_{j_n}}{t_{j_n}} \in \tilde{S}_p + \tilde{S}_p \wedge \overline{\frac{\mathrm{d}p}{p}} + L_0 \quad \text{(using induction hypothesis and (7))}.$$

Thus, $S_{p,e+1} \subseteq \tilde{S}_p + \tilde{S}_p \wedge \overline{\frac{\mathrm{d}p}{p}} + L_0$. Consequently, $S_{p,e} \subseteq \tilde{S}_p + \tilde{S}_p \wedge \overline{\frac{\mathrm{d}p}{p}} + L_0$ for all $e \geq 0$ and this implies that

(8) $$S_p \subseteq \tilde{S}_p + \tilde{S}_p \wedge \overline{\frac{\mathrm{d}p}{p}} + L_0.$$

Since $\tilde{S}_p \wedge \overline{\frac{\mathrm{d}x}{x}} \subseteq \tilde{S}_p$, we have from (8):

(9) $$S_p \wedge \overline{\frac{\mathrm{d}x}{x}} \subseteq \tilde{S}_p + \tilde{S}_p \wedge \overline{\frac{\mathrm{d}p}{p}} + L_0 \wedge \overline{\frac{\mathrm{d}x}{x}}.$$

Finally by (8) and (9),

$$S_p + S_p \wedge \overline{\frac{\mathrm{d}x}{x}} + L_0 + L_0 \wedge \overline{\frac{\mathrm{d}x}{x}} \subseteq \tilde{S}_p + \tilde{S}_p \wedge \overline{\frac{\mathrm{d}p}{p}} + L_0 + L_0 \wedge \overline{\frac{\mathrm{d}x}{x}}.$$

For converse by definition, we have

(10) $$\tilde{S}_p \subseteq S_p + S_p \wedge \overline{\frac{\mathrm{d}x}{x}},$$

and from statement (1) we know that

(11) $$S_p \wedge \overline{\frac{\mathrm{d}p}{p}} \subseteq S_p + S_p \wedge \overline{\frac{\mathrm{d}x}{x}}.$$

Therefore,

(12) $$S_p \wedge \overline{\frac{\mathrm{d}p}{p}} \wedge \overline{\frac{\mathrm{d}x}{x}} \subseteq S_p \wedge \overline{\frac{\mathrm{d}x}{x}} \subseteq S_p \wedge \overline{\frac{\mathrm{d}x}{x}} + L_0 \wedge \overline{\frac{\mathrm{d}x}{x}}.$$

By (10), (11) and (12), we have :

(13) $$\tilde{S}_p \wedge \overline{\frac{\mathrm{d}p}{p}} \subseteq S_p + S_p \wedge \overline{\frac{\mathrm{d}x}{x}} + L_0 \wedge \overline{\frac{\mathrm{d}x}{x}}.$$

So using (10) and (13) we get $\tilde{S}_p + \tilde{S}_p \wedge \overline{\frac{\mathrm{d}p}{p}} + L_0 + L_0 \wedge \overline{\frac{\mathrm{d}x}{x}} \subseteq S_p + S_p \wedge \overline{\frac{\mathrm{d}x}{x}} + L_0 + L_0 \wedge \overline{\frac{\mathrm{d}x}{x}}$.

This gives the desired equality and finishes the proof of statement (3). □

## 4. Milnor's exact sequence

As in Section 3, we take $\mathcal{F}$ a field of characteristic 2 with fixed 2-basis $\{t_i \mid i \in \mathcal{I}\}$, where $\mathcal{I}$ is an ordered set. Let $F = \mathcal{F}(x)$ be the rational function field in one variable $x$ over $\mathcal{F}$. Let $p \in \mathcal{F}[x]$ be a monic irreducible polynomial of degree $d$ or $p = \frac{1}{x}$, and $v_p : F \to \mathbb{Z}$ the $p$-adic valuation. Let $F_p$ and $\overline{F}_p$ be the completion of $F$ and the residue field, respectively. Note that in Section 2, we work with a fixed 2-basis of the residue field. Therefore, now we find that particular 2-basis of $\overline{F}_p$ for all $p$.

By [4, Lemma 3.1], for $p$ separable or $p = \frac{1}{x}$, we take $\{t_i \mid i \in \mathcal{I}\}$ as a 2-basis of $\overline{F}_p$ and $\{p\} \cup \{t_i \mid i \in \mathcal{I}\}$ a 2-basis of $F_p$. For $p$ inseparable, Definition 3.7 suggests that $\{\overline{x}\} \cup \{t_i \mid i \in \mathcal{I} \setminus \{j_n\}\}$ and $\{p, x\} \cup \{t_i \mid i \in \mathcal{I} \setminus \{j_n\}\}$ will be the 2-bases of $\overline{F}_p$ and $F_p$, respectively.

From now on we will use the following notations for the set of products of the elements of 2-bases of $\mathcal{F}$, $\overline{F}_p$ and $F_p$.

**Definition 4.1.** *(1) Let $\{t^I \mid I \in T\}$ be the set of products of different elements of the fixed 2-basis of $\mathcal{F}$.*
*(2) Let $\{t^I \mid I \in \overline{T}_p\}$ be the set of products of distinct elements of the 2-basis of $\overline{F}_p$, which depends on the separability of the polynomial p as stated just before the Definition.*
*(3) Let $\{t^I \mid I \in T_p\}$ be the set of products of the distinct elements of the 2-basis of $F_p$, which is stated before the Definition.*

Now we explain how the elements of $W_1(H_2^{m+1}(F_p))$ look like for any monic irreducible $p \in \mathcal{F}[x]$ and $p = \frac{1}{x}$.

**Theorem 4.2.** *(1) Let $p \in \mathcal{F}[x]$ be monic and irreducible. Then, for each element $\varphi \in W_1(H_2^{m+1}(F_p))$ we have a unique decomposition as $\varphi = \psi \perp \varphi_2 \wedge \overline{\dfrac{dp}{p}}$, where*

$$\varphi_2 \in H_2^m(\mathcal{F}(p)) \quad \text{and} \quad \psi = \sum_{I \in (T_p)_m} \sum_{J \in T_p,\, J+I>I} \overline{t^J u_{I,J}^2 \dfrac{dt_I}{t_I}}$$

*such that $u_{I,J} = \sum_{l \geq 1} \dfrac{u_{l,I,J}}{p^l}$ with $u_{l,I,J} \in \mathcal{F}[x]_{<\deg p}$. Here $H_2^m(\mathcal{F}(p)) \wedge \overline{\dfrac{dp}{p}}$ is viewed as a subgroup of $H_2^{m+1}(F_p)$ and $u_{l,I,J}, \psi$ and $\varphi_2$ are unique for each $\varphi$. (Recall that $(T_p)_m = \{I \in T_p \mid |\tilde{I}| = m\}$, where $\tilde{I}$ is the support of I, see Section 2).*

*(2) For $p = \frac{1}{x}$, any element $\varphi \in W_1(H_2^{m+1}(F_{\frac{1}{x}}))$ has a unique decomposition as $\varphi = \psi + \varphi_2 \wedge \overline{\dfrac{dx}{x}}$, where*

$$\varphi_2 \in H_2^m(\mathcal{F}) \quad \text{and} \quad \psi = \sum_{I \in (T_{\frac{1}{x}})_m} \sum_{J \in T_{\frac{1}{x}},\, J+I>I} \overline{t^J u_{I,J}^2 \dfrac{dt_I}{t_I}} \quad \text{with } u_{I,J} \in x\mathcal{F}[x].$$

*Here $u_{I,J}, \psi, \varphi_2$ are unique for each $\varphi$. Moreover, we have $\sum_{J,\, J+I>I} t^J u_{I,J}^2 \in x\mathcal{F}[x]$.*

*Proof.* (1) Let $p$ be separable and $\alpha : \overline{F}_p \to F_p$ the unique field monomorphism, i.e., the Teichmüller lifting such that $\alpha(t_i) = t_i$ for all $i \in \mathcal{I}$. If $p$ is inseparable, then we take $\alpha : \overline{F}_p \to F_p$ such that $\alpha(t_i) = t_i$ for all $i \in \mathcal{I} \setminus \{j_n\}$ and $\alpha(\overline{x}) = x$. Here, we use the same 2-basis which was defined just before of Definition 4.1 for separable and inseparable cases. Now $\overline{F}_p \simeq \alpha(\overline{F}_p)$ and we can see $\alpha(\overline{F}_p)$ as the residue field of $F_p$ for $p$-adic valuation. Since, $\alpha(\overline{F}_p) \subseteq F_p$, it satisfies the conditions of Theorem 2.3, which implies that each element $\varphi \in H_2^{m+1}(F_p)$ has the following unique decomposition

$$\varphi = \varphi_1' + \psi' + \varphi_2' \wedge \overline{\dfrac{dp}{p}},$$

where $\varphi_1' \in H_2^{m+1}(\alpha(\overline{F}_p))$, $\varphi_2' \in H_2^m(\alpha(\overline{F}_p))$ and $\psi' = \sum_{I \in (T_p)_m} \sum_{\{J \in T_p \mid J+I>I\}} \overline{t^J r_{I,J}^2 \dfrac{dt_I}{t_I}}$ for $r_{I,J} \in p^{-1}\alpha(\overline{F}_p)[p^{-1}]$.

Since $W_1(H_2^{m+1}(F_p)) = \operatorname{coker}(H_2^{m+1}(\overline{F}_p) \xrightarrow{\alpha'} H_2^{m+1}(F_p))$, it follows that $W_1(H_2^{m+1}(F_p)) = H_2^{m+1}(F_p)/H_2^{m+1}(\alpha(\overline{F}_p))$, where $\alpha'$ is induced by $\alpha$. So the class of $\varphi$ in $W_1(H_2^{m+1}(F_p))$ is given as follows: $\varphi = \psi' + \varphi_2' \wedge \overline{\dfrac{dp}{p}}$.

**Step 1:** By definition, $\psi' = \sum_{I \in (T_p)_m} \sum_{\{J \mid J+I>I\}} \overline{t^J r_{I,J}^2 \frac{dt_I}{t_I}}$, where $r_{I,J} = \sum_{i=1}^{k} \frac{\alpha(\overline{f_i})}{p^i}$ for some $f_i \in \mathcal{F}[x]$ and $k \in \mathbb{N}_0$.

So in $F_p$, the element $t^J r_{I,J}^2$ can be written as $t^J r_{I,J}^2 = \sum_{i=1}^{l} \frac{\alpha(\overline{g_i})}{p^i}$ for some $g_i \in \mathcal{F}[x]$ and $l \in \mathbb{N}_0$. Moreover, $\alpha(\overline{g_i}) = \sum_{j \geq 0} g_{i,j} p^j$, where $g_{i,j} \in \mathcal{F}[x]_{<\deg p}$. Hence, for any $i \in \mathbb{N}_0$ and any $\overline{g_i} \in \mathcal{F}(p)$, we have:

$$\frac{\alpha(\overline{g_i})}{p^i} = \sum_{l \leq i} \frac{g_{i,l}}{p^l} \in F_p, \quad g_{i,l} \in \mathcal{F}[x]_{<\deg p}.$$

Therefore, it suffices to prove that, in $W_1(H_2^{m+1}(F_p))$, for any fixed $h \in \mathcal{F}[x]_{<\deg p}, I \in (T_p)_m$ and $j \in \mathbb{Z}$, we have

(14) $$\overline{\frac{h}{p^j} \frac{dt_I}{t_I}} = \sum_{I \in (T_p)_m} \sum_{\{J \in T_p \mid J+I>I\}} \overline{t^J u_{I,J}^2 \frac{dt_I}{t_I}} + \varphi_2 \wedge \overline{\frac{dp}{p}},$$

where $\varphi_2 \in H_2^m(\mathcal{F}(p))$ and $u_{I,J} = \sum_{l \geq 1} \frac{u_{l,I,J}}{p^l}$, with $u_{l,I,J} \in \mathcal{F}[x]_{<\deg p}$.

- If $j < 0$, then $v_p\left(\frac{h}{p^j}\right) > 0$ and thus $\frac{h}{p^j} \in \wp(F_p)$. Hence, $\overline{\frac{h}{p^j} \frac{dt_I}{t_I}} = 0 \in H_2^{m+1}(F_p)$.
- If $j = 0$, then $\overline{h \frac{dt_I}{t_I}} \in H_2^m(\mathcal{F}(p)) \wedge \overline{\frac{dp}{p}}$ or $H_2^{m+1}(\mathcal{F}(p))$, according as $I_\infty = 1$ or $0$. Moreover, for $I_\infty = 0$,

$\overline{h \frac{dt_I}{t_I}} = \overline{\alpha(\overline{h}) \frac{dt_I}{t_I}} \in \text{Im}(H_2^{m+1}(\mathcal{F}(p)) \xrightarrow{\alpha'} H_2^{m+1}(F_p))$. So $\overline{h \frac{dt_I}{t_I}} = 0 \in W_1(H_2^{m+1}(F_p))$.

- Let $j > 0$. Recall that, in Step 1 of the proof of Theorem 2.5, we proved that for any $c \in \overline{F}, k \geq 1$ and $I \in (T)_m$, $\overline{\frac{c}{\pi^k} \frac{dt_I}{t_I}} \in \mathcal{R}$, where $\pi$ is an uniformizer. By exactly the same way, we can prove that for each $h \in \mathcal{F}[x]_{<\deg p}, j \geq 1$ and $I \in (T_p)_m$, $\overline{\frac{h}{p^j} \frac{dt_I}{t_I}}$ can be decomposed as in (14). Finally, we obtain that any class $\varphi \in W_1(H_2^{m+1}(F_p))$ can be written as

$$\varphi = \psi' + \varphi_2' \wedge \overline{\frac{dp}{p}},$$

where $\varphi_2' \in H_2^m(\alpha(\mathcal{F}(p)))$, $\psi' = \sum_{I \in (T_p)_m} \sum_{\{J \mid J+I>I\}} \overline{t^J u_{I,J} \frac{dt_I}{t_I}} + \varphi_2 \wedge \overline{\frac{dp}{p}}$, for $u_{I,J} = \sum_{l \geq 1} \frac{u_{l,I,J}}{p^l}, u_{l,I,J} \in \mathcal{F}[x]_{<\deg p}$ and $\varphi_2 \in H_2^m(\mathcal{F}(p))$.

**Step 2:** Let us write $\varphi_2' = \sum_{I \in (\overline{T}_p)_{m-1}} \overline{f_I \frac{dt_I}{t_I}}$ for some $f_I \in \alpha(\overline{F}_p)$. So $f_I = f_{I,1} + f_{I,2}$, where $f_{I,1} \in \mathcal{F}[x]_{<\deg(p)}$ and $v_p(f_{I,2}) > 0$. Hence, $f_{I,2} \in \wp(F_p)$ and therefore $\varphi_2' = \sum_{I \in (\overline{T}_p)_{m-1}} \overline{f_{I,1} \frac{dt_I}{t_I}} \in H_2^m(\mathcal{F}(p))$.

Taking Step 1 and Step 2 together, we have $\varphi = \psi + \varphi_2 \wedge \overline{\frac{dp}{p}} \in W_1(H_2^{m+1}(F_p))$, where $\varphi_2 \in H_2^m(\mathcal{F}(p))$ and $\psi = \sum_{I \in (T_p)_m} \sum_{\{J \in T_p \mid J+I>I\}} \overline{t^J u_{I,J}^2 \frac{dt_I}{t_I}}$, where $u_{I,J} = \sum_{l \geq 1} \frac{u_{l,I,J}}{p^l}$ with $u_{l,I,J} \in \mathcal{F}[x]_{<\deg p}$.

For the uniqueness, we repeat the same technique of the Uniqueness part of Theorem 2.5.

(2) For $p = \frac{1}{x}$, we have $\overline{F}_p = \mathcal{F} \subseteq F_p$, and thus we are in the conditions of Theorem 2.5. Now the result clearly comes from this theorem. □

The next lemma describes the group $H_2^{m+1}(F)$ by generators and relations.

**Lemma 4.3.** *Let $F$ be a field of characteristic $p > 0$. Then, there exists a natural group homomorphism $\theta : F \otimes \underbrace{F^* \otimes \ldots \otimes F^*}_{m} \to H_p^{m+1}(F)$ defined by: $\theta(a \otimes b_1 \otimes \ldots \otimes b_m) = a\overline{\frac{db_1}{b_1} \wedge \ldots \wedge \frac{db_m}{b_m}}$. Moreover, $\theta$ is surjective and its kernel is generated by the following elements:*

(1) $(a^p - a) \otimes b_1 \otimes b_2 \otimes \ldots \otimes b_m$,
(2) $a \otimes a \otimes b_2 \otimes \ldots \otimes b_m$,
(3) $a \otimes b_1 \otimes b_2 \otimes \ldots \otimes b_m$, *where $b_i = b_j$ for some $i \neq j$.*

*Proof.* See [10, Page 22]. □

Using Theorem 4.2, we introduce in the next lemma a homomorphism $\tau_p : W_1(H_2^{m+1}(F_p)) \to L_d/L_{d-1}$, which is called the Milnor's splitting. For $p \in \mathcal{F}[x]$ a monic irreducible polynomial of degree $d \geq 1$, and $g \in \mathcal{F}[x]$, let $\overline{g}$ denote the unique polynomial of degree $\leq d - 1$ such that $\overline{g} \equiv g \pmod{p}$.

**Lemma 4.4.** *Let $p$ be a monic irreducible polynomial of degree $d \geq 1$. By Theorem 4.2, any element $\varphi \in W_1(H_2^{m+1}(F_p))$ can be written as*

$$\varphi = \psi + \varphi_2 \wedge \overline{\frac{dp}{p}}, \quad \text{where } \psi = \sum_{I \in (T_p)_m} \sum_{\{J \mid J+I>I\}} \overline{t^J r_{I,J}^2 \frac{dt_I}{t_I}},$$

*with unique $r_{I,J} \in p^{-1}\mathcal{F}[x]_{<d}[p^{-1}]$ and $\varphi_2$ belongs to $H_2^m(\mathcal{F}(p))$, considered as a subgroup of $H_2^m(F_p)$ by Corollary 6.9. Moreover, $\varphi_2$ is unique over $\mathcal{F}(p)$. Let $\varphi_2 = \sum_i \overline{r_i \frac{ds_{i_1}}{s_{i_1}} \wedge \ldots \wedge \frac{ds_{i_{m-1}}}{s_{i_{m-1}}}}$, where $r_i, s_{i_j} \in \mathcal{F}[x]_{<d}$. Then, we have a homomorphism $\tau_p : W_1(H_2^{m+1}(F_p)) \to L_d/L_{d-1}$, given by:*

$$\tau_p(\varphi) = \psi \perp \overline{\left(\sum_i r_i \frac{ds_{i_1}}{s_{i_1}} \wedge \ldots \wedge \frac{ds_{i_{m-1}}}{s_{i_{m-1}}}\right) \wedge \frac{dp}{p}} + L_{d-1} \in L_d/L_{d-1}.$$

*Proof.* We will have to prove well definedness of $\tau_p$ only. It suffices to define a map $\eta_p : H_2^m(\overline{F}_p) \to L_d/L_{d-1}$ giving as follows:

$$\eta_p\left(\sum_i \overline{\overline{r_i} \frac{d\overline{s_{i_1}}}{\overline{s_{i_1}}} \wedge \ldots \wedge \frac{d\overline{s_{i_{m-1}}}}{\overline{s_{i_{m-1}}}}}\right) = \overline{\left(\sum_i r_i \frac{ds_{i_1}}{s_{i_1}} \wedge \ldots \wedge \frac{ds_{i_{m-1}}}{s_{i_{m-1}}}\right) \wedge \frac{dp}{p}} + L_{d-1},$$

where $r_i, s_{i_j} \in \mathcal{F}[x]_{<d}$. By Lemma 4.3, we know that $H_2^m(\overline{F}_p) \simeq \overline{F}_p \otimes \underbrace{\overline{F}_p^* \otimes \ldots \otimes \overline{F}_p^*}_{m-1}/J$, where $J$ is the subgroup generated by the following elements:

(1) $(\overline{r}^2 - \overline{r}) \otimes \overline{s_1} \otimes \ldots \otimes \overline{s_{m-1}}$, where $s_i \neq 0, r \in \mathcal{F}[x]_{<d}$,

(2) $\overline{r} \otimes \overline{r} \otimes \overline{s_2} \otimes \ldots \otimes \overline{s_{m-1}}$,

(3) $\overline{r} \otimes \overline{s_1} \otimes \ldots \otimes \overline{s_{m-1}}$, where $\overline{s_i} = \overline{s_j}$, for some $i \neq j$.

Let us consider a map $\zeta_p : \overline{F}_p \times \underbrace{\overline{F}_p^* \times \ldots \times \overline{F}_p^*}_{m-1} \to L_d/L_{d-1}$, defined by:

$$\zeta_p(\overline{r}, \overline{s_1}, \ldots, \overline{s_{m-1}}) = \overline{r \frac{ds_1}{s_1} \wedge \ldots \wedge \frac{ds_{m-1}}{s_{m-1}} \wedge \frac{dp}{p}} + L_{d-1}.$$

Clearly, $\zeta_p$ is linear with respect to the first quadrant. Now we prove the linearity for the other quadrants. It suffices to do it for the second quadrant. Let $t_1 \in \mathcal{F}[x]_{<d}$ and $\overline{s_1 t_1} = \overline{s}$ for some $s \in \mathcal{F}[x]_{<d}$. Therefore, $s_1 t_1 = s + pt$ for some $t \in \mathcal{F}[x]_{<d}$. Let $p = \sum_{j=0}^{d} p_j x^{d-j}$ for $p_j \in \mathcal{F}$, and put $A = \frac{ds_2}{s_2} \wedge \ldots \wedge \frac{ds_{m-1}}{s_{m-1}} \wedge \frac{dp}{p}$. Then, we have:

$$\overline{\frac{rp}{s} dt \wedge A} = \overline{\frac{r}{s} dt \wedge \frac{ds_2}{s_2} \wedge \ldots \wedge \frac{ds_{m-1}}{s_{m-1}} \wedge dp} = \sum_{j=0}^{d} \overline{\frac{p_j x^{d-j} rt}{s} \frac{dt}{t} \wedge \frac{ds_2}{s_2} \wedge \ldots \wedge \frac{ds_{m-1}}{s_{m-1}} \wedge \frac{d(p_j x^{d-j})}{p_j x^{d-j}}}$$

$$= \sum_{d-j \text{ even}} \overline{\frac{p_j x^{d-j} rt}{s} \frac{dt}{t} \wedge \frac{ds_2}{s_2} \wedge \ldots \wedge \frac{ds_{m-1}}{s_{m-1}} \wedge \frac{dp_j}{p_j}} + \sum_{d-j \text{ odd}} \overline{\frac{p_j x^{d-j} rt}{s} \frac{dt}{t} \wedge \frac{ds_2}{s_2} \wedge \ldots \wedge \frac{ds_{m-1}}{s_{m-1}} \wedge \frac{d(p_j x)}{p_j x}}.$$

Clearly, $\overline{\frac{rp}{s} dt \wedge A} \in L_{d-1}$ because $\deg t, \deg s, \deg s_i \le d-1$, $p_j \in \mathcal{F}$ and $\frac{dx}{x}$ exists in second sum where the separable part will be a multiple of $x$ in $d = 1$ case. Similarly we can show that

(15) $$\overline{\frac{r(s_1 t_1 - s)}{s} \frac{d(s_1 t_1)}{s_1 t_1}} \wedge A = \overline{\frac{rpt}{s} \frac{d(s_1 t_1)}{s_1 t_1}} \wedge A \in L_{d-1},$$

where $s_i, t_1, t, s \in \mathcal{F}[x]_{<d}$. Now we compute:

$$\zeta_p\left(\overline{r}, \overline{s_1 t_1}, \overline{s_2}, \ldots, \overline{s_{m-1}}\right) = \overline{r \frac{ds}{s} \wedge \frac{ds_2}{s_2} \wedge \ldots \wedge \frac{ds_{m-1}}{s_{m-1}} \wedge \frac{dp}{p}} + L_{d-1} \quad (\text{since } \overline{s_1 t_1} = \overline{s})$$

$$= \overline{r \frac{d(s_1 t_1 - pt)}{s} \wedge \frac{ds_2}{s_2} \wedge \ldots \wedge \frac{ds_{m-1}}{s_{m-1}} \wedge \frac{dp}{p}} + L_{d-1}$$

$$= \overline{\frac{r}{s}[d(s_1 t_1) + d(pt)] \wedge \frac{ds_2}{s_2} \wedge \ldots \wedge \frac{ds_{m-1}}{s_{m-1}} \wedge \frac{dp}{p}} + L_{d-1}.$$

Using $dp \wedge dp = 0$, we get:

$$\zeta_p\left(\overline{r}, \overline{s_1 t_1}, \overline{s_2}, \ldots, \overline{s_{m-1}}\right) = \overline{r \left[\frac{d(s_1 t_1)}{s} + \frac{pdt}{s}\right] \wedge A} + L_{d-1}$$

$$= \overline{\frac{rs_1 t_1}{s} \frac{d(s_1 t_1)}{s_1 t_1} \wedge A} + L_{d-1} \quad \left(\text{since } \overline{\frac{rp}{s} dt \wedge A} \in L_{d-1}\right)$$

$$= \overline{r \frac{d(s_1 t_1)}{s_1 t_1} \wedge A} + L_{d-1} \quad (\text{ by (15)})$$

$$= \overline{r \left[\frac{ds_1}{s_1} + \frac{dt_1}{t_1}\right] \wedge A} + L_{d-1}$$

$$= \zeta_p\left(\overline{r}, \overline{s_1}, \ldots, \overline{s_{m-1}}\right) + \zeta_p\left(\overline{r}, \overline{t_1}, \ldots, \overline{s_{m-1}}\right).$$

Therefore, we have a homomorphism $\zeta'_p : \overline{F}_p \otimes \underbrace{\overline{F}^*_p \otimes \ldots \otimes \overline{F}^*_p}_{m-1} \to L_d/L_{d-1}$, defined on the generators:

$$\zeta'_p(\overline{r} \otimes \overline{s_1} \otimes \ldots \otimes \overline{s_{m-1}}) = \overline{r\frac{ds_1}{s_1} \wedge \ldots \wedge \frac{ds_{m-1}}{s_{m-1}} \wedge \frac{dp}{p}} + L_{d-1}.$$

Now we will show that $\zeta'_p(J) = 0$, where $J$ is the subgroup of $\overline{F}_p \otimes \underbrace{\overline{F}^*_p \otimes \ldots \otimes \overline{F}^*_p}_{m-1}$ defined previously. For the second type of generators of $J$, we have :

$$\zeta'_p(\overline{r} \otimes \overline{r} \otimes \overline{s_2} \otimes \ldots \otimes \overline{s_{m-1}}) = \overline{r\frac{dr}{r} \wedge \ldots \wedge \frac{ds_{m-1}}{s_{m-1}} \wedge \frac{dp}{p}} + L_{d-1} = 0.$$

For the third type of generators, $\overline{s_i} = \overline{s_j}$ and $\deg s_i, \deg s_j \leq d-1$ implies that $s_i = s_j$. Since $ds_i \wedge ds_i = 0$, therefore $\zeta_p$ sends these generators to 0. Now for the first type of generators, let $\overline{r}^2 = \overline{f}$, so $r^2 = f + ph$ for $f, h \in \mathcal{F}[x]_{<d}$. Then, we have:

$$\zeta'_p\left((\overline{r}^2 - \overline{r}) \otimes \overline{s_1} \otimes \overline{s_2} \otimes \ldots \otimes \overline{s_{m-1}}\right) = \overline{(f-r)\frac{ds_1}{s_1} \wedge \ldots \wedge \frac{ds_{m-1}}{s_{m-1}} \wedge \frac{dp}{p}} + L_{d-1}$$

$$= \overline{(r^2 - r - ph)\frac{ds_1}{s_1} \wedge \ldots \wedge \frac{ds_{m-1}}{s_{m-1}} \wedge \frac{dp}{p}} + L_{d-1} \quad \text{(since } f = r^2 - ph)$$

$$= \overline{h\frac{ds_1}{s_1} \wedge \ldots \wedge \frac{ds_{m-1}}{s_{m-1}} \wedge dp} + L_{d-1} \quad \text{( as } r^2 - r \in \wp(F) \text{ )}$$

$$= \overline{\frac{ds_1}{s_1} \wedge \ldots \wedge \frac{ds_{m-1}}{s_{m-1}} \wedge (hdp + d(ph))} + L_{d-1}$$

$$= \overline{\frac{ds_1}{s_1} \wedge \ldots \wedge \frac{ds_{m-1}}{s_{m-1}} \wedge (pdh)} + L_{d-1}$$

$$= \overline{\frac{p}{h}\frac{ds_1}{s_1} \wedge \ldots \wedge \frac{ds_{m-1}}{s_{m-1}} \wedge \frac{dh}{h}} + L_{d-1} = 0 \quad \text{(since } \deg s_i, \deg h \leq d-1).$$

Therefore, $\zeta'_p$ sends $J$ to 0 and we extend the map $\zeta'_p$ on $\overline{F}_p \otimes \underbrace{\overline{F}^*_p \otimes \ldots \otimes \overline{F}^*_p}_{m-1}/J$. Since $H_2^m(\overline{F}_p) \simeq \overline{F}_p \otimes \underbrace{\overline{F}^*_p \otimes \ldots \otimes \overline{F}^*_p}_{m-1}/J$, it follows that the map $\eta_p$ is well defined.

By Theorem 4.2, each $\varphi \in W_1(H_2^{m+1}(F_p))$ can be written as $\varphi = \psi + \varphi_2 \wedge \frac{dp}{p}$, where

$$\psi = \sum_{I \in (T_p)_m} \sum_{\{J \in T_p \ | \ J+I>I\}} \overline{t^J r_{I,J}^2 \frac{dt_I}{t_I}}, \ r_{I,J} \in p^{-1}\mathcal{F}[x]_{<d}[p^{-1}] \text{ and } \varphi_2 \in H_2^m(\mathcal{F}(p)),$$

with unique $r_{I,J}$ and $\varphi_2$. Moreover, $\varphi_2$ is unique over $\mathcal{F}(p)$ and $F_p$ both. Since $\eta_p$ is well defined,

$$\tau_p(\varphi) = \psi + \eta_p(\varphi_2)$$

is also well defined. Moreover, $\tau_p$ is additive and thus it is a homomorphism. □

The following lemma is crucial to prove the surjectivity of $\bigoplus_p \tau_p$.

**Lemma 4.5.** *Let $p = \sum_{j=0}^{d} p_j x^{d-j} \in \mathcal{F}[x]$ be a monic irreducible polynomial of degree d.*
*(1) If p is separable, then $S_p + S_p \wedge \overline{\frac{dp}{p}} + L_{d-1} \subseteq \mathrm{Im}(\tau_p)$.*
*(2) If p is inseparable, then $\tilde{S}_p + \tilde{S}_p \wedge \overline{\frac{dp}{p}} + L_{d-1} \subseteq \mathrm{Im}(\tau_p)$.*

*Proof.* (1) Suppose that $p$ is separable. Recall that $S_p$ is generated by the forms $\overline{\frac{h}{p^e} \frac{dc_1}{c_1} \wedge \ldots \wedge \frac{dc_m}{c_m}}$, where $h \in \mathcal{F}[x], c_i \in \mathcal{F}^*$ and $e \geq 0$. Similarly $S_p \wedge \overline{\frac{dp}{p}}$ is generated by the elements $\overline{\frac{g}{p^e} \frac{da_1}{a_1} \wedge \ldots \wedge \frac{da_{m-1}}{a_{m-1}} \wedge \frac{dp}{p}}$, where $g \in \mathcal{F}[x], a_i \in \mathcal{F}^*$ and $e \geq 0$. Moreover, $S_{p,e}$ is the subgroup of $S_p$ generated by $\overline{\frac{h}{p^i} \frac{dc_1}{c_1} \wedge \ldots \wedge \frac{dc_m}{c_m}}$, for $h \in \mathcal{F}[x], c_i \in \mathcal{F}^*$ and $i \leq e \in \mathbb{N}$. Similarly we have the subgroup $S_{p,e} \wedge \overline{\frac{dp}{p}} \subseteq S_p \wedge \overline{\frac{dp}{p}}$.

To get (1) it suffices to prove that $S_{p,e} + S_{p,e} \wedge \overline{\frac{dp}{p}} + L_{d-1} \subseteq \mathrm{Im}(\tau_p)$ for all $e \in \mathbb{N}$. To this end, we will proceed by induction on $e$. If $e = 0$, then $h\overline{\frac{dc_1}{c_1} \wedge \ldots \wedge \frac{dc_m}{c_m}} \in L_{d-1}$. Hence, $S_{p,0} + L_{d-1} = L_{d-1} \in \mathrm{Im}(\tau_p)$.

For any $g \in \mathcal{F}[x]$, we can write $g = \sum_{i=0}^{l} g_i p^i$ for suitable $g_i \in \mathcal{F}[x]_{<d}$. Now $g_0 \overline{\frac{da_1}{a_1} \wedge \ldots \wedge \frac{da_{m-1}}{a_{m-1}} \wedge \frac{dp}{p}} + L_{d-1} \in \mathrm{Im}(\tau_p)$, and for $i \geq 1$ we have:

$$\overline{g_i p^i \frac{da_1}{a_1} \wedge \ldots \wedge \frac{da_{m-1}}{a_{m-1}} \wedge \frac{dp}{p}} = \overline{p^{i-1} g_i \frac{da_1}{a_1} \wedge \ldots \wedge \frac{da_{m-1}}{a_{m-1}} \wedge dp}$$

$$= \sum_{j=0}^{d} \overline{p^{i-1} g_i p_j x^{d-j} \frac{da_1}{a_1} \wedge \ldots \wedge \frac{da_{m-1}}{a_{m-1}} \wedge \frac{d(p_j x^{d-j})}{p_j x^{d-j}}}$$

$$= \sum_{\{j \mid d-j=\text{even}\}} \overline{p^{i-1} g_i p_j x^{d-j} \frac{da_1}{a_1} \wedge \ldots \wedge \frac{da_{m-1}}{a_{m-1}} \wedge \frac{dp_j}{p_j}} +$$

$$\sum_{\{j \mid d-j=\text{odd}\}} \overline{p^{i-1} g_i p_j x^{d-j} \frac{da_1}{a_1} \wedge \ldots \wedge \frac{da_{m-1}}{a_{m-1}} \wedge \frac{d(p_j x)}{p_j x}} \in L_{d-1}.$$

Hence, $S_{p,0} \wedge \overline{\frac{dp}{p}} + L_{d-1} \subseteq \mathrm{Im}(\tau_p)$ and thus $S_{p,0} + S_{p,0} \wedge \overline{\frac{dp}{p}} + L_{d-1} \subseteq \mathrm{Im}(\tau_p)$. To complete the induction, let us assume $S_{p,e} + S_{p,e} \wedge \overline{\frac{dp}{p}} + L_{d-1} \subseteq \mathrm{Im}(\tau_p)$ for some $e \geq 0$. Now we prove the inclusion for $e+1$.

Let $h, g \in \mathcal{F}[x]$, and take $\overline{h}, \overline{g}$ their classes in $\mathcal{F}(p)$, respectively. Using the 2-basis of $\mathcal{F}(p)$, we have $\overline{h} = \sum_{J \in \overline{T}_p} t^J \overline{h_J}^2$, where $h_J \in \mathcal{F}[x]_{<d}$. Therefore, $h = \sum_{J \in \overline{T}_p} t^J h_J^2 + pk$ for some $k \in \mathcal{F}[x]$. Similarly, we find some $g_J \in \mathcal{F}[x]_{<d}$ and $k' \in \mathcal{F}[x]$ such that $g = \sum_{J \in \overline{T}_p} t^J g_J^2 + pk'$. So putting the value of $h, g$ and using the fact that the scalars $c_i, a_j$ can be written using the 2-basis $T = \overline{T}_p$, we have modulo $S_{p,e} + S_{p,e} \wedge \overline{\frac{\mathrm{d}p}{p}}$:

$$\overline{\frac{h}{p^{e+1}} \frac{\mathrm{d}c_1}{c_1} \wedge \ldots \wedge \frac{\mathrm{d}c_m}{c_m}} + \overline{\frac{g}{p^{e+1}} \frac{\mathrm{d}a_1}{a_1} \wedge \ldots \wedge \frac{\mathrm{d}a_{m-1}}{a_{m-1}} \wedge \frac{\mathrm{d}p}{p}} \equiv \sum_{I \in (T_p)_m} \sum_{J \in \overline{T}_p} \overline{t^J \frac{g_{I,J}^2}{p^{e+1}} \frac{\mathrm{d}t_I}{t_I}}, \; g_{I,J} \in \mathcal{F}[x]_{<d}.$$

By induction hypothesis, $S_{p,e} + S_{p,e} \wedge \overline{\frac{\mathrm{d}p}{p}} + L_{d-1} \subseteq \mathrm{Im}(\tau_p)$. So we will have to show that $\sum_{I \in (T_p)_m} \sum_{J \in \overline{T}_p} \overline{t^J \frac{g_{I,J}^2}{p^{e+1}} \frac{\mathrm{d}t_I}{t_I}} +$
$L_{d-1} \in \mathrm{Im}(\tau_p)$, where $g_{I,J} \in \mathcal{F}[x]_{<d}$. We distinguish between two cases:

**Case 1.** Suppose $e + 1$ odd.

$$\sum_{I \in (T_p)_m} \sum_{J \in \overline{T}_p} \overline{t^J \frac{g_{I,J}^2}{p^{e+1}} \frac{\mathrm{d}t_I}{t_I}} = \sum_{I \in (T_p)_m} \sum_{J \in \overline{T}_p} \overline{t^J p \frac{g_{I,J}^2}{p^{e+2}} \frac{\mathrm{d}t_I}{t_I}} = \sum_{I \in (T_p)_m, \{J' \in T_p \mid J' \neq 0\}} \overline{t^{J'} \left(\frac{g_{I,J}}{p^{\frac{e+2}{2}}}\right)^2 \frac{\mathrm{d}t_I}{t_I}} \quad (\text{where } t^{J'} = t^J p)$$

$$= \sum_{I' \in (T_p)_m} \sum_{\{J' \in T_p \mid J'+I' > I'\}} \overline{t^{J'} \left(\frac{g_{I,J}}{p^{\frac{e+2}{2}}}\right)^2 \frac{\mathrm{d}t_{I'}}{t_{I'}}} \quad (\text{by Remark 2.4}).$$

It is clear by definition that this last term, modulo $L_{d-1}$, belongs to $\mathrm{Im}(\tau_p)$.

**Case 2.** Suppose $e + 1$ even.

$$\sum_{I \in (T_p)_m} \sum_{J \in \overline{T}_p} \overline{t^J \frac{g_{I,J}^2}{p^{e+1}} \frac{\mathrm{d}t_I}{t_I}} = \sum_{I \in (T_p)_m} \sum_{J \neq 0} \overline{t^J \frac{g_{I,J}^2}{p^{e+1}} \frac{\mathrm{d}t_I}{t_I}} + \sum_{I \in (T_p)_m} \overline{\frac{g_{I,0}^2}{p^{e+1}} \frac{\mathrm{d}t_I}{t_I}}$$

$$= \sum_{I \in (T_p)_m} \sum_{J \neq 0} \overline{t^J \left(\frac{g_{I,J}}{p^{\frac{e+1}{2}}}\right)^2 \frac{\mathrm{d}t_I}{t_I}} + \sum_{I \in (T_p)_m} \overline{\frac{g_{I,0}}{p^{\frac{e+1}{2}}} \frac{\mathrm{d}t_I}{t_I}}$$

$$\equiv \sum_{I \in (T_p)_m} \sum_{J \neq 0} \overline{t^J \left(\frac{g_{I,J}}{p^{\frac{e+1}{2}}}\right)^2 \frac{\mathrm{d}t_I}{t_I}} \quad (\mathrm{mod}\; S_{p,e} + S_{p,e} \wedge \overline{\frac{\mathrm{d}p}{p}}) \quad \left(\text{since } \frac{e+1}{2} \leq e\right)$$

$$= \sum_{I' \in (T_p)_m} \sum_{\{J \in T_p \mid J+I' > I'\}} \overline{t^J \left(\frac{g_{I,J}}{p^{\frac{e+1}{2}}}\right)^2 \frac{\mathrm{d}t_{I'}}{t_{I'}}} \quad (\text{by Remark 2.4}).$$

This last term, modulo $L_{d-1}$, belongs to $\mathrm{Im}(\tau_p)$, and by induction hypothesis $S_{p,e} + S_{p,e} \wedge \overline{\frac{\mathrm{d}p}{p}} + L_{d-1} \subseteq \mathrm{Im}(\tau_p)$.

Therefore, $\sum_{I \in (T_p)_m} \sum_{J \in \overline{T}_p} \overline{t^J \frac{g_{I,J}^2}{p^{e+1}} \frac{\mathrm{d}t_I}{t_I}} + L_{d-1} \subseteq \mathrm{Im}(\tau_p)$, and consequently $S_{p,e+1} + S_{p,e+1} \wedge \overline{\frac{\mathrm{d}p}{p}} + L_{d-1} \subseteq \mathrm{Im}(\tau_p)$.

Thus, $S_p + S_p \wedge \overline{\frac{\mathrm{d}p}{p}} + L_{d-1} \subseteq \mathrm{Im}(\tau_p)$ when $p$ is separable.

(2) Suppose that $p$ is inseparable. Recall that $\tilde{S}_p$ is generated by the elements $\overline{\frac{h}{p^e} \frac{\mathrm{d}t_I}{t_I}}$, where $h \in \mathcal{F}[x]$, $\frac{\mathrm{d}t_I}{t_I} = \frac{\mathrm{d}t_{i_1}}{t_{i_1}} \wedge \ldots \wedge \frac{\mathrm{d}t_{i_m}}{t_{i_m}}$ such that $\{t_{i_1}, \ldots, t_{i_m}\} \subseteq \{x\} \cup \{t_i \mid i \in \mathcal{I} \setminus \{j_n\}\}$ and $e \geq 0$. Let $\tilde{S}_{p,e}$ be the subgroup of $\tilde{S}_p$ generated by $\overline{\dfrac{h}{p^i}\dfrac{dt_I}{t_I}}$, where $h \in \mathcal{F}[x]$ and $i \le e$. Moreover, $\tilde{S}_p \wedge \overline{\dfrac{dp}{p}}$ is generated by the elements $\overline{\dfrac{g}{p^e}\dfrac{dt_L}{t_L}}$, where $g \in \mathcal{F}[x]$, $e \ge 0$ and $\dfrac{dt_L}{t_L} = \dfrac{dt_{l_1}}{t_{l_1}} \wedge \ldots \wedge \dfrac{dt_{l_m}}{t_{l_m}}$, such that $\{t_{l_1}, \ldots, t_{l_{m-1}}\} \subseteq \{x\} \cup \{t_i \mid i \in \mathcal{I} \setminus \{j_n\}\}$ and $t_{l_m} = p$. Similarly we can define $\tilde{S}_{p,e} \wedge \overline{\dfrac{dp}{p}}$ for all $e \ge 0$. It is sufficient to prove that $\tilde{S}_{p,e} + \tilde{S}_{p,e} \wedge \overline{\dfrac{dp}{p}} + L_{d-1} \subseteq \mathrm{Im}(\tau_p)$ for all $e \in \mathbb{N}$. We will prove this by induction on $e$.

Let us consider $e = 0$. Since $p$ is inseparable, $\deg p = d \ge 2$. Therefore, $\tilde{S}_{p,0} \subseteq L_{d-1}$ since its generators belong to $L_{d-1}$. We will prove that $\tilde{S}_{p,0} \wedge \overline{\dfrac{dp}{p}} + L_{d-1} \subseteq \mathrm{Im}(\tau_p)$. In fact, let $g = \sum\limits_{i=0}^{l} g_i p^i \in \mathcal{F}[x]$, where $g_i \in \mathcal{F}[x]_{<d}$. It is clear that $\overline{g_0 \dfrac{dt_L}{t_L}} + L_{d-1} \in \mathrm{Im}(\tau_p)$. Using $p = \sum\limits_{j=0}^{d} p_j x^{d-j}$, for all $i \ge 1$, as in $p$ separable case we have:

$$\overline{g_i p^i \dfrac{dt_L}{t_L}} = \sum_{\{j \mid d-j=\text{even}\}} \overline{g_i p^{i-1} p_j x^{d-j} \dfrac{dt_{l_1}}{t_{l_1}} \wedge \ldots \wedge \dfrac{dt_{l_{m-1}}}{t_{l_{m-1}}} \wedge \dfrac{dp_j}{p_j}} +$$
$$\sum_{\{j \mid d-j=\text{odd}\}} \overline{g_i p^{i-1} p_j x^{d-j} \dfrac{dt_{l_1}}{t_{l_1}} \wedge \ldots \wedge \dfrac{dt_{l_{m-1}}}{t_{l_{m-1}}} \wedge \dfrac{dp_j x}{p_j x}} \in L_{d-1}.$$

Hence, $\tilde{S}_{p,0} \wedge \overline{\dfrac{dp}{p}} + L_{d-1} \subseteq \mathrm{Im}(\tau_p)$. Therefore, $\tilde{S}_{p,0} + \tilde{S}_{p,0} \wedge \overline{\dfrac{dp}{p}} + L_{d-1} \subseteq \mathrm{Im}(\tau_p)$.

For induction, suppose that $\tilde{S}_{p,e} + \tilde{S}_{p,e} \wedge \overline{\dfrac{dp}{p}} + L_{d-1} \subseteq \mathrm{Im}(\tau_p)$ for some $e \ge 0$. Let $h, g \in \mathcal{F}[x]$ and $\overline{h}, \overline{g}$ their classes in $\mathcal{F}(p)$, respectively. As in the separable case, here also we have:
$$h = \sum_{J \in \overline{T}_p} t^J h_J^2 + pk, \quad \text{for some } h_J \in \mathcal{F}[x]_{<d}, k \in \mathcal{F}[x],$$
and
$$g = \sum_{J \in \overline{T}_p} t^J g_J^2 + pk', \quad \text{for some } g_J \in \mathcal{F}[x]_{<d}, k' \in \mathcal{F}[x].$$

Using these expressions, we get:
$$\overline{\dfrac{h}{p^{e+1}}\dfrac{dt_I}{t_I}} + \overline{\dfrac{g}{p^{e+1}}\dfrac{dt_L}{t_L}} + L_{d-1} \equiv \sum_{I \in (T_p)_m} \sum_{J \in \overline{T}_P} \overline{\dfrac{t^J g_{I,J}^2}{p^{e+1}}\dfrac{dt_I}{t_I}} + L_{d-1} \quad (\mathrm{mod}\ \tilde{S}_{p,e} + \tilde{S}_{p,e} \wedge \overline{\dfrac{dp}{p}} + L_{d-1}),$$
where $g_{I,J} \in \mathcal{F}[x]_{<d}$.

By induction hypothesis, $\tilde{S}_{p,e} + \tilde{S}_{p,e} \wedge \overline{\dfrac{dp}{p}} + L_{d-1} \subseteq \mathrm{Im}(\tau_p)$. Now similarly, as in the separable case, we can prove that $\sum\limits_{I \in (T_p)_m} \sum\limits_{J \in \overline{T}_P} \overline{\dfrac{t^J g_{I,J}^2}{p^{e+1}}\dfrac{dt_I}{t_I}} + L_{d-1} \subseteq \mathrm{Im}(\tau_p)$, for $g_{I,J} \in \mathcal{F}[x]_{<d}$. Finally, $\tilde{S}_{p,e+1} + \tilde{S}_{p,e+1} \wedge \overline{\dfrac{dp}{p}} + L_{d-1} \subseteq \mathrm{Im}(\tau_p)$. Therefore, $\tilde{S}_p + \tilde{S}_p \wedge \overline{\dfrac{dp}{p}} + L_{d-1} \subseteq \mathrm{Im}(\tau_p)$. □

The next lemma will be used to establish the injectivity of Milnor's splitting $\bigoplus_{\deg p = d} \tau_p$.

**Lemma 4.6.** *If $p$ and $q$ are two distinct monic irreducible polynomials in $\mathcal{F}[x]$ of degree $d$, then the composition $\partial_q \circ \tau_p : W_1(H_2^{m+1}(F_p)) \to W_1(H_2^{m+1}(F_q))$ is the zero map.*

*Proof.* By Lemma 3.3, the map $\partial_q : L_d/L_{d-1} \to W_1(H_2^{m+1}(F_q))$ is well defined. Moreover, by Theorem 4.2, any element $\varphi \in W_1(H_2^{m+1}(F_p))$ has the unique decomposition as $\varphi = \psi + \varphi_2 \wedge \overline{\dfrac{\mathrm{d}p}{p}}$, where $\varphi_2 \in H_2^m(\mathcal{F}(p))$ and $\psi = \displaystyle\sum_{I \in (T_p)_m} \sum_{\{J \in T_p \ | \ J+I > I\}} \overline{t^J r_{I,J}^2 \dfrac{\mathrm{d}t_I}{t_I}}$ such that $r_{I,J} \in p^{-1}\mathcal{F}[x]_{<d}[p^{-1}]$.

Now for simplicity assume $\varphi_2 = \overline{f\dfrac{\mathrm{d}\overline{g_1}}{\overline{g_1}} \wedge \ldots \wedge \dfrac{\mathrm{d}\overline{g_{m-1}}}{\overline{g_{m-1}}}} \in H_2^m(\mathcal{F}(p)) \subset H_2^m(F_p)$, where $f, g_i \in \mathcal{F}[x]_{<d}$ and $g_i \neq 0$. Then, we have:

$$\partial_q \circ \tau_p\left(\varphi_2 \wedge \overline{\dfrac{\mathrm{d}p}{p}}\right) = \partial_q\left(\overline{f\dfrac{\mathrm{d}g_1}{g_1} \wedge \ldots \wedge \dfrac{\mathrm{d}g_{m-1}}{g_{m-1}} \wedge \dfrac{\mathrm{d}p}{p}} + L_{d-1}\right)$$
$$= \partial_q\left(\overline{f\dfrac{\mathrm{d}g_1}{g_1} \wedge \ldots \wedge \dfrac{\mathrm{d}g_{m-1}}{g_{m-1}} \wedge \dfrac{\mathrm{d}p}{p}}\right)$$
$$= 0 \quad \text{(by Lemma 3.3)}.$$

Finally, we will show that $\partial_q(\tau_p(\psi)) = 0$. It suffices to suppose $\psi = \overline{t^J r_{I,J}^2 \dfrac{\mathrm{d}t_I}{t_I}}$, where $r_{I,J} = \dfrac{r'}{p^k}$ such that $r' \in \mathcal{F}[x]_{<d}$ for some $k \geq 1$, and $I \in (T_p)_m, J \in T_p$ satisfy $I + J > I$. Since $p$ is a $q$-adic unit, it follows that $v_q(t^J r_{I,J}^2) \geq 0$ for any $I, J$.

- If $p$ is separable, then $\dfrac{\mathrm{d}t_I}{t_I} = \dfrac{\mathrm{d}t_{i_1}}{t_{i_1}} \wedge \ldots \wedge \dfrac{\mathrm{d}t_{i_m}}{t_{i_m}}$, where $t_{i_l} \in \{p\} \cup \{t_i \ | \ i \in \mathcal{I}\}$ for all $l = 1, 2, \ldots, m$. Therefore, each $t_{i_l}$ is a $q$-adic unit.

- If $p$ is inseparable, then $\dfrac{\mathrm{d}t_I}{t_I} = \dfrac{\mathrm{d}t_{i_1}}{t_{i_1}} \wedge \ldots \wedge \dfrac{\mathrm{d}t_{i_m}}{t_{i_m}}$, where $t_{i_l} \in \{p, x\} \cup \{t_i \ | \ i \in \mathcal{I} \setminus \{j_n\}\}$ for all $l = 1, 2, \ldots, m$. But, in this case, $\deg p = \deg q \geq 2$. So $x$ is a $q$-adic unit. Moreover, each $t_{i_l}$ is also a $q$-adic unit. Therefore, by Lemma 3.3, we get

$$\partial_q \circ \tau_p(\psi) = \partial_q(\psi + L_{d-1}) = \partial_q(\psi) = \partial_q\left(\overline{t^J r_{I,J}^2 \dfrac{\mathrm{d}t_I}{t_I}}\right) = 0.$$

Thus, $\partial_q \circ \tau_p(\varphi) = 0$ for all $\varphi \in W_1(H_2^{m+1}(F_p))$. □

Now we prove that the Milnor's splitting $\bigoplus_{\deg p = d} \tau_p$ is an isomorphism with inverse the second residue map $\bigoplus_{\deg p = d} \partial_p$.

**Theorem 4.7.** *For any positive integer $d$, the Milnor's splitting $\displaystyle\bigoplus_{\deg p = d} \tau_p : \bigoplus_{\deg p = d} W_1(H_2^{m+1}(F_p)) \to L_d/L_{d-1}$ is an isomorphism.*

*Proof.* **(1) The injectivity:** By Theorem 4.2, each $\varphi \in W_1(H_2^{m+1}(F_p))$ has the unique decomposition as $\varphi = \psi + \varphi_2 \wedge \overline{\dfrac{\mathrm{d}p}{p}}$, where

$$\psi = \sum_{I \in (T_p)_m} \sum_{\{J \in T_p \ | \ J+I > I\}} \overline{t^J r_{I,J}^2 \dfrac{\mathrm{d}t_I}{t_I}} \quad \text{and} \quad \varphi_2 = \sum_i \overline{f_i \dfrac{\mathrm{d}\overline{g_{i_1}}}{\overline{g_{i_1}}} \wedge \ldots \wedge \dfrac{\mathrm{d}\overline{g_{i_{m-1}}}}{\overline{g_{i_{m-1}}}}},$$

such that $r_{I,J} \in p^{-1}\mathcal{F}[x]_{<d}[p^{-1}]$, $f_i, g_{i_l} \in \mathcal{F}[x]_{<d}$ for $l = 1, 2, \ldots, m-1$ and $g_{i_l} \neq 0$. Therefore,

$$\partial_p \circ \tau_p(\varphi) = \partial_p \left( \psi + \sum_i \overline{f_i \frac{dg_{i_1}}{g_{i_1}} \wedge \ldots \wedge \frac{dg_{i_{m-1}}}{g_{i_{m-1}}} \wedge \frac{dp}{p}} + L_{d-1} \right)$$

$$= \psi + \sum_i \overline{f_i \frac{dg_{i_1}}{g_{i_1}} \wedge \ldots \wedge \frac{dg_{i_{m-1}}}{g_{i_{m-1}}} \wedge \frac{dp}{p}} \quad \text{(because } \deg f_i, \deg g_{i_l} \leq d-1\text{)}$$

$$= \varphi$$

Thus, $\partial_p \circ \tau_p = \text{Id}$. Now Lemma 4.6 concludes that $\bigoplus_{\deg p = d} \partial_p \circ \bigoplus_{\deg p = d} \tau_p = \text{Id}$. Consequently, $\bigoplus_{\deg p = d} \tau_p$ is injective.

**(2) The surjectivity:** We discuss on the degree $d$ of a monic irreducible polynomial $p$.

(1) Suppose $d > 1$. If $p$ is separable, then by Proposition 3.8, we have

$$S_p + S_p \wedge \overline{\frac{dx}{x}} + L_0 \wedge \overline{\frac{dx}{x}} = S_p + S_p \wedge \overline{\frac{dp}{p}} + L_0 \wedge \overline{\frac{dx}{x}},$$

and by Lemma 4.5, we have

$$S_p + S_p \wedge \overline{\frac{dp}{p}} + L_{d-1} \subseteq \text{Im}(\tau_p).$$

Since $L_0 \wedge \overline{\frac{dx}{x}} \subseteq L_{d-1}$ for all $d > 1$, it follows that $S_p + S_p \wedge \overline{\frac{dx}{x}} + L_{d-1} \subseteq \text{Im}(\tau_p)$.

If $p$ is inseparable, then by Lemma 4.5

$$\tilde{S}_p + \tilde{S}_p \wedge \overline{\frac{dp}{p}} + L_{d-1} \subseteq \text{Im}(\tau_p)$$

and by Proposition 3.8

$$S_p + S_p \wedge \overline{\frac{dx}{x}} + L_0 + L_0 \wedge \overline{\frac{dx}{x}} = \tilde{S}_p + \tilde{S}_p \wedge \overline{\frac{dp}{p}} + L_0 + L_0 \wedge \overline{\frac{dx}{x}}.$$

Since $L_0 + L_0 \wedge \overline{\frac{dx}{x}} \subseteq L_{d-1}$ for $d > 1$, it follows that $S_p + S_p \wedge \overline{\frac{dx}{x}} + L_{d-1} \subseteq \text{Im}(\tau_p)$.

Finally, we obtain $S_p + S_p \wedge \overline{\frac{dx}{x}} + L_{d-1} \subseteq \text{Im}(\tau_p)$ for all polynomials $p$ of degree $d > 1$. Moreover, by Lemma 3.5, we know that $L_d = \sum_{\deg p = d} \left( S_p + S_p \wedge \overline{\frac{dx}{x}} \right) + L_{d-1}$. Therefore, $\bigoplus_{\deg p = d > 1} \tau_p$ is surjective.

(2) Suppose $d = 1$. Then, $p$ is separable. By Proposition 3.8,

(16) $$S_p + S_p \wedge \overline{\frac{dx}{x}} + L_0 \wedge \overline{\frac{dx}{x}} = S_p + S_p \wedge \overline{\frac{dp}{p}} + L_0 \wedge \overline{\frac{dx}{x}}.$$

Moreover, by Lemma 4.5, we have

(17) $$S_p + S_p \wedge \overline{\frac{dp}{p}} + L_0 \subseteq \text{Im}(\tau_p).$$

Putting $p = x$ in (17), we have :

(18) $$S_x \wedge \overline{\frac{dx}{x}} + L_0 \subseteq \text{Im}(\tau_x).$$

By definition,

(19) $$L_0 \wedge \overline{\frac{dx}{x}} \subseteq S_x \wedge \overline{\frac{dx}{x}}.$$

From (16), (17), (18) and (19), we have :

(20) $$S_p + S_p \wedge \overline{\frac{dx}{x}} + L_0 \subseteq \operatorname{Im}(\tau_p) + \operatorname{Im}(\tau_x) \quad \text{for all } p \text{ such that } \deg p = 1.$$

Therefore,

$$S_p + S_p \wedge \overline{\frac{dx}{x}} + L_0 \subseteq \operatorname{Im}\left(\bigoplus_{\deg p=1} \tau_p\right) \quad \text{for all } p \text{ such that } \deg p = 1.$$

Since $L_1 = \sum_{\deg p=1}\left(S_p + S_p \wedge \overline{\frac{dx}{x}}\right) + L_0$, it follows that $L_1/L_0 \subseteq \operatorname{Im}\left(\bigoplus_{\deg p=1} \tau_p\right)$. Hence, $\bigoplus_{\deg p=1} \tau_p$ is surjective and the theorem is proved. □

Using Theorem 4.7, we prove the following short exact sequence, which is analogous to Milnor's exact sequence.

**Theorem 4.8.** *The sequence* $0 \to L_0 \xrightarrow{i} H_2^{m+1}(\mathcal{F}(x)) \xrightarrow{\oplus \partial_p} \bigoplus_p W_1(H_2^{m+1}(F_p)) \to 0$ *is exact, where $i$ is the inclusion map and $p$ varies over all monic irreducible polynomials.*

*Proof.* It is sufficient to prove that $H_2^{m+1}(\mathcal{F}(x))/L_0 \simeq \bigoplus_p W_1(H_2^{m+1}(F_p))$. Recall that $H_2^{m+1}(\mathcal{F}(x)) = \bigcup_{d \geq 0} L_d$. Using an induction on $d$, we prove that $L_d/L_0 \simeq \bigoplus_{\deg p \leq d} W_1(H_2^{m+1}(F_p))$ for all $d \geq 1$. In fact, by Theorem 4.7, we have $L_1/L_0 \simeq \bigoplus_{\deg p \leq 1} W_1(H_2^{m+1}(F_p))$. Suppose that $L_{d-1}/L_0 \simeq \bigoplus_{\deg p \leq d-1} W_1(H_2^{m+1}(F_p))$ for some $d \geq 2$. To prove the isomorphism for $d$, we draw the following diagram where $A_p$ and $d(p)$ denote the group $W_1(H_2^{m+1}(F_p))$ and the degree of $p$, respectively.

$$\begin{array}{ccccccccc}
0 & \longrightarrow & L_{d-1}/L_0 & \longrightarrow & L_d/L_0 & \longrightarrow & L_d/L_{d-1} & \longrightarrow & 0 \\
& & \downarrow \oplus \partial_p & & \downarrow \oplus \partial_p & & \downarrow \oplus \partial_p & & \\
0 & \longrightarrow & \oplus_{d(p) \leq d-1} A_p & \longrightarrow & \oplus_{d(p) \leq d} A_p & \longrightarrow & \oplus_{d(p)=d} A_p & \longrightarrow & 0
\end{array}$$

By induction hypothesis, the first vertical map is an isomorphism. By Theorem 4.7, the last vertical map is an isomorphism. Then, by snake lemma, the middle vertical map is an isomorphism. Consequently, $L_d/L_0 \simeq \bigoplus_{\deg p \leq d} W_1(H_2^{m+1}(F_p))$ for all $d$. Now taking $d \to \infty$, we get

$$H_2^{m+1}(\mathcal{F}(x))/L_0 \simeq \bigoplus_p W_1(H_2^{m+1}(F_p)),$$

where $p$ varies over all monic irreducible polynomials. Thus, the theorem is proved. □

The decomposition of $W_1(H_2^{m+1}(F_{\frac{1}{x}}))$ of Theorem 4.2 and $L_0$ provide us the following result.

**Theorem 4.9.** *The following sequence*

$$0 \to H_2^{m+1}(\mathcal{F}) \xrightarrow{i_1} L_0 \xrightarrow{i_2} W_1(H_2^{m+1}(F_{\frac{1}{x}})) \xrightarrow{\eta} H_2^m(\mathcal{F}) \to 0$$

*is exact, where $i_1$ is induced by inclusion and $i_2$ is the restriction of $\partial_{\frac{1}{x}}$ to $L_0$. Moreover, $\eta(\varphi) = \varphi_2$, where $\varphi_2$ comes from the decomposition of $\varphi$, stated in Theorem 4.2(2).*

*Proof.* The injectivity of $i_1$ follows from [2, Lemma 2.17]. Here we give an alternative proof. Let $\chi_x$ be the residue map given in Definition 6.12, with respect to the $x$-adic valuation. Now we compute $\chi_x \circ i_1$. Let $a_j \in \mathcal{F}$ and $c_{j_l} \in \mathcal{F}^*$. Then we have :

$$\chi_x \circ i_1 \left( \overline{\sum_j a_j \frac{dc_{j_1}}{c_{j_1}} \wedge \ldots \wedge \frac{dc_{j_m}}{c_{j_m}}} \right) = \chi_x \left( \overline{\sum_j a_j \frac{dc_{j_1}}{c_{j_1}} \wedge \ldots \wedge \frac{dc_{j_m}}{c_{j_m}}} \right)$$

$$= f_{m+1}^{-1} \circ \Delta_x \circ g_{m+1} \left( \overline{\sum_j a_j \frac{dc_{j_1}}{c_{j_1}} \wedge \ldots \wedge \frac{dc_{j_m}}{c_{j_m}}} \right)$$

$$= f_{m+1}^{-1} \circ \Delta_x \left( \overline{\sum_j \langle\langle c_{j_1}, \ldots, c_{j_m}; a_j]]} \right)$$

$$= f_{m+1}^{-1} \left( \overline{\sum_j \langle\langle \overline{c_{j_1}}, \ldots, \overline{c_{j_m}}; \overline{a_j}]]} \right)$$

$$= \overline{\sum_j \overline{a_j} \frac{d\overline{c_{j_1}}}{\overline{c_{j_1}}} \wedge \ldots \wedge \frac{d\overline{c_{j_m}}}{\overline{c_{j_m}}}} \in H_2^{m+1}(\overline{F}_x)$$

$$= \overline{\sum_j a_j \frac{dc_{j_1}}{c_{j_1}} \wedge \ldots \wedge \frac{dc_{j_m}}{c_{j_m}}} \in H_2^{m+1}(\mathcal{F}).$$

Here $f_{m+1}, g_{m+1}$ are Kato's isomorphisms. Therefore, $\chi_x \circ i_1 = \text{Id}$. Consequently, $i_1$ is injective.

Moreover, by the decomposition of Theorem 4.2, $\eta$ is onto. Therefore, it remains to prove that $\ker(i_2) = H_2^{m+1}(\mathcal{F})$ and $\text{Im}(i_2) = \ker(\eta)$. To do so, we need to prove that any element $\varphi \in L_0$ can be decomposed as $\varphi = \varphi_1 + \psi$, where $\varphi_1 \in H_2^{m+1}(\mathcal{F})$ and

(21) $$\psi = \overline{\sum_{I \in (T_{\frac{1}{x}})_m} \sum_{\{J \in T_{\frac{1}{x}} \mid I+J > I\}} t^J r_{I,J}^2 \frac{dt_I}{t_I}}$$

such that $r_{I,J} \in x\mathcal{F}[x]$. Recall that by definition, $L_0$ is generated by two types of forms: $\overline{h \frac{dc_1}{c_1} \wedge \ldots \wedge \frac{dc_m}{c_m}}$ and $\overline{xf \frac{de_1}{e_1} \wedge \ldots \wedge \frac{de_{m-1}}{e_{m-1}} \wedge \frac{dx}{x}}$, for $h, f \in \mathcal{F}[x]$ and $c_i, e_j \in \mathcal{F}^*$. Now using the fact $\frac{dx}{x} = \frac{dx^{-1}}{x^{-1}}$ we have :

$$\overline{h \frac{dc_1}{c_1} \wedge \ldots \wedge \frac{dc_m}{c_m}} + \overline{xf \frac{de_1}{e_1} \wedge \ldots \wedge \frac{de_{m-1}}{e_{m-1}} \wedge \frac{dx}{x}} = \sum_{I \in (T)_m} \sum_{J \in T} \overline{t^J a_{I,J}^2 \frac{dt_I}{t_I}} + \sum_{k=1}^{l} \sum_{I \in (T_{\frac{1}{x}})_m} \sum_{J \in T} \overline{x^k t^J b_{I,J}^2 \frac{dt_I}{t_I}},$$

where $a_{I,J}, b_{I,J} \in \mathcal{F}$ and $l \in \mathbb{N}$. Now the first sum of right side is in $H_2^{m+1}(\mathcal{F})$.

Using the same argument as in the Cases 1 and 2 in the proof of Lemma 4.5(1), we prove that for $k \geq 1$ and $b_{I,J} \in \mathcal{F}$, the element $\sum_{I \in (T_{\frac{1}{x}})_m} \sum_{J \in T} \overline{x^k t^J b_{I,J}^2 \frac{dt_I}{t_I}}$ can be expressed like $\psi$ of (21). Thus, each element $\varphi \in L_0$ has an expression $\varphi = \varphi_1 + \psi$, as defined previously and $\psi$ is unique by Theorem 4.2. Consequently, $\ker(i_2) = H_2^{m+1}(\mathcal{F})$ and $\mathrm{Im}(i_2) = W_1(H_2^{m+1}(F_{\frac{1}{x}}))/H_2^m(\mathcal{F}) \wedge \overline{\frac{dx}{x}} = \ker(\eta)$. Hence, the theorem is proved. □

Now we prove our main exact sequence, which is an extension of Milnor's exact sequence to Kato-Milne cohomology in characteristic 2.

**Theorem 4.10.** *Let $\mathcal{F}$ be a field of characteristic 2 and $F = \mathcal{F}(x)$ the rational function field in one variable $x$ over $\mathcal{F}$. Then, the following sequence is short exact:*

$$0 \to H_2^{m+1}(\mathcal{F}) \xrightarrow{i_3} H_2^{m+1}(F) \xrightarrow{\oplus_p \partial_p \oplus \partial'_{\frac{1}{x}}} \oplus_p W_1(H_2^{m+1}(F_p)) \oplus \left( W_1(H_2^{m+1}(F_{\frac{1}{x}}))/H_2^m(\mathcal{F}) \wedge \overline{\frac{dx}{x}} \right) \to 0,$$

*where $i_3$ is induced by the inclusion $\mathcal{F} \hookrightarrow F$ and the direct sum is taken over $\frac{1}{x}$ and all monic irreducible polynomials $p \in \mathcal{F}[x]$. Moreover, $\partial'_{\frac{1}{x}}$ is the composition of $\partial_{\frac{1}{x}}$ and the projection map.*

*Proof.* By definition $i_3 = i \circ i_1$, where $i$ and $i_1$ come from Theorems 4.8 and 4.9 respectively. Since $i$ and $i_1$ are injective, it follows that $i_3$ is injective. Moreover, $\mathrm{Im}(i_3) = \mathrm{Im}(i_1)$ as $i$ is inclusion map. Now we compute

$$\begin{aligned}
\ker(\oplus_p \partial_p \oplus \partial'_{1/x}) &= \ker(\oplus_p \partial_p) \cap \ker(\partial'_{1/x}) \\
&= L_0 \cap \ker(\partial'_{1/x}) \quad \text{(By Theorem 4.8)} \\
&= \ker(i_2) \quad \text{(By Theorem 4.9)} \\
&= \mathrm{Im}(i_1) \\
&= \mathrm{Im}(i_3).
\end{aligned}$$

Hence, it remains to prove that $\oplus_p \partial_p \oplus \partial'_{\frac{1}{x}}$ is onto. Let $(u, \bar{v})$ be an arbitrary element of $\oplus_p W_1(H_2^{m+1}(F_p)) \oplus \left( W_1(H_2^{m+1}(F_{\frac{1}{x}}))/H_2^m(\mathcal{F}) \wedge \overline{\frac{dx}{x}} \right)$ such that $u \in \oplus_p W_1(H_2^{m+1}(F_p))$ and $v \in W_1(H_2^{m+1}(F_{\frac{1}{x}}))$. By Theorem 4.8, the map $\oplus_p \partial_p$ is onto and it follows that there exists an element $\varphi \in H_2^{m+1}(\mathcal{F}(x))$ such that

$$\oplus_p \partial_p(\varphi) = u.$$

Let $\partial_{\frac{1}{x}}(\varphi) = v' + v''$, where $v'' \in H_2^m(\mathcal{F}) \wedge \overline{\frac{dx}{x}}$ and $v' = \sum_{I \in (T_{\frac{1}{x}})_m} \sum_{\{J \in T_{\frac{1}{x}} \mid J+I>I\}} \overline{t^J r_{I,J}^2 \frac{dt_I}{t_I}}$, $r_{I,J} \in x\mathcal{F}[x]$. Then, we get

(22) $$\partial'_{\frac{1}{x}}(\varphi) = \overline{v'} \in W_1(H_2^{m+1}(F_{\frac{1}{x}}))/H_2^m(\mathcal{F}) \wedge \overline{\frac{dx}{x}}.$$

Now by Theorem 4.2, $v = v_1 + v_2$, where $v_2 \in H_2^m(\mathcal{F}) \wedge \overline{\frac{dx}{x}}$ and $v_1 = \sum_{I \in (T_{\frac{1}{x}})_m} \sum_{J \in T_{\frac{1}{x}} \mid J+I>I} \overline{t^J r_{I,J}'^2 \frac{dt_I}{t_I}}$ with $r_{I,J}' \in x\mathcal{F}[x]$. By Theorem 4.9, there exists $\varphi' \in L_0$ such that $i_2(\varphi') = v_1 + v'$. Since $i_2$ is the restriction of $\partial'_{\frac{1}{x}}$ on $L_0$, we have

$$\partial'_{\frac{1}{x}}(\varphi') = \overline{v_1 + v'} \in W_1(H_2^{m+1}(F_{\frac{1}{x}}))/H_2^m(\mathcal{F}) \wedge \overline{\frac{dx}{x}}.$$

Now using (22) we get,

$$\partial'_{\frac{1}{x}}(\varphi + \varphi') = \overline{v}_1 = \overline{v} \in W_1(H_2^{m+1}(F_{\frac{1}{x}}))/H_2^m(\mathcal{F}) \wedge \overline{\frac{dx}{x}}.$$

Since $\varphi' \in L_0$, it follows from Theorem 4.8 that $\oplus_p \partial_p(\varphi') = 0$. Therefore, $\oplus_p \partial_p \oplus \partial'_{\frac{1}{x}}(\varphi + \varphi') = (u, \overline{v})$. Thus, the theorem is proved. $\square$

## 5. The norm theorem for Kato-Milne cohomology

Our aim in this section is to give an important application of Theorem 4.10 concerning the norm theorem in the setting of Kato-Milne cohomology. First of all we define the notions of hyperbolicity and norm of a differential form, which were first introduced by Mukhija in [17, Definition A.1, A.2].

**Definition 5.1.** *Let $F$ be a field of characteristic 2 and $w \in \Omega_F^m$, a differential form of degree some $m \geq 1$. The form $w$ is said to be hyperbolic if $w \in d\Omega_F^{m-1} + \wp(\Omega_F^m)$, which means that $\overline{w} = 0 \in H_2^{m+1}(F)$.*

This definition was made doing a parallel between the Witt group of nonsingular quadratic forms of $F$ and the group $H_2^{m+1}(F)$ of differential forms. Similarly, we may extend to differential forms the notion of norm. Recall that a scalar $\alpha \in F^*$ is a called a norm of a nonsingular quadratic form $\varphi$ over $F$ if $\varphi \simeq \alpha\varphi$, which means that $\langle 1, \alpha \rangle_b \otimes \varphi$ is hyperbolic. The class of the bilinear form $\langle 1, \alpha \rangle_b$ in $I(F)/I^2(F)$ corresponds to $\frac{d\alpha}{\alpha} \in \nu_F(1)$ by Kato's isomorphism, and the module action of $W(F)$ on $W_q(F)$, induced by tensor product, has a simialr action of $\nu_F(n)$ on $H_2^{m+1}(F)$ induced by exterior product. This motivates the following definition:

**Definition 5.2.** *Let $w \in \Omega_F^m$ and $\alpha \in F^*$. The scalar $\alpha$ is said to be a norm of $w$ if $w \wedge \frac{d\alpha}{\alpha}$ is hyperbolic.*

The norm theorem started by Knebusch in the seventies who proved the following: *If $k$ is a field of arbitrary characteristic, $p \in A := k[x_1, \cdots, x_n]$ a normed[1] irreducible polynomial in $n$ variables, then an anisotropic bilinear form $B$ over $k$ becomes metabolic over $k(p)$ iff $B \simeq pB$ over $k(x_1, \cdots, x_n)$, where $k(p)$ is the function field of the affine hypersurface given by $p$, i.e., the quotient field of the ring $k[x_1, \cdots, x_n]/(p)$* [12, Theorem 4.2]. In particular, we get the norm theorem for quadratic forms in characteristic not 2. Later, in characteristic 2, Baeza proved the norm theorem for nonsingular quadratic forms [6]. For totally singular quadratic forms, the norm theorem was completed by the first author [13], and independently by Hoffmann [9]. Recently, the first author and Mukhija proved the norm theorem for semisingular quadratic forms [14], thus completing the picture of the norm theorem for quadratic forms. Based on the notion of transfer defined by Aravire, Laghribi and O'Ryan in the setting of Kato-Milne cohomology for purely inseparable extensions [5], and using Definitions 5.1 and 5.2, Mukhija proved the norm theorem for Kato-Milne cohomology only for inseparable irreducible polynomials [17, Theorem A.3]. Our aim in this section is to complete the open case that concerns normed irreducible separable polynomials. First, in the following theorem, we treat the case of an arbitrary irreducible polynomial in one variable, and then the case of many variables in Theorem 5.5.

**Theorem 5.3.** *Let $\mathcal{F}$ be a field of characteristic 2, $F = \mathcal{F}(x)$ and $w \in \Omega_\mathcal{F}^m$. Let $p \in \mathcal{F}[x]$ be monic irreducible polynomial (separable or not). Then, the following statements are equivalent:*
*(1) $w$ is hyperbolic over $\mathcal{F}(p)$.*
*(2) $p$ is a norm of $w_F$.*

*Proof.* (2) $\Longrightarrow$ (1) This is done in [17, Proposition A.7].

(1) $\Longrightarrow$ (2) Let $w \in \Omega_F^m$ be hyperbolic over $\mathcal{F}(p)$. Then, $\overline{w} = 0 \in H_2^{m+1}(\mathcal{F}(p))$. If $\deg p$ is odd, then $\overline{w} = 0 \in H_2^{m+1}(\mathcal{F})$ because $\ker(H_2^{m+1}(\mathcal{F}) \to H_2^{m+1}(\mathcal{F}(p))) = 0$ [17, Theorem A.16]. Therefore, $\overline{w \wedge \frac{dp}{p}} =$

---
[1]The polynomial $p$ is normed means that the leading coefficient of $p$ with respect to the lexicographic ordering is 1.

$0 \in H_2^{m+2}(F)$, and thus $p$ is a norm of $w$ over $F$. So suppose that deg $p$ is even. We will show that $\overline{w \wedge \frac{dp}{p}} \in \ker\left(\oplus_q \partial_q \oplus \partial_{1/x}\right)$, where $q$ varies over all monic irreducible polynomials.

Since $\overline{w} = 0 \in H_2^{m+1}(\mathcal{F}(p))$, it follows from Corollary 6.9 that $\overline{w} = 0$ in $H_2^{m+1}(F_p)$. Therefore, $\overline{w \wedge \frac{dp}{p}} = 0 \in H_2^{m+2}(F_p)$. By Theorem 4.2, we have

(23) $$\partial_p\left(\overline{w \wedge \frac{dp}{p}}\right) = \overline{w \wedge \frac{dp}{p}} = 0.$$

Let $q \in \mathcal{F}[x]$ be monic and irreducible such that $p \neq q$. Then, $p$ is a $q$-adic unit. Let us write $w = \sum_j a_j \frac{dc_{j_1}}{c_{j_1}} \wedge \ldots \wedge \frac{dc_{j_m}}{c_{j_m}}$, where $a_j, c_{j_l} \in \mathcal{F}^*$. Then, for each $j, l$, we have $v_q(c_{j_l}) = 0 = v_q(a_j)$. Therefore, by Lemma 3.3, we get:

(24) $$\partial_q\left(\overline{w \wedge \frac{dp}{p}}\right) = 0 \quad \text{for all monic irreducible } q \neq p.$$

To compute $\partial_{\frac{1}{x}}$, let us write $p = x^{2k} + p_{2k-1}x^{2k-1} + \cdots + p_1 + p_0$, where $p_i \in \mathcal{F}$ and $k \geq 1$. We rewrite $p$ as follows:

$$p = x^{2k}\left(\sum_{i=0}^{2k} p_{2k-i} u^i\right), \quad \text{where } u = \frac{1}{x} \text{ and } p_{2k} = 1.$$

Then, we have:

$$\overline{w \wedge \frac{dp}{p}} = \overline{w \wedge \frac{d(\sum_{i=0}^{2k} p_{2k-i} u^i)}{\sum_{i=0}^{2k} p_{2k-i} u^i}}$$

$$= \sum_{i=0}^{2k} \overline{\frac{p_{2k-i}u^i}{p'} w \wedge \frac{d(p_{2k-i} u^i)}{p_{2k-i} u^i}} \quad \left(\text{where } p' = \sum_{i=0}^{2k} p_{2k-i} u^i\right)$$

$$= \sum_j \sum_{i=0}^{2k} \overline{\frac{p_{2k-i} u^i a_j}{p'} \frac{dc_{j_1}}{c_{j_1}} \wedge \ldots \wedge \frac{dc_{j_m}}{c_{j_m}} \wedge \frac{d(p_{2k-i} u^i)}{p_{2k-i} u^i}}.$$

Since $u = \frac{1}{x}$ and $p'$ is $\frac{1}{x}$-unit, $v_{\frac{1}{x}}\left(\frac{p_{2k-i} u^i a_j}{p'}\right) \geq 1$ for $i \geq 1$. Therefore, $\frac{p_{2k-i} u^i a_j}{p'} \in \wp(F_{\frac{1}{x}})$ for each $i \geq 1$. Hence, we have in $H_2^{m+2}(F_{\frac{1}{x}})$:

$$\overline{w \wedge \frac{dp}{p}} = \sum_j \overline{\frac{p_{2k} a_j}{p'} \frac{dc_{j_1}}{c_{j_1}} \wedge \ldots \wedge \frac{dc_{j_m}}{c_{j_m}} \wedge \frac{dp_{2k}}{p_{2k}}}$$

$$= 0 \quad (\text{as } p_{2k} = 1 \text{ and } d1 = 0).$$

Consequently,

(25) $$\partial_{1/x}\left(\overline{w \wedge \frac{dp}{p}}\right) = 0 = \partial'_{\frac{1}{x}}\left(\overline{w \wedge \frac{dp}{p}}\right).$$

By (23), (24) and (25), $\overline{w \wedge \frac{dp}{p}} \in \ker\left(\oplus_q \partial_q \oplus \partial'_{1/x}\right)$, and Theorem 4.10 implies that $\overline{w \wedge \frac{dp}{p}} \in \text{Im}(i_3)$.

Let $\alpha = \sum_k \overline{s_k \frac{dh_{k_1}}{h_{k_1}} \wedge \ldots \wedge \frac{dh_{k_{m+1}}}{h_{k_{m+1}}}} \in H_2^{m+2}(\mathcal{F})$ be such that $\overline{w \wedge \frac{dp}{p}} = i_3(\alpha)$. By applying the residue map $\chi_\pi$ from Definition 6.12, for $\pi = 1/x$, we get $\chi_\pi\left(\overline{w \wedge \frac{dp}{p}}\right) = \chi_\pi(i_3(\alpha))$. This gives

$$f_{m+2}^{-1} \circ \Delta_\pi \circ g_{m+2}\left(\sum_j \overline{a_j \frac{dc_{j_1}}{c_{j_1}} \wedge \ldots \wedge \frac{dc_{j_m}}{c_{j_m}} \wedge \frac{dp}{p}}\right) = f_{m+2}^{-1} \circ \Delta_\pi \circ g_{m+2}\left(\sum_k \overline{s_k \frac{dh_{k_1}}{h_{k_1}} \wedge \ldots \wedge \frac{dh_{k_{m+1}}}{h_{k_{m+1}}}}\right)$$

$$f_{m+2}^{-1} \circ \Delta_\pi\left(\sum_j \overline{\langle\langle c_{j_1}, \ldots, c_{j_m}, p; a_j]]}\right) = f_{m+2}^{-1} \circ \Delta_\pi\left(\sum_k \overline{\langle\langle h_{k_1}, \ldots h_{k_{m+1}}; s_k]]}\right)$$

$$f_{m+2}^{-1}\left(\sum_j (\partial_\pi^1 + \partial_\pi^2)(\langle\langle c_{j_1}, \ldots, c_{j_m}, p\rangle\rangle_b) \otimes [1, \overline{a_j}]\right) = f_{m+2}^{-1}\left(\sum_k (\partial_\pi^1 + \partial_\pi^2)(\langle\langle h_{k_1}, \ldots, h_{k_{m+1}}\rangle\rangle_b) \otimes [1, \overline{s_k}]\right).$$

By Definition 6.12, $\partial_\pi^1(\langle c_{j_l}\rangle_b) = \langle \overline{c_{j_l}}\rangle_b$, $\partial_\pi^1(\langle h_{k_l}\rangle_b) = \langle \overline{h_{k_l}}\rangle_b$ and $\partial_\pi^1(\langle p\rangle_b) = \langle \overline{p'}\rangle_b = \langle 1\rangle_b$. Moreover, as each of $p', c_{j_l}, h_{k_l}$ is $\frac{1}{x}$-adic unit, so $\partial_\pi^2(\langle c_{j_l}\rangle_b) = 0 = \partial_\pi^2(\langle h_{k_l}\rangle_b)$ and $\partial_\pi^2(\langle p\rangle_b) = \partial_\pi^2(\langle p'\rangle_b) = 0$. Using these relations, we get $(\partial_\pi^1 + \partial_\pi^2)(\langle\langle c_{j_1}, \ldots, c_{j_m}, p\rangle\rangle_b) = 2\langle\langle \overline{c_{j_1}}, \ldots, \overline{c_{j_m}}\rangle\rangle_b = 0$ and $(\partial_\pi^1 + \partial_\pi^2)(\langle\langle h_{k_1}, \ldots, h_{k_{m+1}}\rangle\rangle_b) = \langle\langle \overline{h_{k_1}}, \ldots, \overline{h_{k_{m+1}}}\rangle\rangle_b$. Hence, we get:

$$f_{m+2}^{-1}\left(\sum_k \overline{\langle\langle \overline{h_{k_1}}, \ldots, \overline{h_{k_{m+1}}}; \overline{s_k}]]}\right) = 0$$

$$\sum_k \overline{\overline{s_k} \frac{d\overline{h_{k_1}}}{\overline{h_{k_1}}} \wedge \ldots \wedge \frac{d\overline{h_{k_{m+1}}}}{\overline{h_{k_{m+1}}}}} = \alpha = 0.$$

The last equation holds over $\overline{F}_{\frac{1}{x}} \simeq \mathcal{F}$. Hence, over $F = \mathcal{F}(x)$, we have:

$$\overline{w \wedge \frac{dp}{p}} = i_3(\alpha) = 0.$$

Consequently, the polynomial $p$ is a norm of $w$ over $F$. Thus, the theorem is proved. □

To generalize Theorem 5.3 in any number of variables we need the following lemma:

**Lemma 5.4.** *Let $\mathcal{F}$ be a field of characteristic 2, $F = \mathcal{F}(x_1, \ldots, x_n)$ the rational function field in the variables $x_1, \cdots, x_n$ over $\mathcal{F}$, and $w \in \Omega_\mathcal{F}^m$. Let $p \in \mathcal{F}[x_1, \ldots, x_n]$ be a normed irreducible polynomial and $f \in \mathcal{F}[x_1, \ldots, x_n]$. If $f$ is a norm of $w$ over $F$, and $p$ divides $f$ in exactly odd power, then $p$ is also a norm of $w$ over $F$.*

*Proof.* We use the same proof as given in [17, Proposition A.15], just in Step 1 of this proof one has to use Theorem 5.3 instead of [17, Proposition A.12]. □

Finally, we prove the norm theorem for normed irreducible polynomials in arbitrary number of variables. The method is similar to that used by Mukhija in [17], which is inspired from the induction argument used by Knebusch in the case of bilinear forms [12].

**Theorem 5.5.** *Let $\mathcal{F}$ be a field of characteristic 2, $F = \mathcal{F}(x_1, \ldots, x_n)$ and $w \in \Omega_\mathcal{F}^m$. Let $p \in \mathcal{F}[x_1, \ldots, x_n]$ be a normed irreducible polynomial. Then, the following statements are equivalent:*
*(1) $w$ is hyperbolic over $\mathcal{F}(p)$.*
*(2) $p$ is a norm of $w_F$.*

*Proof.* The case $n = 1$ is given by Theorem 5.3. So suppose $n > 1$.

$(2) \Longrightarrow (1)$ This is done in [17, Proposition A.7].

$(1) \Longrightarrow (2)$ If needed, by renaming the variables, we consider $p \in \mathcal{F}(x_2, \ldots, x_n)[x_1]$ such that the leading coefficient of $p$ in $\mathcal{F}(x_2, \ldots, x_n)$ is nonzero. Let $K = \mathcal{F}(x_2, \ldots, x_n)$ and $g$ the leading coefficient of $p \in K[x_1]$. Since $g$ is a unit in $K$, $\mathcal{F}(p) = K(p) = K(g^{-1}p)$. Moreover, $w$ is hyperbolic over $\mathcal{F}(p) = K(g^{-1}p)$ and $g^{-1}p$ is a monic irreducible polynomial in $K[x_1]$. Therefore, by the case $n = 1$, $g^{-1}p$ is a norm of $w_{K(x_1)}$, and this implies that $gp$ is a norm of $w_{K(x_1)}$. Now Lemma 5.4 implies that $p$ is a norm of $w_{K(x_1)} = w_F$. Thus, the theorem is proved. □

## 6. Appendix

Our aim in this appendix is to give an explicit description of the Teichmüller lifting map that played an important role in the previous sections. The results we present here are from [18, chapter 7]. But to make easy to understand we are providing proofs. Throughout this section, $F$ is a complete field of characteristic 2 with respect to a discrete valuation $v : F \to \mathbb{Z}$, except stated. Let $\overline{F}$ be the residue field of $F$ and $A_F$ the valuation ring of $F$. Moreover, $A_F^\times$ denotes the group of units of $A_F$. To define the Teichmüller lifting, we need the following results.

**Theorem 6.1.** *If the residue field $\overline{F}$ of the complete field $F$ is perfect, then there exists a unique field monomorphism $\beta : \overline{F} \to F$ such that $\mathrm{Im}(\beta) \subseteq A_F$, and $\overline{\beta(a)} = a$ for all $a \in \overline{F}$.*

*Proof.* **Existence:** Let $a$ be an arbitrary element of $\overline{F}$. Since $\overline{F}$ is perfect, $a^{2^{-n}} \in \overline{F}$ for all $n \in \mathbb{N}_0$. Suppose $a_n \in A_F$ such that

$$\overline{a}_n = a^{2^{-n}} \quad \forall n \in \mathbb{N}_0. \tag{26}$$

Consequently, we have the following relations:

$$\overline{a}_{n+1}^2 = \overline{a}_n \quad \forall n \in \mathbb{N}_0. \tag{27}$$

$$a_{n+1}^2 \equiv a_n \pmod{P} \quad (P \text{ is the prime ideal of } A_F) \tag{28}$$

$$a_{n+1}^{2^{n+1}} \equiv a_n^{2^n} \pmod{P^{n+1}} \quad \text{(by [18, Lemma 33])}. \tag{29}$$

Let $k, m, n \in \mathbb{N}_0$ be such that $n \geq m \geq k$. Then, we have:

$$a_n^{2^n} - a_m^{2^m} = a_n^{2^n} - a_{n-1}^{2^{n-1}} + a_{n-1}^{2^{n-1}} - \ldots + a_{m+1}^{2^{m+1}} - a_m^{2^m} \in P^k \quad \text{(by (29))}.$$

Therefore,

$$v\left(a_n^{2^n} - a_m^{2^m}\right) \geq k \quad \forall n \geq m \geq k.$$

Thus, $(a_n^{2^n})_n$ is a Cauchy sequence. Since $F$ is complete, there exists an element $x_a \in A_F$ such that $\lim_{n \to \infty} a_n^{2^n} = x_a$. Moreover, (26) implies that $\overline{x}_a = a \in \overline{F}$. Now we prove that this $x_a$ does not depend on the choice of the elements $a_n$. If possible, let $a'_n \in A_F$ be such that $\overline{a'}_n = a^{2^{-n}}$ for all $n \in \mathbb{N}_0$. Then

$$\overline{a'}_n = \overline{a}_n, \quad \forall n \in \mathbb{N}_0, \quad \text{i.e.,} \quad a'_n \equiv a_n \pmod{P} \; \forall n \in \mathbb{N}_0.$$

Therefore, $a'^{2^n}_n \equiv a_n^{2^n} \pmod{P^{n+1}}$ by [18, Lemma 33]. Hence, $\lim_{n \to \infty} a'^{2^n}_n = \lim_{n \to \infty} a_n^{2^n} = x_a$. Now for each $a \in \overline{F}$, we can define such $x_a \in A_F$. Therefore, we have the map $\beta : \overline{F} \to F$ defined by:

$$\beta(a) = x_a \quad \forall a \in \overline{F}.$$

In particular, if $a_n \in A_F$ such that $(\overline{a}_n)^{2^n} = a$ for all $n \in \mathbb{N}_0$, then we have

$$\beta(a) = \lim_{n \to \infty} a_n^{2^n} \quad \forall a \in \overline{F}. \tag{30}$$

Now we prove that $\beta$ is a field homomorphism. Let $a, b \in \overline{F}$ and $a_n, b_n \in A_F$ such that $(\overline{a_n})^{2^n} = a$ and $(\overline{b_n})^{2^n} = b$ for all $n \in \mathbb{N}_0$. Then, $(\overline{a_n b_n})^{2^n} = ab$ and $(\overline{a_n + b_n})^{2^n} = a + b$, for all $n \in \mathbb{N}_0$. Therefore, we have

$$\begin{aligned} \beta(ab) &= \lim_{n \to \infty} (a_n b_n)^{2^n} \\ &= \lim_{n \to \infty} a_n^{2^n} \lim_{n \to \infty} b_n^{2^n} \quad \text{(the multiplication is a continuous function in } F) \\ &= \beta(a)\beta(b). \end{aligned}$$

Again we get

$$\begin{aligned} \beta(a+b) = \lim_{n \to \infty} (a_n + b_n)^{2^n} &= \lim_{n \to \infty} (a_n^{2^n} + b_n^{2^n}) \quad \text{(as char}(F) = 2) \\ &= \lim_{n \to \infty} a_n^{2^n} + \lim_{n \to \infty} b_n^{2^n} \quad \text{(the addition is a continuous function in } F) \\ &= \beta(a) + \beta(b). \end{aligned}$$

Moreover, $\beta(1) = 1$. Thus, $\beta$ is a field homomorphism such that $\overline{\beta(a)} = a$ for all $a \in \overline{F}$. Hence, $\beta$ is injective.

**Uniqueness:** Let $\beta' : \overline{F} \to F$ be another field monomorphism such that $\operatorname{Im}(\beta') \subseteq A_F$ and $\overline{\beta'(a)} = a$ for all $a \in \overline{F}$. We will prove that $\beta = \beta'$. In fact, let $a \in \overline{F}$. Let $\beta'(a) = y$ and $\beta'(a^{2^{-n}}) = y_n$ for all $n \in \mathbb{N}_0$. Hence,

(31) $$\beta'(a) = (y_n)^{2^n} \quad \forall n \in \mathbb{N}_0$$
(32) $$= y.$$

Now by the definition of $\beta'$, we have $\overline{y} = a$. Using (31) and (32) we get

$$\overline{y}_n = \overline{\beta'(a)}^{2^{-n}} = \overline{y}^{2^{-n}} = a^{2^{-n}} \quad \forall n \in \mathbb{N}_0.$$

Now by the definition of $\beta$ from (30), we have:

$$\begin{aligned} \beta(a) &= \lim_{n \to \infty} (y_n)^{2^n} \quad \text{(since } a = (\overline{y_n})^{2^n} \forall n \in \mathbb{N}_0) \\ &= y \quad \text{( using (31) and (32))} \\ &= \beta'(a) \quad \forall a \in \overline{F}. \end{aligned}$$

Thus, the Theorem is proved. $\square$

The next theorem gives a description of $F$ using the map $\beta$ of Theorem 6.1.

**Theorem 6.2.** *Let $F$, $\overline{F}$ and $\beta$ be as defined in Theorem 6.1. Then, $F = \beta(\overline{F})((\pi))$, where $\pi$ is an uniformizer.*

*Proof.* Clearly, $\beta(\overline{F})((\pi))$ is a complete subfield of $F$. Moreover, the value groups and residue fields of $F$ and $\beta(\overline{F})((\pi))$ are the same. Let $a \in F^*$. Then, there exists an element $l \in \beta(\overline{F})((\pi))$ such that $v(a) = v(l)$. Therefore, $v(a/l) = 0$. Then, there exists an element $c_0 \in \beta(\overline{F})$ such that :

$$a/l \equiv c_0 \pmod{P_F^{n_1}}, \text{ where } P_F = \langle \pi \rangle \text{ is the prime ideal of } A_F \text{ and } v(a/l - c_0) = n_1 \geq 1.$$

Then, $v\left(\frac{a/l - c_0}{\pi^{n_1}}\right) = 0$. Again there exists an element $c_1 \in \beta(\overline{F})$ such that

$$\frac{a/l - c_0}{\pi^{n_1}} \equiv c_1 \pmod{P_F^{n_2}}, \text{ where } v\left(\frac{a/l - c_0}{\pi^{n_1}} - c_1\right) = n_2 \geq 1.$$

Then, $v\left(\frac{a/l - c_0 - c_1 \pi^{n_1}}{\pi^{n_1+n_2}}\right) = 0$. Similarly, there exists an element $c_2 \in \beta(\overline{F})$ such that:

$$\frac{a/l - c_0 - c_1 \pi^{n_1}}{\pi^{n_1+n_2}} \equiv c_2 \pmod{P_F^{n_3}}, \text{ where } v\left(\frac{a/l - c_0 - c_1 \pi^{n_1}}{\pi^{n_1+n_2}} - c_2\right) = n_3 \geq 1.$$

Thus, by induction we have a Cauchy sequence $\left( \sum_{k=0}^{r} c_k \pi^{n_0+n_1+\ldots+n_k} \right)_{r \in \mathbb{N}}$, where $n_0 = 0$ and $c_i \in \beta(\overline{F})$. Since $\beta(\overline{F})((\pi))$ is a complete field, there exists an element $b \in \beta(\overline{F})((\pi))$ such that

$$b = \sum_{k=0}^{\infty} c_k \pi^{n_0+n_1+n_2+\ldots+n_k}.$$

Since $a/l \equiv c_0 + c_1 \pi^{n_1} + \ldots + c_k \pi^{n_1+n_2+\ldots+n_k} \pmod{P_F^{n_1+\ldots+n_{k+1}}}$ for all $k \in \mathbb{N}$, it follows that $a/l = b$. Hence, $a \in \beta(\overline{F})((\pi))$ since $l, b \in \beta(\overline{F})((\pi))$. Therefore, $F = \beta(\overline{F})((\pi))$. □

In the sequel, we will generalize this notion for a complete valued field $F$ whose residue field $\overline{F}$ is imperfect. If $\overline{F}$ is imperfect, then it has a nontrivial 2-basis. Let $S$ be a 2-independent subset of the residue field $\overline{F}$. We consider a subset $S' \subseteq A_F$ which is made of exactly one representative for each element in $S$. Then, $S'$ will be 2-independent in $F$.

**Lemma 6.3.** *Suppose $K = F(b_{1,1}, b_{1,2}, \ldots, b_{1,r})$, where $b_{1,i}^2 = b_i$ for $b_i \in S'$. So $K$ is a finite extension of $F$. Hence, we can extend the valuation $v$ uniquely on $K$. If $\overline{K}$ is the residue field of $K$, then we have:*
*(1) $[K : F] = 2^r$.*
*(2) $[\overline{K} : \overline{F}] = 2^r$ and $\mathcal{T}(K) = \mathcal{T}(F)$, where $\mathcal{T}(K)$ and $\mathcal{T}(F)$ are the value groups of $K$ and $F$, respectively.*

*Proof.* Since $K = F(b_{1,1}, b_{1,2}, \ldots, b_{1,r})$ is generated by $r$ elements over $F$ and each $b_{1,i}$ satisfies a quadratic equation over $F$, it folows that $[K : F] \leq 2^r$. Now $\overline{K} = \overline{F}(\overline{b}_{1,1}, \overline{b}_{1,2}, \ldots, \overline{b}_{1,r})$ and $\overline{b}_{1,i}^2 = \overline{b}_i \in S$. Since $\{\overline{b}_1, \overline{b}_2, \ldots, \overline{b}_r\}$ is a 2-independent set in $\overline{F}$, we get $[\overline{K} : \overline{F}] = 2^r$. Moreover, by [18, Theorem 11, chapter 2], $[\mathcal{T}(K) : \mathcal{T}(F)][\overline{K} : \overline{F}] = [K : F]$. Therefore, $[K : F] = 2^r$ and $[\mathcal{T}(K) : \mathcal{T}(F)] = 1$. So the lemma is proved. □

**Definition 6.4.** *Let $F$ be a field of characteristic 2, and $F^{2^{-1}}$ the extension of $F$ which is obtained by adjoining the 2-th roots of all elements of $F$. Then, $\cup_{n \in \mathbb{N}_0} F^{2^{-n}}$ is the minimal perfect extension of $F$. We denote it by $F^{2^{-\infty}}$.*

**Lemma 6.5.** *Let $F$, $\overline{F}$, $S$ and $S'$ be as defined above. Suppose that $\overline{F}$ is imperfect. Then, there exists a complete extension $L/F$ such that the residue field of $L$ is $\overline{L} = \overline{F}\left(S^{2^{-\infty}}\right)$. Moreover, if $S$ is a 2-basis of $\overline{F}$, then $\overline{L} = \overline{F}^{2^{-\infty}}$, the minimal perfect extension of $\overline{F}$.*

*Proof.* Let $T$ be a finite subset of $S'$. Then, we can form the field $F\left(T^{\frac{1}{2}}\right)$ as in Lemma 6.3. Consider $K_{1,1} = \bigcup_T F\left(T^{\frac{1}{2}}\right) = F\left(S'^{\frac{1}{2}}\right)$, where the union is over all finite subsets $T \subseteq S'$. Since $K_{1,1}$ is an algebraic extension of the complete field $F$, we extend the valuation $v$ uniquely on $K_{1,1}$, say $v_{K_{1,1}}$. Now by Lemma 6.3, $\mathcal{T}\left(F\left(T^{\frac{1}{2}}\right)\right) = \mathcal{T}(F)$, for all finite subsets $T \subseteq S'$. So $\mathcal{T}(K_{1,1}) = \mathcal{T}(F)$. Let $K_1$ denote the completion of $K_{1,1}$ with respect to the valuation $v_{K_{1,1}}$. Then, the valuation $v_{K_{1,1}}$ can be extended on $K_1$ uniquely say, $v_{K_1}$ and $\mathcal{T}(K_1) = \mathcal{T}(K_{1,1}) = \mathcal{T}(F)$. Moreover, by Lemma 6.3, the residue field of $K_1$ is given by $\overline{K}_1 = \overline{F}\left(S^{\frac{1}{2}}\right)$.

Note that $S^{\frac{1}{2}}$ is a 2-independent subset of $\overline{K}_1 = \overline{F}\left(S^{\frac{1}{2}}\right)$. Take $M_1 = S'^{\frac{1}{2}}$ and we repeat the same process. Let $K_{2,1} = K_1\left(M_1^{\frac{1}{2}}\right)$ and the completion of $K_{2,1}$ is $K_2$. Again the residue field of $K_2$ is $\overline{K}_2 = \overline{K}_1\left(S^{\frac{1}{2^2}}\right) = \overline{F}\left(S^{\frac{1}{2^2}}\right)$. By repeating this process, we get a chain of fields $F = K_0 \subseteq K_1 \subseteq K_2 \subseteq \ldots$ satisfying the two conditions:

(1) The value groups $\mathcal{T}(K_n) = \mathcal{T}(F)$ for all $n \in \mathbb{N}$.

(2) The residue field of $K_n$ is $\overline{K}_n = \overline{F}\left(S^{2^{\frac{1}{n}}}\right)$ for all $n \in \mathbb{N}_0$.

Let us consider the field $L_0 = \bigcup_{n \in \mathbb{N}} K_n$. Then, its residue field is $\overline{L}_0 = \bigcup_{n \in \mathbb{N}} \overline{K}_n = \overline{F}\left(S^{2^{-\infty}}\right)$ and $\mathcal{T}(L_0) = \mathcal{T}(F)$. Finally, let $L$ be the completion of $L_0$, with the residue field $\overline{L}$. Then $\overline{L} = \overline{F}\left(S^{2^{-\infty}}\right)$ and $\mathcal{T}(L) = \mathcal{T}(F)$. Moreover, if $S$ is a 2-basis of $\overline{F}$, then $\overline{K}_1 = \overline{F}\left(S^{\frac{1}{2}}\right) = \left(\overline{F}\right)^{\frac{1}{2}}$. Moreover, $S^{\frac{1}{2}}$ is a 2-basis of $\overline{K}_1$. Therefore, $\overline{K}_2 = \overline{F}\left(S^{\frac{1}{2^2}}\right) = \overline{K}_1\left(S^{\frac{1}{2^2}}\right) = \left(\overline{K}_1\right)^{\frac{1}{2}} = \left(\overline{F}\right)^{\frac{1}{2^2}}$. Similarly $\overline{K}_n = \left(\overline{F}\right)^{\frac{1}{2^n}}$ for all $n \in \mathbb{N}$. So $\overline{L} = \bigcup_n \overline{K}_n = \left(\overline{F}\right)^{2^{-\infty}}$, which is the minimal perfect extension of $\overline{F}$. Thus, the lemma is proved. □

Using Lemmas 6.3 and 6.5, we prove the next theorem which extends Theorem 6.1 to the case of imperfect residue field.

**Theorem 6.6.** *Suppose that the residue field $\overline{F}$ of $F$ is imperfect. Let $S$ be a 2-basis of $\overline{F}$ and $S' \subseteq A_F$, made of exactly one representative for each element of $S$. Then, there exists a unique field monomorphism $\alpha : \overline{F} \to F$ such that $\mathrm{Im}(\alpha) \subseteq A_F$, $\alpha(S) = S'$ and $\overline{\alpha(a)} = a$ for all $a \in \overline{F}$.*

*Proof.* **Existence:** By Lemma 6.5, there exists a complete field extension $L$ of $F$ whose residue field is $\overline{L} = \left(\overline{F}\right)^{2^{-\infty}}$, the minimal perfect extension of $\overline{F}$, and $L$ contains $F\left((S')^{2^{-\infty}}\right)$. Since $\overline{L}$ is perfect, it follows from Theorem 6.1 that there exists a unique field monomorphsim $\beta : \overline{L} \to L$ such that $\overline{\beta(x)} = x$ for all $x \in \overline{L}$.

**Step 1:** Let $a \in S$ and $x \in S'$ be the representative of $a$. Now for each $n \in \mathbb{N}$, $a^{2^{-n}} \in \overline{L}$ and $x^{2^{-n}} \in L$ with $\overline{x^{2^{-n}}} = a^{2^{-n}}$. Therefore
$$a = \overline{x^{2^{-n}}}^{2^n} = \overline{a_n}^{2^n}, \quad \text{where} \quad a_n = x^{2^{-n}} \; \forall n \in \mathbb{N}.$$
Then, we have
$$\beta(a) = \lim_{n \to \infty} a_n^{2^n} = x, \quad \text{(by definition of } \beta \text{ from (30))}.$$
It implies that $\beta(S) = S'$.

**Step 2:** Now we have to prove that $\beta(\overline{F}) \subseteq F$. Let $a \in \overline{F}^*$. Then $a^{2^{-n}} \in \overline{F}^{2^{-n}} = \overline{K}_n$ for all $n \in \mathbb{N}$, where $\overline{K}_n$ was defined in Lemma 6.5. Therefore, there exist $a_n \in K_n$ such that
$$v_{K_n}(a_n) = 0 \quad \text{and} \quad \overline{a}_n = a^{2^{-n}}, \quad \forall n \in \mathbb{N}.$$
Recall that $K_n$ was constructed by adjoining $2^n$-th roots of $S'$ with $F$, and then by completion. So for each $n \in \mathbb{N}$, there exists a sequence $\{a_{k,n} \mid k \in \mathbb{N}\} \subset F$ such that
$$a_n = \lim_{k \to \infty} a_{k,n}^{2^{-n}}.$$
Since $F$ is complete, we have for all $n \in \mathbb{N}$
$$(a_n)^{2^n} = \lim_{k \to \infty} a_{k,n} \in F.$$
By the definition of $\beta$ from (30), $\beta(a) = \lim_{n \to \infty}(a_n)^{2^n} \in F$. Thus $\beta(\overline{F}) \subseteq F \cap A_L = A_F$, where $A_L, A_F$ are the valuation rings of $L$ and $F$ respectively. Therefore, we define the map $\alpha : \overline{F} \to F$, such that $\alpha = \beta_{\overline{F}}$. Now $\alpha$ is a field monomorphism with $\alpha(S) = S'$ and $\overline{\alpha(a)} = a$, for all $a \in \overline{F}$.

**Uniqueness:** Let $\alpha_1 : \overline{F} \to F$ be a field monomorphism, such that $\mathrm{Im}(\alpha_1) \subseteq A_F$, $\alpha_1(S) = S'$ and $\overline{\alpha_1(a)} = a$ for all $a \in \overline{F}$. Recall that $\overline{L} = \overline{F}\left(S^{2^{-\infty}}\right)$ and $F\left((S')^{2^{-\infty}}\right) \subseteq L$. Therefore, we define a homomorphism $\beta' : \overline{L} \to L$, such that :
$$\beta'_{|\overline{F}} = \alpha_1 \quad \text{and} \quad \beta'\left(s^{2^{-n}}\right) = (\alpha_1(s))^{2^{-n}}, \quad \text{for all } n \in \mathbb{N} \text{ and for all } s \in S.$$

Since $S$ is 2-independent in $\overline{F}$, $\beta'$ is a well defined field homomorphism. Clearly $\text{Im}(\beta') \subseteq A_L$ and $\overline{\beta'(x)} = x$ for all $x \in \overline{L}$. So $\beta'$ satisfies all the conditions of $\beta$. By Theorem 6.1, $\beta$ is unique and it implies that $\beta = \beta'$. Therefore, $\alpha = \alpha_1$. Thus the theorem is proved. □

The following theorem is a generalization of Theorem 6.2 to the case of imperfect residue field.

**Theorem 6.7.** *Let $F$, $\overline{F}$, $S$, $S'$ and $\alpha$ be as defined in Theorem 6.6. Then, $F = \alpha(\overline{F})((\pi))$, where $\pi$ is an uniformizer.*

*Proof.* The same proof as that of Theorem 6.2. □

**Definition 6.8.** *The monomorphism $\alpha$ of Theorem 6.6 is called the Teichmüller lifting.*

**Corollary 6.9.** *The monomorphism $\alpha$ of Theorem 6.6 induces a natural homomorphism $\alpha' : H_2^{m+1}(\overline{F}) \to H_2^{m+1}(F)$. Moreover, $\alpha'$ is a monomorphism and for any $x \in A_F$ and $x_i \in A_F \setminus \pi A_F$, we have:*

$$\alpha'\left(\overline{x}\frac{d\overline{x}_1}{\overline{x}_1} \wedge \ldots \wedge \frac{d\overline{x}_m}{\overline{x}_m}\right) = \overline{x\frac{dx_1}{x_1} \wedge \ldots \wedge \frac{dx_m}{x_m}}.$$

*Proof.* Let $\alpha' : H_2^{m+1}(\overline{F}) \to H_2^{m+1}(F)$ be the homomorphism induced by the monomorphism $\alpha : \overline{F} \longrightarrow F$. By Theorem 6.6, for each $x \in A_F$ we have $b_x \in A_F$ such that $v(b_x) > 0$ and

$$\alpha(\overline{x}) = x + b_x.$$

Moreover, for any $x \in A_F \setminus \pi A_F$, we have

$$(\alpha(\overline{x}))^{-1} = (x + b_x)^{-1} = (y_x + c_x),$$

where $\overline{xy}_x = 1 \in \overline{F}$ and $v(c_x) > 0$. Therefore, for each $x \in A_F \setminus \pi A_F$ we get

(33) $$\frac{x}{\alpha(\overline{x})} = x(x + b_x)^{-1} = x(y_x + c_x) = 1 + d_x,$$

where $v(d_x) > 0$. Now for $x \in A_F$ and $x_i \in A_F \setminus \pi A_F$, we have:

$$\begin{aligned}
\alpha'\left(\overline{x}\frac{d\overline{x}_1}{\overline{x}_1} \wedge \ldots \wedge \frac{d\overline{x}_m}{\overline{x}_m}\right) &= \overline{\alpha(\overline{x})\frac{d(\alpha(\overline{x}_1))}{\alpha(\overline{x}_1)} \wedge \ldots \wedge \frac{d(\alpha(\overline{x}_m))}{\alpha(\overline{x}_m)}} \\
&= \overline{(x + b_x)\frac{d(\alpha(\overline{x}_1))}{\alpha(\overline{x}_1)} \wedge \ldots \wedge \frac{d(\alpha(\overline{x}_m))}{\alpha(\overline{x}_m)}} \\
&= \overline{x\frac{d(x_1 + b_{x_1})}{x_1 + b_{x_1}} \wedge \ldots \wedge \frac{d(\alpha(\overline{x}_m))}{\alpha(\overline{x}_m)}} \quad \text{(as } b_x \in \wp(F)\text{)} \\
&= \overline{\frac{xx_1}{\alpha(\overline{x}_1)}\frac{dx_1}{x_1} \wedge \ldots \wedge \frac{d(\alpha(\overline{x}_m))}{\alpha(\overline{x}_m)}} + \overline{\frac{xb_{x_1}}{\alpha(\overline{x}_1)}\frac{d(b_{x_1})}{b_{x_1}} \wedge \ldots \wedge \frac{d(\alpha(\overline{x}_m))}{\alpha(\overline{x}_m)}} \\
&= \overline{\frac{xx_1}{\alpha(\overline{x}_1)}\frac{dx_1}{x_1} \wedge \ldots \wedge \frac{d(\alpha(\overline{x}_m))}{\alpha(\overline{x}_m)}} \quad \left(\text{as } \frac{xb_{x_1}}{\alpha(\overline{x}_1)} \in \wp(F)\right) \\
&= \overline{x(1 + d_{x_1})\frac{dx_1}{x_1} \wedge \ldots \wedge \frac{d(\alpha(\overline{x}_m))}{\alpha(\overline{x}_m)}} \quad \text{(by (33))} \\
&= \overline{x\frac{dx_1}{x_1} \wedge \ldots \wedge \frac{d(\alpha(\overline{x}_m))}{(\alpha(\overline{x}_m))}} \quad \text{(since } xd_{x_1} \in \wp(F)\text{).}
\end{aligned}$$

Repeating the process for $x_2, \ldots, x_m$, we get: $\alpha'\left(\overline{x}\dfrac{\overline{\mathrm{d}\overline{x}_1}}{\overline{x}_1} \wedge \ldots \wedge \dfrac{\mathrm{d}\overline{x}_m}{\overline{x}_m}\right) = \overline{x\dfrac{\mathrm{d}x_1}{x_1} \wedge \ldots \wedge \dfrac{\mathrm{d}x_m}{x_m}}$.

Now we prove that $\alpha'$ is injective. Let $x_i \in A_F$ and $x_{i_j} \in A_F \setminus \pi A_F$ be such that $\sum_i \overline{x}_i \dfrac{\overline{\mathrm{d}\overline{x}_{i_1}}}{\overline{x}_{i_1}} \wedge \ldots \wedge \dfrac{\mathrm{d}\overline{x}_{i_m}}{\overline{x}_{i_m}} \in \ker(\alpha')$. Then, we have:

$$\alpha'\left(\sum_i \overline{x}_i \dfrac{\overline{\mathrm{d}\overline{x}_{i_1}}}{\overline{x}_{i_1}} \wedge \ldots \wedge \dfrac{\mathrm{d}\overline{x}_{i_m}}{\overline{x}_{i_m}}\right) = \sum_i \overline{x_i \dfrac{\mathrm{d}x_{i_1}}{x_{i_1}} \wedge \ldots \wedge \dfrac{\mathrm{d}x_{i_m}}{x_{i_m}}} = 0.$$

Applying the residue map $\chi_\pi$ from Definition 6.12, we get $\chi_\pi\left(\sum_i \overline{x_i \dfrac{\mathrm{d}x_{i_1}}{x_{i_1}} \wedge \ldots \wedge \dfrac{\mathrm{d}x_{i_m}}{x_{i_m}}}\right) = 0$. Using the definition of $\chi_\pi$, one gets that $\chi_\pi\left(\sum_i \overline{x_i \dfrac{\mathrm{d}x_{i_1}}{x_{i_1}} \wedge \ldots \wedge \dfrac{\mathrm{d}x_{i_m}}{x_{i_m}}}\right) = \sum_i \overline{x}_i \dfrac{\overline{\mathrm{d}\overline{x}_{i_1}}}{\overline{x}_{i_1}} \wedge \ldots \wedge \dfrac{\mathrm{d}\overline{x}_{i_m}}{\overline{x}_{i_m}} = 0 \in H_2^{m+1}(\overline{F})$. Thus $\alpha'$ is injective. Consequently, we consider $H_2^{m+1}(\overline{F})$ as a subgroup of $H_2^{m+1}(F)$. $\square$

**Corollary 6.10.** *By Theorem 6.6, $\alpha$ depends on $S, S'$. However, Corollary 6.9 implies that $\alpha'$ is independent of the choices of $S$ or $S'$.*

**Corollary 6.11.** *From Corollary 6.9, we get $H_2^{m+1}(\overline{F}) \simeq \alpha'(H_2^{m+1}(\overline{F})) \subseteq H_2^{m+1}(F)$.*

Now we define the residue maps $\zeta$ and $\chi_\pi$ from [7, Section 6], which are used several times.

**Definition 6.12.** *Let $F$ be a field of characteristic 2 with a discrete valuation $v : F \to \mathbb{Z}$. Let $\pi$ be an uniformizer, $A_F$ the valuation ring, $A_F^\times$ its group of units, and $\overline{F}$ the residue field of $F$. Let $H_2^{m+1}(F)'$ be the subgroup of $H_2^{m+1}(F)$ generated by the elements $a\dfrac{\mathrm{d}b_1}{b_1} \wedge \ldots \wedge \dfrac{\mathrm{d}b_m}{b_m}$, where $a \in A_F$ and $b_i \in F^*$. Similarly, $\overline{I}_q^{m+1}(F)' \subseteq \overline{I}_q^{m+1}(F)$ is generated by the $(m+1)$-fold Pfister quadratic forms $\overline{\langle\langle b_1 \ldots, b_m; a]]}$, where $a \in A_F$ and $b_i \in F^*$.*

*If $x \in F^*$, then we have $\langle x \rangle_b = \langle \pi^i y \rangle_b$ for some $y \in A_F^\times$ and $i = 0$ or $1$. Now we have two homomorphisms $\partial_\pi^1, \partial_\pi^2 : W(F) \to W(\overline{F})$, defined on the generators as follows:*

$$\partial_\pi^k(\langle \pi^i y \rangle_b) = \begin{cases} \langle \overline{y} \rangle_b & \text{if } k \not\equiv i \pmod{2} \\ 0 & \text{if } k \equiv i \pmod{2}. \end{cases}$$

*See [16, page 85]. The maps $\partial_\pi^1$ and $\partial_\pi^2$ are said to be the first and the second residue homomorphism associated with the valuation $v$.*

By [7, Lemma 6.4], there exist two group homomorphisms $\delta : \overline{I}_q^{m+1}(F)' \to \overline{I}_q^m(\overline{F})$ and $\Delta_\pi : \overline{I}_q^{m+1}(F)' \to \overline{I}_q^{m+1}(\overline{F})$, defined on generators as follows:

$$\delta(\overline{\langle\langle b_1 \ldots, b_m; a]]}) = \overline{\partial_\pi^1(\langle\langle b_1, \ldots, b_m \rangle\rangle_b) \otimes [1, \overline{a}]}, \quad \text{for } a \in A_F, b_i \in F^*,$$
$$\Delta_\pi(\overline{\langle\langle b_1 \ldots, b_m; a]]}) = \overline{(\partial_\pi^1 + \partial_\pi^2)(\langle\langle b_1, \ldots, b_m \rangle\rangle_b) \otimes [1, \overline{a}]}, \quad \text{where } a \in A_F, b_i \in F^*.$$

We have an isomorphism $g_{m+1} : H_2^{m+1}(F)' \to \overline{I}_q^{m+1}(F)'$, induced from Kato's isomorphism $f_{m+1}$. As an immediate consequence, we obtain the following residue maps

$$\zeta : H_2^{m+1}(F)' \to H_2^m(\overline{F}) \quad \text{and} \quad \chi_\pi : H_2^{m+1}(F)' \to H_2^{m+1}(\overline{F})$$

*such that the following diagrams commute :*

$$\begin{array}{ccc}
H_2^{m+1}(F)' \xrightarrow{g_{m+1}} \overline{I}_q^{m+1}(F)' & \quad & H_2^{m+1}(F)' \xrightarrow{g_{m+1}} \overline{I}_q^{m+1}(F)' \\
\downarrow \zeta \quad\quad \downarrow \delta & \quad & \downarrow \chi_\pi \quad\quad \downarrow \Delta_\pi \\
H_2^m(\overline{F}) \xrightarrow{f_m} \overline{I}_q^m(\overline{F}) & \quad & H_2^{m+1}(\overline{F}) \xrightarrow{f_{m+1}} \overline{I}_q^{m+1}(\overline{F})
\end{array}$$

i.e, $\zeta = f_m^{-1} \circ \delta \circ g_{m+1}$ and $\chi_\pi = f_{m+1}^{-1} \circ \Delta_\pi \circ g_{m+1}$.

Univ. Artois, UR 2462, Laboratoire de Mathématiques de Lens (LML), F-62300 Lens, France
*Email address*: ahmed.laghribi@univ-artois.fr, trisha.maiti@univ-artois.fr